\documentclass[screen]{aomart}

\usepackage[utf8]{inputenc}
\usepackage[T1]{fontenc}
\usepackage{lmodern}
\usepackage[protrusion=true,expansion=false]{microtype}

\usepackage{amsmath,amssymb,amsthm}

\usepackage{tabularray}

\usepackage{enumitem}
\usepackage{aliascnt}

\usepackage[capitalize,nameinlink]{cleveref}

\makeatletter
\newcommand{\withouttoc}[1]{%
  \begingroup
  \let\addcontentsline\@gobblethree
  #1%
  \endgroup
}
\makeatother

\fancypagestyle{firstpage}{%
  \fancyhf{}%
  \cfoot{\footnotesize\thepage}%
}
\doinumber{}

\allowdisplaybreaks[2]

\newtheorem{theorem}{Theorem}[section]

\newaliascnt{proposition}{theorem}
\newtheorem{proposition}[proposition]{Proposition}
\aliascntresetthe{proposition}

\newaliascnt{lemma}{theorem}
\newtheorem{lemma}[lemma]{Lemma}
\aliascntresetthe{lemma}

\newaliascnt{corollary}{theorem}
\newtheorem{corollary}[corollary]{Corollary}
\aliascntresetthe{corollary}

\newaliascnt{conjecture}{theorem}

\aliascntresetthe{conjecture}

\theoremstyle{definition}

\newaliascnt{definition}{theorem}
\newtheorem{definition}[definition]{Definition}
\aliascntresetthe{definition}

\newaliascnt{example}{theorem}
\newtheorem{example}[example]{Example}
\aliascntresetthe{example}

\newaliascnt{construction}{theorem}

\aliascntresetthe{construction}

\newaliascnt{assumption}{theorem}

\aliascntresetthe{assumption}

\theoremstyle{remark}

\newaliascnt{remark}{theorem}
\newtheorem{remark}[remark]{Remark}
\aliascntresetthe{remark}

\crefname{theorem}{Theorem}{Theorems}
\crefname{proposition}{Proposition}{Propositions}
\crefname{lemma}{Lemma}{Lemmas}
\crefname{corollary}{Corollary}{Corollaries}
\crefname{conjecture}{Conjecture}{Conjectures}
\crefname{definition}{Definition}{Definitions}
\crefname{example}{Example}{Examples}
\crefname{construction}{Construction}{Constructions}
\crefname{assumption}{Assumption}{Assumptions}
\crefname{remark}{Remark}{Remarks}

\numberwithin{equation}{section}

\newcommand{\cO}{\mathcal{O}}
\newcommand{\cI}{\mathcal{I}}
\newcommand{\cF}{\mathcal{F}}

\newcommand{\cH}{\mathcal{H}}

\newcommand{\XhatY}{\widehat{X}_Y}
\AtBeginDocument{}

\newcommand{\Coh}{\operatorname{Coh}}
\newcommand{\Perf}{\operatorname{Perf}}
\newcommand{\Db}{\mathrm{D}^{b}}

\newcommand{\MF}{\operatorname{MF}}

\DeclareMathOperator{\Supp}{Supp}
\DeclareMathOperator{\Sing}{Sing}

\DeclareMathOperator{\CritVal}{CritVal}
\DeclareMathOperator{\CritLoc}{CritLoc}
\DeclareMathOperator{\Tot}{Tot}
\DeclareMathOperator{\Kar}{Kar}

\newcommand{\Dtwo}{\mathrm{D}^{\mathrm{coh}}_{\mathbf Z/2}}
\newcommand{\Cliff}{\operatorname{Cl}}
\newcommand{\MorseCat}[1]{\mathsf{Q}_{#1}}
\newcommand{\MBcat}[1]{\mathsf{MB}(#1)}
\newcommand{\Aper}[2]{\mathsf{A}^{\mathrm{per}}_{#2}(#1)}
\newcommand{\QAper}[2]{\mathsf{QA}^{\mathrm{per}}_{#2}(#1)}
\newcommand{\QAtwist}[3]{\mathsf{QA}^{\mathrm{per}}_{#2}(#1;[#3])}
\newcommand{\Dper}[3]{\mathsf{D}^{\mathrm{per}}_{#2}(#1;[#3])}
\newcommand{\FormalCat}[1]{\mathsf{F}(#1)}

\title[Formal Completions and Matrix Factorization Decompositions]{Formal completions and decompositions of matrix factorization categories of Fano mirrors}
\author[Victor Ibrahim Santos El Adji]{Victor Ibrahim Santos El Adji}
\copyrightnote{}
\address{
Departamento de Matem\'atica, PUC-Rio,
Rua Marqu\^es de S\~ao Vicente, 225, Rio de Janeiro, RJ 22451-900, Brazil
}
\email{victoribrahimsantos1997@gmail.com}
\begin{document}
\begin{abstract}
We study the derived categories of matrix factorizations for certain toric
Landau--Ginzburg models that are mirror dual to del Pezzo surfaces and Fano
threefolds. Some of these models have non-isolated singularities: their
critical loci may contain curves and surfaces, possibly with intersecting
irreducible components. We handle these cases by studying the potential in
the formal neighborhood of each irreducible component of the critical locus.
Away from $\Sing(\CritLoc(W)_{\mathrm{red}})$, we construct global formal
coordinates and determine the local normal form of the singularity. Near the
singular or intersection strata, we analyze how the pieces of the whole
connected formal neighborhood are related by localization. This gives
explicit decompositions of the derived categories of matrix factorizations
after idempotent completion.  For the five rank-one pairs
$(B_d,V_{4d+2})$, we also compare projective compactifications.  The isolated
sectors and central-component counts agree for $d=2,3,4,5$ and differ for
$d=1$. Equivalence of the distinguished compactified contributions for the remaining matched pairs requires compatible global gluing data or compactified Homological Mirror
Symmetry.
\end{abstract}
\maketitle
\setcounter{tocdepth}{2}
\tableofcontents

\section{Introduction}
\label{sec:introduction}

Mirror symmetry, developed over the last four decades, predicts deep
correspondences between the symplectic geometry of a space and the algebraic
geometry of its mirror. In particular, Homological Mirror Symmetry, proposed by Kontsevich
\cite{kontsevich1994homological,kontsevich1998lectures}, predicts an equivalence of the respective categories. A Fano variety $F$ is a smooth projective variety with an ample
anticanonical line bundle. Its mirror partner is expected to be a so-called
\emph{Landau--Ginzburg model}, that is, a pair $(X,W)$ of a smooth algebraic variety $X$, usually non-proper, and a regular function $W:X\to\mathbf A^1$, known as the
\emph{superpotential} or \emph{Landau--Ginzburg potential}, or, in what follows, just a \emph{potential}.

Both sides of mirror symmetry come in families. Let
$\rho=\operatorname{rk}\operatorname{Pic}(F)$ and choose coordinates
$\mathbf q=(q_1,\dots,q_\rho)$ on the complexified K\"ahler moduli space of
$F$, the torus
\[
\mathcal K_F=\operatorname{Spec}
\mathbf C[q_1^{\pm1},\dots,q_\rho^{\pm1}].
\]
A complexified K\"ahler class $\boldsymbol\omega$ determines these coordinates
by $q_i=\exp(-\int_{\beta_i}\boldsymbol\omega)$, where $\beta_1,\dots,\beta_\rho$
is a basis of curve classes. On the mirror side, the coefficients $\mathbf y$
of a family $(X_{\mathbf y},W_{\mathbf y})$ are coordinates on a complex
moduli space. The basic fact we need is that the \emph{mirror map}
\[
\mathbf y=\mathfrak m(\mathbf q)
\]
identifies the K\"ahler moduli space of $F$ with the complex moduli space of
the mirror. The categories below come in families over these two spaces, and
mirror symmetry compares them at points matched by $\mathfrak m$. In practice
the mirror map is found by matching the quantum connection and periods of $F$
with the Picard--Fuchs system and periods of $W$
\cite{givental1998mirror,hori2000mirror}; here this matching serves only to
identify the two moduli spaces. The expected correspondences are
\[
\begin{array}{ccc}
\text{Fano variety }F && \text{Mirror }(X,W) \\[4pt]
\Db\Coh(F)
&\overset{\mathrm{HMS}}{\longleftrightarrow}&
\mathrm D^\pi\operatorname{FS}(X,W) \\[12pt]
\mathrm D^\pi\operatorname{Fuk}(F)
&\overset{\mathrm{HMS}}{\longleftrightarrow}&
\overline{\MF(X,W)}.
\end{array}
\]
The category $\Db\Coh(F)$ is the bounded derived category of coherent sheaves
on $F$; it depends on the choice of the complex structure on $F$, but does not depend on any choice of a symplectic structure. On the mirror side, $\operatorname{FS}(X,W)$ is the Fukaya--Seidel
category built from the vanishing cycles of the map $W$;
in turn, it depends on the choice of a symplectic structure on $(X,W)$, but does not depend on any choice of complex structure.
When $F$ is a del Pezzo surface, this type of equivalence of categories was constructed by
Auroux, Katzarkov, and Orlov in
\cite{auroux2006mirror}.

The lower row, which is the one relevant to this paper, decomposes further.
Since $F$ is Fano, the small quantum product on $H^{\mathrm{ev}}(F;\mathbf C)$
has structure constants that are polynomials in $\mathbf q$, so it can be
evaluated at every $\mathbf q\in\mathcal K_F$; write
$\mathrm{QH}_{\mathbf q}(F)$ for the resulting algebra. Quantum
multiplication by $c_1(F)=-c_1(K_F)$ gives an operator
\[
c_1(F)\star_{\mathbf q}(-):
\mathrm{QH}_{\mathbf q}(F)
\longrightarrow
\mathrm{QH}_{\mathbf q}(F).
\]
Its eigenvalues form a finite cover of $\mathcal K_F$. Locally they are
algebraic functions $\lambda_i(\mathbf q)$, equivalently functions of
$\boldsymbol\omega$. Where two of these functions meet, the eigenvalues
collide and $\mathrm{QH}_{\mathbf q}(F)$ may cease to be semisimple.

The Fukaya category splits into summands
\[
\mathrm D^\pi\operatorname{Fuk}(F)
=
\bigoplus_{\lambda}
\mathrm D^\pi\operatorname{Fuk}(F)_\lambda,
\]
where an object belongs to $\mathrm D^\pi\operatorname{Fuk}(F)_\lambda$ when
its obstruction term $m_0$ is $\lambda$ times the identity; geometrically,
$m_0$ counts holomorphic discs of Maslov index two. When $\boldsymbol\omega$
is monotone, that is, $[\boldsymbol\omega]$ is a positive multiple of
$c_1(F)$, a nonzero summand can occur only when $\lambda$ is an eigenvalue of
$c_1(F)\star_{\mathbf q}$ \cite[Lemma~3.2]{RitterSmith2017}.

On the mirror side the same value selects the fiber
\[
Y_{\mathbf y,\lambda}=W_{\mathbf y}^{-1}(\lambda).
\]
It is singular precisely when $\lambda$ is a critical value of
$W_{\mathbf y}$. The derived category of matrix factorizations satisfies
\[
\MF_\lambda(X_{\mathbf y},W_{\mathbf y})
\simeq
\Db_{\mathrm{sg}}(Y_{\mathbf y,\lambda})
\]
\cite{eisenbud1980homological,Orlov2004,Orlov_2011}, and the total category at
the parameter $\mathbf y$ is
\[
\MF(X_{\mathbf y},W_{\mathbf y})
=
\bigoplus_{\lambda\in\CritVal(W_{\mathbf y})}
\MF_\lambda(X_{\mathbf y},W_{\mathbf y}).
\]

After applying the mirror map $\mathbf y=\mathfrak m(\mathbf q)$, Homological
Mirror Symmetry predicts
\[
\begin{aligned}
\lambda_i(\mathbf q)\text{ an eigenvalue of }c_1(F)\star_{\mathbf q}
&\quad\longleftrightarrow\quad
\lambda_i(\mathbf q)\text{ a critical value of }
W_{\mathfrak m(\mathbf q)},\\
\mathrm D^\pi\operatorname{Fuk}(F)_{\lambda_i(\mathbf q)}
&\quad\simeq\quad
\overline{\MF_{\lambda_i(\mathbf q)}
(X_{\mathfrak m(\mathbf q)},W_{\mathfrak m(\mathbf q)})}.
\end{aligned}
\]
This relation between critical values and eigenvalues is established in
several mirror constructions \cite{Yuan2025CriticalValues}. For toric Fano
manifolds with a monotone symplectic form, the categorical equivalence
follows by combining the split-generation result of Cho, Hong, and Lau
\cite[Theorem~9.1]{ChoHongLau2019} with the generation criterion of Evans and
Lekili \cite{EvansLekili2019}, as observed in \cite{ChoHongLau2019}; for
toric Fano manifolds the potential $W_{\mathbf y}$ is Morse for $\mathbf y$ in a dense
Zariski-open subset of the parameter space.

The decomposition varies with $\boldsymbol\omega$. It is finest away from the
discriminant, where distinct critical values label distinct summands, and it
coarsens when the functions $\lambda_i(\boldsymbol\omega)$ collide. If the
critical points over a common critical value stay disjoint, the summand is
still a direct sum over the connected components of the critical locus in
that fiber. If components of the critical locus meet, the summand is the
category of the formal neighborhood of their union, which in general is not
a direct sum of contributions of the separate components.

The Laurent polynomials in this paper are taken at the anticanonical point,
$\mathbf q=(1,\dots,1)$ in standard coordinates, and at the corresponding
value of $\mathbf y$. For the models with non-isolated critical points this
is the most singular member of the family: at a generic parameter the
colliding eigenvalues and critical points separate, and the critical locus is
easier to study. All decompositions below concern this specialization.

The main difficulty addressed in this paper is that the critical locus of a
Landau--Ginzburg potential need not be finite. It may have
positive-dimensional components, and distinct components may intersect. A
decomposition then cannot be obtained by listing isolated critical points. A
positive-dimensional component contributes a category that depends on the
formal neighborhood of the whole component, not only on the singularity type
at each of its points; and where components intersect, their contributions
cannot be separated, so one has to work with the formal neighborhood of their
union.

Many mirrors of Fano varieties are given by Laurent polynomials on an
algebraic torus $(\mathbf C^*)^n$
\cite{givental1998mirror,hori2000mirror,auroux2007mirror}. This makes the
critical equations explicit, but it also introduces a limitation: the torus
is only an open chart of a larger mirror, and a birational change of
coordinates can remove critical points from that chart. The simplest example
already appears for $\mathbf P^1\times\mathbf P^1$. Before taking the
canonical specialization, consider the two-parameter family
\[
W_{q_1,q_2}(x,y)
=x+\frac{q_1}{x}+y+\frac{q_2}{y}.
\]
After choosing square roots, its four critical points and values are
\[
(x,y)=(\varepsilon_1\sqrt{q_1},\varepsilon_2\sqrt{q_2}),
\qquad
\lambda_{\varepsilon_1,\varepsilon_2}(\mathbf q)
=2\varepsilon_1\sqrt{q_1}+2\varepsilon_2\sqrt{q_2},
\]
where $\varepsilon_1,\varepsilon_2\in\{\pm1\}$. For generic parameters the
four values are distinct. On the canonical diagonal $q_1=q_2$, the two
values with $\varepsilon_1=-\varepsilon_2$ collide at zero; thus two Morse
points lie in the same singular fiber even though the points themselves
remain disjoint. The specialization $q_1=q_2=1$ gives
\[
W_1(x,y)=x+y+\frac1x+\frac1y,
\qquad
W_2(z,w)=\frac{(1+z)^2}{zw}+w.
\]
The birational change of coordinates
\[
z=\frac{x}{y},\quad w=x+y,
\qquad
x=\frac{zw}{1+z},\quad y=\frac{w}{1+z}
\]
identifies $W_1$ with $W_2$ where $x+y\ne0$, equivalently where $z+1\ne0$.
The first potential has four Morse critical points $(\pm1,\pm1)$. The two
diagonal points map to $(1,\pm2)$, which are the only critical points of
$W_2$. The other two points lie on $x+y=0$, so their image would have $w=0$
and is not contained in the second torus.

We fix the quadratic notation before using it. If $E$ is a complex vector
space, its tensor algebra is
\[
T(E)=\bigoplus_{m\geq0}E^{\otimes m},
\qquad E^{\otimes0}=\mathbf C,
\]
with multiplication given by concatenation of tensors.  We regard it as a
$\mathbf Z/2$-graded algebra by declaring a tensor in $E^{\otimes m}$ to have
parity $m$ modulo $2$.  If $q$ is a quadratic form on $E$, its
$\mathbf Z/2$-graded \emph{Clifford algebra} is
\[
\Cliff(E,q)
:=
T(E)\big/\langle v\otimes v-q(v)1:v\in E\rangle.
\]
Here the brackets denote the two-sided ideal generated by the displayed
relations; equivalently, if $b_q$ is the symmetric bilinear form obtained by
polarizing $q$, the quotient satisfies
$uv+vu=2b_q(u,v)1$.  Thus the Clifford algebra replaces the free tensor
algebra by imposing the quadratic relations determined by $q$.
For a quadratic vector bundle $(E,q)$ on a scheme $Z$, this construction gives
a sheaf of graded algebras $\Cliff(E,q)$ on $Z$. Its two-periodic category of
perfect modules is the model used below for Morse and Morse--Bott
singularities.

Let $\MorseCat{2}$ denote the two-periodic local category contributed by one
nondegenerate critical point on a surface; its Clifford-algebra description
is given in Definition~\ref{def:quadratic_categories}. At the anticanonical
point $q_1=q_2=1$ of $\mathbf P^1\times\mathbf P^1$, the eigenvalues of
$c_1(F)\star$ and the critical values of $W_1$ are $-4,0,4$. The
correspondence can be seen directly:
\[
\begin{array}{c|c|c}
\lambda & \text{critical points in }W_1^{-1}(\lambda)
& \overline{\MF_\lambda(W_1)} \\ \hline
-4 & (-1,-1) & \MorseCat{2} \\
0 & (1,-1),\ (-1,1) & \MorseCat{2}^{\oplus2} \\
4 & (1,1) & \MorseCat{2}.
\end{array}
\]
Toric Homological Mirror Symmetry therefore gives
\[
\overline{\MF(W_1)}
\simeq
\mathrm D^\pi\operatorname{Fuk}(\mathbf P^1\times\mathbf P^1)
\simeq
\MorseCat{2}^{\oplus4}
\]
\cite{ChoHongLau2019,EvansLekili2019}. By contrast, $W_2$ has only the
critical values $-4$ and $4$, and hence
\[
\overline{\MF(W_2)}\simeq\MorseCat{2}^{\oplus2}.
\]
Thus $W_1$ recovers the whole Fukaya category, whereas $W_2$ misses the
$\lambda=0$ piece, which consists of two independent Morse summands. This
example explains why our results refer to a specified Laurent polynomial on
a specified torus. Mutations and birational transformations between Laurent
polynomials are studied in
\cite{akhtar2012minkowski,galkin2010mutations}, but they do not automatically
identify the derived categories of matrix factorizations visible on the two
open charts.

There is also a geometric reason to study critical loci that are curves or
surfaces rather than finite sets. Katzarkov, Kontsevich, and Pantev
conjectured that the Hodge numbers of $F$ equal Hodge-type numbers
$f^{p,q}(Y,w)$ of a tame compactified Landau--Ginzburg model $(Y,w)$
\cite{KatzarkovKontsevichPantev2008}. For smooth Fano threefolds and the
compactified toric Landau--Ginzburg models of
\cite{CoatesCortiGalkinKasprzyk2016}, Cheltsov and Przyjalkowski proved
\[
h^{1,2}(F)=\sum_{P\in\mathbf C}(\rho_P-1),
\]
where $\rho_P$ is the number of irreducible components of the fiber
$w^{-1}(P)$ \cite[Main Theorem]{CheltsovPrzyjalkowski2025}. Thus
$h^{1,2}(F)>0$ forces a reducible fiber. Two components of a connected
fiber in the smooth threefold $Y$ that meet do so along a curve; this curve
lies in the singular locus of the fiber, hence in the critical locus of $w$.
A torus chart sees such a curve only if it meets it.
Section~\ref{sec:hodge_boundary_comparison} compares this with the critical
loci that actually occur in our models.

An isolated critical point can be studied in the completed local ring at that
point. For a critical curve or surface $Z$, completing at one point at a time
is not enough: the coordinates chosen at different points may fail to agree,
and the quadratic form in the directions normal to $Z$ may be globally
twisted. The appropriate object is the formal completion $\widehat X_Z$. If
$\mathcal I_Z$ is the ideal sheaf of $Z$, then $\widehat X_Z$ remembers all
infinitesimal neighborhoods defined by the powers $\mathcal I_Z^m$. It
therefore records both the singularity transverse to $Z$ and the way that
singularity varies along $Z$.

The potential $W$ restricts compatibly to every infinitesimal neighborhood of
$Z$. We denote the induced function on the formal completion by
\[
\widehat W:\widehat X_Z\longrightarrow\mathbf A^1_{\mathbf C}
\]
and call it the \emph{completed} or \emph{formal potential}. Concretely, if
$X=\operatorname{Spec}A$, $Z=V(I)$, and $W$ is represented by $f\in A$, then
\[
\widehat X_Z=\operatorname{Spf}\widehat A,
\qquad
\widehat A=\varprojlim_m A/I^m,
\qquad
\widehat W=(f\bmod I^m)_m\in\widehat A.
\]
Thus $\widehat W$ is the image, or formal Taylor expansion, of $W$ in the
completed structure sheaf. The formal fiber $\widehat W^{-1}(w_0)$ is the
completion of $Y=W^{-1}(w_0)$ along
$Z\cap Y$; when $Z\subset Y$, this formal fiber is simply $\widehat Y_Z$.

This leads to the basic question of the paper. Given an explicit potential
$W$ and a closed subset $Z$ of the critical locus of a fiber, can the
contribution supported on $Z$ be computed from the restriction of $W$ to
$\widehat X_Z$? Results of Orlov, Chen, Preygel, and others provide the
categorical comparison needed for this passage
\cite{ORLOV2011206,Chen2010Orlov,Preygel2011,WuXuZhou2025}. We use them in
the following form.

\begin{theorem}[Formal comparison; Proposition~\ref{prop:supported_formal_excision}
and Theorem~\ref{thm:general_formal}]
Let $X$ be a smooth quasi-projective variety, let
$W:X\to\mathbf A^1_{\mathbf C}$ be flat,
fix $w_0\in\mathbf C$, and put $Y=W^{-1}(w_0)$. Let $Z\hookrightarrow X$ be a
closed embedding and let $\widehat X_Z$ be the formal completion of $X$ along
$Z$. Completion always induces an equivalence
\[
\MF^Z_{w_0}(X,W)
\simeq
\MF^t_{w_0}(\widehat X_Z,\widehat W)
\]
with the torsion category on the formal completion. Suppose in addition that
the formal fiber is affine,
\[
\widehat W^{-1}(w_0)=\operatorname{Spf}B,
\]
where $B$ is its complete coordinate ring and $J$ is the ideal cutting out
$Z\cap Y$. Under the finiteness hypotheses of
Theorem~\ref{thm:general_formal}, if
\[
\operatorname{Sing}(\operatorname{Spec}B)\subseteq V(J),
\]
then
\[
\overline{\MF^Z_{w_0}(X,W)}
\simeq
\overline{\MF_{w_0}(\widehat X_Z,\widehat W)}.
\]
If $\operatorname{Sing}(Y)\subseteq Z\cap Y$, then the supported category on
the left may be replaced by the global category $\MF_{w_0}(X,W)$.
\end{theorem}

The first equivalence says that no information supported on $Z$ is lost by
passing to its formal neighborhood. The word \emph{torsion} means that the
objects are killed by some power of the ideal of $Z$, so they live on a
finite infinitesimal neighborhood of $Z$. The additional condition on
$\operatorname{Spec}B$ says that the completed fiber has no singular point
away from the closed subset defined by $J$. Under this condition, the torsion
category already generates the full category of the completed fiber.

When two critical components are disjoint and closed, their supported
categories form a direct sum. When they intersect, their punctured branches
are not closed in the original fiber, and such a direct sum is generally
false. We then complete along the whole union and use localization to compare
the crossing with the complement of the crossing.

A second useful principle is invariance under formal changes of coordinates.
Two completed potentials are \emph{formally contact-equivalent} if a formal
coordinate change carries one to a unit multiple of the other. Multiplication
by a unit does not change the zero fiber. The classical Mather--Yau theorem
motivates this notion for an isolated hypersurface singularity
\cite{MatherYau1982}; along a curve or surface, however, the coordinate change
must be defined over the whole component, not separately at each point.

\begin{theorem}[Contact invariance; Corollary~\ref{mather-yau}]
For $i=1,2$, let $X_i$ be a smooth quasi-projective variety over $\mathbf C$,
let $W_i:X_i\to\mathbf A^1_{\mathbf C}$ be flat, and put
$Y_i=W_i^{-1}(w_i)$. Suppose that each reduced singular locus
$Z_i=\operatorname{Sing}(Y_i)_{\mathrm{red}}$ is a finite union of smooth
closed subvarieties that are pairwise disjoint as closed subschemes. If the
components are matched by a bijection and the corresponding completed
potentials are formally contact-equivalent, then
\[
\overline{\MF_{w_1}(X_1,W_1)}
\simeq
\overline{\MF_{w_2}(X_2,W_2)}.
\]
\end{theorem}

We apply these principles to specified Laurent-polynomial models for all del
Pezzo surfaces and for the $105$ deformation families of smooth Fano
threefolds. The threefold polynomials are the models listed by
Coates--Corti--Galkin--Kasprzyk
\cite{CoatesCortiGalkinKasprzyk2016}, except for four toric families, for
which we use the Hori--Vafa polynomial of the fan
(Section~\ref{sec:fano_threefold_decompositions}). For each polynomial we determine the
critical scheme and its critical values. Along every smooth
positive-dimensional component for which we give a standard category, we
construct formal coordinates valid on the entire component and record the
equation in the normal directions. When the normal quadratic form is not
globally split, the contribution is the category of modules over its
Clifford algebra. When components of the critical locus meet, we take the
category of the formal neighborhood of their union and describe it by
localization.

For example, the del Pezzo potential
\[
W=\frac{(1+x)^2(1+y)^2}{xy}-4
\]
has critical locus consisting of the isolated point $(1,1)$ and the two lines
$\{x=-1\}$ and $\{y=-1\}$, which meet at $(-1,-1)$. Away from the
intersection, the singularity is quadratic in the direction normal to each
line. At the intersection the completed equation is $u^2v^2$. Since the lines
meet, their contribution is one category on the completed union, not two
curve categories plus a point category. A localization sequence relates this
category to the two punctured lines and to the category supported at their
intersection. The calculation appears for $\mathrm{Bl}_5\mathbf P^2$ in
Section~\ref{sec:dp_decompositions}.

The main computational result is the following.

\begin{theorem}[Open-torus decompositions]
\label{thm:open_torus_fano_decomposition}
For each of the $105$ deformation families of smooth Fano threefolds, let
$W_F:(\mathbf C^*)^3\to\mathbf A^1$ be the specified Laurent-polynomial model
listed in Section~\ref{sec:fano_threefold_decompositions}. Then its critical
locus and critical values on the open torus $(\mathbf C^*)^3$ are those listed
there, and
\[
\overline{\MF((\mathbf C^*)^3,W_F)}
\]
is equivalent to the category displayed for $F$ in
Table~\ref{tab:fano_threefold_decompositions}.
\end{theorem}

\begin{remark}[Scope and interpretation of Theorem~\ref{thm:open_torus_fano_decomposition}]
For a smooth critical curve or surface, the entry in the table comes from a
formal normal form defined along the whole component. It therefore remembers
any coefficient or Clifford algebra that cannot be removed by a global change
of coordinates. Disjoint closed components give independent direct summands.
If critical components meet, or if their connected union is singular, we keep
the formal neighborhood of the whole connected set. Its category is described
by localization and is not separated artificially into contributions from
the individual branches or strata.

The theorem concerns exactly the regular function
$W_F:(\mathbf C^*)^3\to\mathbf A^1$. It does not include critical points that
may appear on the boundary of a compactification. It also does not claim that
two mutation-equivalent Laurent polynomials have equivalent categories: such
an equivalence requires a separate argument. The direct sums in the table
come from critical subsets that are disjoint inside the torus; they are not
semiorthogonal decompositions transported from the Fano threefold.
\end{remark}

\withouttoc{\subsection{Organization of the paper}}
The paper is organized as follows. Section~\ref{sec:preliminaries} reviews
idempotent completion, singularity categories, formal completions, and matrix
factorizations. Section~\ref{sec:support_mf} develops the localization
formalism for categories supported on critical loci, and
Section~\ref{sec:mather_yau} proves the contact-invariance results used in
the computations. Section~\ref{sec:dp_decompositions} treats the del Pezzo
models, while Section~\ref{sec:fano_threefold_decompositions} gives the
open-torus decompositions for the Fano threefold models and the comparison
with their compactifications.

All computations in this paper (critical loci, critical values, Hessians,
global normal coordinates, and the Hodge and Milnor census) are available in
the public repository \cite{FanoMirrorsComputations}, together with the
Python and Singular scripts that reproduce them and the companion appendix
with the model-by-model calculations.

\withouttoc{\subsection*{Acknowledgements}}
The author would like to thank his advisor, Sergey Galkin, and Alex Junior
Gomez Saltachin for valuable discussions and guidance throughout the
development of this work. The author also gratefully acknowledges the
support of CNPq and PUC-Rio.

\section{Preliminaries}
\label{sec:preliminaries}

\subsection{Idempotent completion}\label{sec:idempotent_completion}
Idempotent completion is the categorical operation of adding the direct summands determined by projectors which need not already be represented by objects of the category. We shall use it throughout the paper, especially when comparing supported derived categories of matrix factorizations with categories on formal completions. For the triangulated and Karoubi-envelope background, see \cite{happel1988triangulated,balmer2001idempotent,maclane1998categories}.

\begin{definition}
Let $\mathcal C$ be a category. An endomorphism $e:X\to X$ is called an \emph{idempotent} if $e^2=e$. We say that $e$ \emph{splits} if there exist an object $Y$ and morphisms $i:Y\to X$ and $p:X\to Y$ such that $p\circ i=\mathrm{id}_Y$ and $i\circ p=e$.
\end{definition}

Thus splitting an idempotent means that the image of a projector is represented by an actual object of the category.

\begin{definition}
The \emph{idempotent completion}, or \emph{Karoubi envelope}, of a category $\mathcal C$ is denoted by $\overline{\mathcal C}$. Its objects are pairs $(X,e)$, where $X\in\mathcal C$ and $e:X\to X$ is an idempotent. A morphism $(X,e)\to (Y,e')$ is a morphism $f:X\to Y$ in $\mathcal C$ satisfying
\[
        f=e'\circ f\circ e.
\]
Composition is inherited from $\mathcal C$.
\end{definition}

There is a canonical fully faithful functor $\mathcal C\to \overline{\mathcal C}$ sending $X$ to $(X,\mathrm{id}_X)$. The category $\overline{\mathcal C}$ is idempotent complete, and every functor from $\mathcal C$ to an idempotent complete category factors uniquely through $\overline{\mathcal C}$ up to unique isomorphism.

\begin{definition}
A category $\mathcal C$ is called \emph{idempotent complete} if every idempotent in $\mathcal C$ splits. If $\mathcal T$ is triangulated, we write $\overline{\mathcal T}$ for its idempotent completion.
\end{definition}

For a triangulated category, passing to the idempotent completion is often forced by localization and formal-completion arguments, since local categories may contain direct summands which are not represented before taking the Karoubi envelope.

\begin{theorem}[Balmer--Schlichting, \cite{balmer2001idempotent}]
Let $\mathcal T$ be a triangulated category. Then its idempotent completion $\overline{\mathcal T}$ carries a natural triangulated structure for which the canonical functor $\mathcal T\to\overline{\mathcal T}$ is exact.
\end{theorem}

\subsection{The derived category of singularities}\label{sec:derived_singularity}
We now recall the derived category of singularities. Let $X$ be a separated Noetherian scheme of finite type over $\mathbf C$, and write $\Db(X)=\Db(\Coh(X))$. Following Orlov, we say that $X$ satisfies $(\mathrm{ELF})$ if it is separated, Noetherian, of finite Krull dimension, and every coherent sheaf on $X$ is a quotient of a locally free sheaf of finite rank; every quasi-projective scheme satisfies $(\mathrm{ELF})$.

\begin{definition}
A complex $E\in \Db(X)$ is called \emph{perfect} if it is locally quasi-isomorphic to a bounded complex of locally free sheaves of finite rank. The full triangulated subcategory of perfect complexes is denoted by $\operatorname{Perf}(X)$.
\end{definition}

Perfect complexes form the part of the derived category which is already visible on a smooth space. The singularity category is obtained by quotienting out this smooth contribution.

\begin{definition}[Orlov, \cite{Orlov2004}]
The \emph{derived category of singularities} of $X$ is the Verdier quotient
\[
        \Db_{\mathrm{sg}}(X):=\frac{\Db(X)}{\operatorname{Perf}(X)}.
\]
We also write $\Db_{\mathrm{Sing}}(X)$ for the same category.
\end{definition}

If $X$ is smooth, then every bounded coherent complex is perfect, hence $\Db_{\mathrm{sg}}(X)=0$. In general, $\Db_{\mathrm{sg}}(X)$ records the contribution of the singular locus to $\Db(X)$.

\begin{lemma}[{\cite[Lemma 1.6]{Orlov2004}}]
Any object of $\Db(X)$ which is isomorphic to a bounded complex of locally free sheaves is perfect.
\end{lemma}

The functoriality needed later follows from the preservation of perfect complexes.

\begin{lemma}[{\cite[Tag 09UA]{stacks-project}}]
Let $f:X\to Y$ be a morphism of finite Tor-dimension. If $E\in \Db(Y)$ is perfect, then $\mathbf L f^*E$ is perfect in $\Db(X)$. Hence $\mathbf L f^*$ induces a functor
\[
        \widetilde{\mathbf L f^*}:\Db_{\mathrm{sg}}(Y)\longrightarrow \Db_{\mathrm{sg}}(X).
\]
\end{lemma}

\begin{lemma}[{\cite[Tag 0B6G]{stacks-project}}]
Let $f:X\to Y$ be a proper morphism of finite Tor-dimension. (For schemes of finite type over $\mathbf C$ such a morphism is perfect, which is the hypothesis of the cited lemma.) If $E\in \Db(X)$ is perfect, then $\mathbf R f_*E$ is perfect in $\Db(Y)$. Hence $\mathbf R f_*$ induces a functor
\[
        \widetilde{\mathbf R f_*}:\Db_{\mathrm{sg}}(X)\longrightarrow \Db_{\mathrm{sg}}(Y).
\]
\end{lemma}

We shall use the following consequence of Orlov's formal-local principle. The notation $\overline{(-)}$ denotes the idempotent completion introduced in Section~\ref{sec:idempotent_completion}.

\begin{theorem}[Orlov, \cite{ORLOV2011206}]
Let $X$ and $Y$ be Noetherian schemes satisfying the hypotheses of Orlov's formal localization theorem. If their formal completions along the singular loci are isomorphic,
\[
        \widehat X_{\Sing(X)}\cong \widehat Y_{\Sing(Y)},
\]
then there is an equivalence of idempotent-complete singularity categories
\[
        \overline{\Db_{\mathrm{sg}}(X)}\simeq \overline{\Db_{\mathrm{sg}}(Y)}.
\]
\end{theorem}

This explains why idempotent completion naturally appears when global singularity categories are compared through formal neighborhoods.

\subsection{Coherent sheaves with torsion on a formal completion}\label{sec:formal_torsion}
We next introduce the formal-completion categories used in the supported setting. Let $X$ be a Noetherian scheme and let $Y\hookrightarrow X$ be a closed subscheme defined by an ideal sheaf $\cI\subset \cO_X$.

\begin{definition}[Formal completion, cf. \cite{atiyah1969introduction,EGAIII1,stacks-project}]
The \emph{formal completion} of $X$ along $Y$ is the formal scheme $\XhatY$ whose underlying topological space is $Y$ and whose structure sheaf is
\[
        \cO_{\XhatY}=\varprojlim_n \cO_X/\cI^n.
\]
We denote by $\kappa:\XhatY\to X$ the canonical morphism of locally ringed spaces.
\end{definition}

The formal completion remembers all infinitesimal neighborhoods of $Y$ inside $X$. Since the categories used later are supported on $Y$, we also need the corresponding torsion category on $\XhatY$.

\begin{definition}
Let $\Coh_Y(X)$ be the abelian category of coherent sheaves on $X$ whose support is contained in $Y$. A coherent sheaf $\cF\in\Coh(\XhatY)$ is called \emph{torsion} if, locally on $Y$, there exists $N\gg0$ such that $(\cI\cO_{\XhatY})^N\cF=0$. The full subcategory of torsion coherent sheaves on $\XhatY$ is denoted by $\Coh_t(\XhatY)$.
\end{definition}

\begin{definition}
We define the supported and torsion derived categories by
\[
        \Db_Y(X):=\{E\in\Db(X)\mid \Supp \cH^i(E)\subseteq Y \text{ for all }i\}
\]
and
\[
        \Db_t(\XhatY):=\Db\bigl(\Coh_t(\XhatY)\bigr).
\]
We also set
\[
        \operatorname{Perf}_Y(X):=\operatorname{Perf}(X)\cap \Db_Y(X).
\]
\end{definition}

\begin{theorem}[Formal torsion equivalence]
\label{thm:formal_torsion_equiv}
By \cite[Proposition~2.8]{ORLOV2011206}, the completion
functor induces an equivalence of abelian categories
\[
        \Coh_Y(X)\simeq \Coh_t(\XhatY).
\]
Since $\Db(\Coh_Y(X))\simeq\Db_Y(X)$ by \cite[Lemma~2.1]{ORLOV2011206}, it
induces an equivalence of triangulated categories
\[
        \Db_Y(X)\simeq \Db_t(\XhatY).
\]
\end{theorem}

\begin{definition}[Perfect complexes with torsion]
\label{def:formal_perfect_torsion}
Assume that $X$ satisfies $(\mathrm{ELF})$. Under the equivalence of
Theorem~\ref{thm:formal_torsion_equiv}, define
$\operatorname{Perf}_t(\XhatY)$ to be the essential image of
$\operatorname{Perf}_Y(X)$. By Orlov's intrinsic characterization
\cite[Lemma~2.4]{ORLOV2011206}, this subcategory is determined inside the
torsion derived category and does not depend on the chosen algebraic
presentation of the formal completion.
\end{definition}

This definition avoids identifying $\operatorname{Perf}_t(\XhatY)$ with the
intersection of $\Db_t(\XhatY)$ and a separately defined category of perfect
complexes on the formal scheme; such an identification is not needed here.

\begin{corollary}[Supported quotient under completion]
\label{cor:formal_supported_quotient}
Completion induces an equivalence
\[
\frac{\Db_Y(X)}{\operatorname{Perf}_Y(X)}
\simeq
\frac{\Db_t(\XhatY)}{\operatorname{Perf}_t(\XhatY)}.
\]
\end{corollary}

Thus the supported singularity quotient can be read entirely from the formal
neighborhood of $Y$. This is the formal mechanism behind the supported matrix
factorization category used in Section~\ref{sec:support_mf}.

\subsection{Matrix factorizations}\label{sec:matrix_factorizations}
We now define the derived category of matrix factorizations associated with a Landau--Ginzburg potential, following the standard two-periodic formalism of Eisenbud and Orlov and its dg-categorical refinements \cite{Dyckerhoff2011,eisenbud1980homological,Orlov2004,PolishchukVaintrob2012}. The graded, equivariant dimensional-reduction equivalence of \cite{isik2013equivalence} is closely related, but it is not invoked below: our categories are ungraded and two-periodic, and the potentials are treated directly rather than through a regular-section presentation. Let $X$ be a scheme over $\mathbf C$, let $W:X\to\mathbf A^1_{\mathbf C}$ be a regular function, and let $w_0\in\mathbf C$. We write $X_{w_0}:=W^{-1}(w_0)$.

\begin{definition}
A matrix factorization of $(X,W-w_0)$ is a diagram of locally free sheaves of finite rank
\[
        E=\left(E_1 \overset{e_1}{\underset{e_0}{\rightleftarrows}} E_0\right)
\]
such that
\[
        e_1e_0=(W-w_0)\mathrm{id}_{E_0},
        \qquad
        e_0e_1=(W-w_0)\mathrm{id}_{E_1}.
\]
The differential $d_E=e_0+e_1$ is odd with respect to the $\mathbf Z/2$-grading.
\end{definition}

\begin{definition}
The dg category $\operatorname{DG}_{w_0}(X,W)$ has matrix factorizations as objects. For two objects $E$ and $F$, the morphism complex is the $\mathbf Z/2$-graded complex
\[
        \operatorname{Hom}_{\operatorname{DG}_{w_0}(X,W)}(E,F)
        =
        \bigoplus_{i,j\in\mathbf Z/2}\operatorname{Hom}(E_i,F_j),
\]
with parity $j-i\pmod 2$ and differential
\[
        d_{E,F}(\phi)=d_F\circ\phi-(-1)^{|\phi|}\phi\circ d_E.
\]
\end{definition}

We shall use two associated categories. The category $\operatorname{Pair}_{w_0}(X,W)$ has the same objects as $\operatorname{DG}_{w_0}(X,W)$ and degree-zero closed morphisms. The homotopy category $H^0\operatorname{DG}_{w_0}(X,W)$ has the same objects and morphisms given by degree-zero cohomology of the morphism complexes.

Concretely, a morphism $p:E\to F$ in $\operatorname{Pair}_{w_0}(X,W)$ is a pair $p=(p_1,p_0)$ with $p_i:E_i\to F_i$ satisfying $p_1e_0=f_0p_0$ and $f_1p_1=p_0e_1$. Such a morphism is null-homotopic if there exist maps $s_0:E_0\to F_1$ and $s_1:E_1\to F_0$ such that
\[
        p_1=f_0s_1+s_0e_1,
        \qquad
        p_0=s_1e_0+f_1s_0.
\]

The translation functor on $H^0\operatorname{DG}_{w_0}(X,W)$ is defined by
\[
        E[1]=\left(E_0\overset{-e_0}{\underset{-e_1}{\rightleftarrows}}E_1\right).
\]
For a closed degree-zero morphism $p:E\to F$, its cone is
\[
        \operatorname{Cone}(p)=
        \left(F_1\oplus E_0\overset{c_1}{\underset{c_0}{\rightleftarrows}}F_0\oplus E_1\right),
\]
where
\[
        c_0=
        \begin{pmatrix}
        f_0 & p_1\\
        0 & -e_1
        \end{pmatrix},
        \qquad
        c_1=
        \begin{pmatrix}
        f_1 & p_0\\
        0 & -e_0
        \end{pmatrix}.
\]
The standard triangles are the triangles isomorphic to
\[
        E\xrightarrow{p}F\longrightarrow\operatorname{Cone}(p)\longrightarrow E[1].
\]
With these triangles, $H^0\operatorname{DG}_{w_0}(X,W)$ is triangulated.

\begin{definition}
Let $\operatorname{Ac}_{w_0}(X,W)$ be the smallest full triangulated subcategory of $H^0\operatorname{DG}_{w_0}(X,W)$ generated by totalizations of acyclic complexes in the exact category $\operatorname{Pair}_{w_0}(X,W)$. The triangulated category of matrix factorizations is
\[
        \MF_{w_0}(X,W):=
        \frac{H^0\operatorname{DG}_{w_0}(X,W)}
        {\operatorname{Ac}_{w_0}(X,W)}.
\]
We also write
\[
        \MF(X,W):=\bigoplus_{w_0\in\mathbf A^1}\MF_{w_0}(X,W),
\]
where only singular fibers contribute.
\end{definition}

The bridge with the previous sections is Orlov's equivalence between matrix factorizations and singularity categories.

\begin{theorem}[Orlov, \cite{Orlov2004,Orlov_2011}]
Assume that $X$ is regular and that $W-w_0$ satisfies the usual flatness hypothesis near the fiber $X_{w_0}$. Then there is an equivalence of triangulated categories
\[
        \MF_{w_0}(X,W)\simeq \Db_{\mathrm{sg}}(X_{w_0}).
\]
In particular, if $X_{w_0}$ is smooth, then $\MF_{w_0}(X,W)=0$.
\end{theorem}

Thus the derived category of matrix factorizations is controlled by the singularities of the fibers of $W$. Combined with the formal-completion results of Section~\ref{sec:formal_torsion}, this allows us to replace supported pieces of $\MF_{w_0}(X,W)$ by categories defined on formal neighborhoods.

\section{Supported Matrix Factorizations}\label{sec:support_mf}
In this section, we study matrix factorizations supported on a closed subscheme
$Z$. Formal excision works for arbitrary closed support. The stronger
assumption that $Z$ contains the singular locus of the relevant fiber is needed
only when one wants the supported category to recover the global category. Our
main objective is to isolate $\MF_{w_0}^Z(X,W)$, identify it with the torsion
category on the formal completion, and then give a separate criterion for
passing to the unrestricted formal category. This use of support, Verdier
localization, and formal completions follows the framework of Verdier,
Thomason--Trobaugh, Neeman, Orlov, and the matrix-factorization localization
formalism developed for Landau--Ginzburg models by Polishchuk--Vaintrob,
Preygel, and Hirano
\cite{neeman2001triangulated,ORLOV2011206,Orlov_2011,thomason1990higher,verdier1996categories,Preygel2011,Hirano2017}.
The supported quotient used here is the category denoted
$D^b_{\mathrm{Sing}}(Y,T)$ by Efimov--Positselski. Related localization and
torsion-formal formulations appear in
\cite{Chen2010Orlov,EfimovPositselski2015,WuXuZhou2025}.
\begin{definition}
Let $D^b_Z(X)\subset D^b(X)$ denote the full triangulated subcategory consisting of complexes whose cohomology sheaves are supported on $Z$. Equivalently, $D^b_Z(X)$ consists of complexes whose restriction to the open subset $U=X\setminus Z$ is acyclic. We also set
\[
\operatorname{Perf}_Z(X):=\operatorname{Perf}(X)\cap D^b_Z(X).
\]
\end{definition}

The supported category isolates the contribution of the singular fiber near $Z$. Thus, replacing $\MF_{w_0}(X,W)$ by $\MF^Z_{w_0}(X,W)$ is a way of focusing on the part of the Landau--Ginzburg model that is actually responsible for the singular contribution.

\begin{definition}
 Let $X$ be a regular $\mathbf{C}$-scheme and $W:X\to \mathbf{A}^1_{\mathbf{C}}$ be a flat morphism. We define the supported derived category of matrix factorizations by
    \[
    \MF^{Z}_{w_0}(X,W)=\frac{D^b_{Z\cap W^{-1}(w_0)}(W^{-1}(w_0))}{\Perf_{Z\cap W^{-1}(w_0)}(W^{-1}(w_0))}.
    \]
\end{definition}

This notation makes the support condition part of the definition of the category. The following proposition explains how this supported piece behaves under restriction to the complement of $Z$.

\begin{proposition}[Chen; Hirano]\label{propSupMF}
Let $X$ be a regular $\mathbf C$-scheme satisfying the condition
$(\mathrm{ELF})$, let $Z\subset X$ be a closed subscheme, and let
$U=X\setminus Z$ be its open complement. Let
$W:X\to\mathbf A^1_{\mathbf C}$ be a flat morphism, fix $w_0\in\mathbf C$,
and set $U_0=U\cap W^{-1}(w_0)$. Let
\[
\langle \MF^Z_{w_0}(X,W)\rangle\subset \MF_{w_0}(X,W)
\]
denote the thick triangulated subcategory generated by the essential image of the natural functor from the supported category. Then the restriction functor induces an equivalence
\[
\frac{\MF_{w_0}(X,W)}{\langle \MF^Z_{w_0}(X,W)\rangle}
\cong
\MF_{w_0}(U,W_{|U}).
\]
\end{proposition}
\begin{proof}
Let $X_0=W^{-1}(w_0)$, let $Z_0=Z\cap X_0$, and let
$j_0:U_0\hookrightarrow X_0$ be the induced open immersion. Put
\[
\mathcal A=D^b(X_0),\qquad
\mathcal B=\operatorname{Perf}(X_0),\qquad
\mathcal C=D^b_{Z_0}(X_0).
\]
Then
$\mathcal B\cap\mathcal C=\operatorname{Perf}_{Z_0}(X_0)$. Let
$\mathcal D=\langle\mathcal B,\mathcal C\rangle$ be the smallest thick
subcategory of $\mathcal A$ containing $\mathcal B$ and $\mathcal C$. If
$Q_{\mathcal B}:\mathcal A\to\mathcal A/\mathcal B$ is the quotient functor,
Verdier's quotient theorem gives
\[
\frac{\mathcal A}{\mathcal D}
\cong
\frac{\mathcal A/\mathcal B}
{\langle Q_{\mathcal B}(\mathcal C)\rangle}.
\]
By Orlov's equivalence,
$\mathcal A/\mathcal B\simeq \MF_{w_0}(X,W)$. Moreover, the essential image of
\[
\mathcal C/(\mathcal B\cap\mathcal C)=\MF^Z_{w_0}(X,W)
\]
in $\mathcal A/\mathcal B$ generates
$\langle Q_{\mathcal B}(\mathcal C)\rangle$. Thus the right-hand side is the
quotient of $\MF_{w_0}(X,W)$ by
$\langle \MF^Z_{w_0}(X,W)\rangle$.

On the other hand, coherent localization for $j_0$ gives
$\mathcal A/\mathcal C\simeq D^b(U_0)$. Applying the same quotient theorem
after first quotienting by $\mathcal C$ yields
\[
\frac{\mathcal A}{\mathcal D}
\cong
\frac{D^b(U_0)}{\langle j_0^*\operatorname{Perf}(X_0)\rangle}.
\]
This follows from the proof of the localization theorem
\cite[Theorem~1.3]{Chen2010Orlov}, which states
$\Db_{\mathrm{sg}}(X_0)/\langle\Db_{Z_0}(X_0)\rangle\simeq\Db_{\mathrm{sg}}(U_0)$;
compare the analogue for coherent factorizations in
\cite[Proposition~2.26(2)]{Hirano2017}. For completeness, the remaining
perfect-complex step follows from Thomason--Trobaugh: Key
Proposition~5.2.2(a) implies that the image of
$j_0^*: \operatorname{Perf}(X_0)\to\operatorname{Perf}(U_0)$ is dense: every
perfect complex $F$ on $U_0$ is a direct summand of an object in the image,
because $[F\oplus F[1]]=0$ lies in the image of
$K_0(X_0)\to K_0(U_0)$ \cite[Key Proposition~5.2.2(a)]{thomason1990higher}.
Consequently,
\[
\langle j_0^*\operatorname{Perf}(X_0)\rangle
=\operatorname{Perf}(U_0).
\]
Therefore, using Orlov's equivalence also on $U$, we obtain
\[
\frac{\mathcal A}{\mathcal D}
\cong
\frac{D^b(U_0)}{\operatorname{Perf}(U_0)}
\simeq \MF_{w_0}(U,W_{|U}).
\]
Combining the two descriptions proves the claim.
\end{proof}

Using this approach, we can analyze the category $\MF_{w_0}(X,W)$ by
considering closed supports that contain the singular locus of the relevant
fiber. If the complement of such a support has smooth fiber over $w_0$, then
the quotient above is zero, so the supported piece thickly generates the
corresponding derived category of matrix factorizations. This does not by itself give a
direct-sum decomposition; for that, the supported components must be genuinely
disjoint.

\begin{proposition}[Disjoint supported components]\label{prop:disjoint_supports}
Let $X$ be a regular $\mathbf C$-scheme satisfying $(\mathrm{ELF})$, let
$W:X\to\mathbf A^1_{\mathbf C}$ be a flat morphism, fix $w_0\in\mathbf C$,
and set $Y=W^{-1}(w_0)$. Let
\[
S=\coprod_{\alpha=1}^m Z_\alpha\subset Y
\]
be a finite disjoint union of closed subsets. Then
\[
\MF^S_{w_0}(X,W)
\simeq
\bigoplus_{\alpha=1}^m \MF^{Z_\alpha}_{w_0}(X,W).
\]
If $\operatorname{Sing}(Y)\subseteq S$ set-theoretically, then
\[
\overline{\MF_{w_0}(X,W)}
\simeq
\bigoplus_{\alpha=1}^m
\overline{\MF^{Z_\alpha}_{w_0}(X,W)}.
\]
\end{proposition}

\begin{proof}
Since the $Z_\alpha$ are pairwise disjoint and closed, the category of coherent
sheaves with set-theoretic support on $S$ decomposes as
\[
\operatorname{Coh}_S(Y)
\simeq
\bigoplus_{\alpha=1}^m\operatorname{Coh}_{Z_\alpha}(Y).
\]
This induces compatible decompositions of $D^b_S(Y)$ and
$\operatorname{Perf}_S(Y)$, and hence of their Verdier quotient. This proves
the first assertion.

If $\operatorname{Sing}(Y)\subseteq S$, then
\cite[Lemma~2.6]{ORLOV2011206} shows that the natural functor
\[
\MF^S_{w_0}(X,W)\longrightarrow \MF_{w_0}(X,W)
\]
is fully faithful. Moreover, \cite[Proposition~2.7]{ORLOV2011206} shows that
every object of $D^b_{\mathrm{sg}}(Y)$ is a direct summand of an object
represented by a complex supported on $\operatorname{Sing}(Y)$, hence on $S$.
The functor is therefore dense as well as fully faithful. Passing to
idempotent completions and using the first assertion gives the required
direct-sum decomposition.
\end{proof}

\begin{lemma}[Supported reduction to the singular locus]
\label{lem:supported_reduction_singular}
Let $Y$ be a separated scheme of finite type over $\mathbf C$ satisfying
$(\mathrm{ELF})$, let $T\subset Y$ be a closed subset, and put
$R=T\cap\operatorname{Sing}(Y)$. Then $R$ is closed, and the natural functor
\[
\frac{D^b_R(\operatorname{coh}Y)}{\operatorname{Perf}_R(Y)}
\longrightarrow
\frac{D^b_T(\operatorname{coh}Y)}{\operatorname{Perf}_T(Y)}
\]
is fully faithful and dense. Equivalently, every object of the target is a
direct summand of an object in the source. In particular, the idempotent
completions of these categories are equivalent.
\end{lemma}

\begin{proof}
Since $Y$ is of finite type over $\mathbf C$, its regular locus is open, so
$\operatorname{Sing}(Y)$ and $R$ are closed.
By \cite[Lemma~2.6]{ORLOV2011206}, both supported quotient categories embed
fully faithfully in $D^b_{\mathrm{sg}}(Y)$, so it remains to prove density.
Let $E\in D^b_T(\operatorname{coh}Y)$ and choose a bounded-above locally free
resolution $P^\bullet\to E$. Set $V=Y\setminus R$. The restriction $E|_V$
is perfect: it vanishes away from $T$, while every point of
$T\cap V=T\setminus\operatorname{Sing}(Y)$ is regular. Hence we may choose
$n$ so large that
\[
\operatorname{Hom}_{D^b(V)}(E|_V,H[n+1])=0
\]
for every coherent sheaf $H$ on $V$. The bound is uniform in $H$: cover $V$
by finitely many opens on which $E|_V$ is a complex of locally free sheaves in
degrees $[a,b]$; then $\operatorname{RHom}(E|_V,H)$ has cohomology sheaves in
degrees $[-b,-a]$, and Grothendieck vanishing on the Noetherian space $V$ of
finite Krull dimension gives
$\operatorname{Hom}(E|_V,H[m])=H^m(V,\operatorname{RHom}(E|_V,H))=0$ for
$m>-a+\dim V$. Enlarging $n$ further, we may also assume that $-n-1$ is below
the lowest degree in which $E$ has cohomology, so that the brutal truncation
$\sigma_{\le -n-1}P^\bullet$ is quasi-isomorphic to a single coherent sheaf
placed in degree $-n-1$. Brutal truncation at this fixed $n$ gives a morphism
\[
\alpha:E\longrightarrow G[n+1]
\]
to a coherent sheaf $G$, whose cone is perfect; hence $\alpha$ is an
isomorphism in $D^b_{\mathrm{sg}}(Y)$.

Since $E|_V$ is perfect, $\alpha|_V=0$. By the localization argument of
\cite[Lemma~2.2]{ORLOV2011206}, $\alpha$ factors through an object
$A\in D^b_R(\operatorname{coh}Y)$:
\[
E\longrightarrow A\longrightarrow G[n+1].
\]
Since $\alpha$ is invertible in $D^b_{\mathrm{sg}}(Y)$, the first morphism has
a left inverse there. Thus $E$ is a direct summand of the image of $A$.
This proves density, and the assertion about idempotent completions follows.
\end{proof}

\paragraph{Choice of support.}
For applications of the formal-completion theorem below, one should use the
reduced support of the singular locus, or a smooth reduced closed support
containing it. For a fiber with local equation $u^n=0$, with $n\geq3$, the
scheme-theoretic singular locus is defined by $u^{n-1}=0$ and is therefore
nonreduced, so it is not a smooth subvariety. Proposition~\ref{prop:disjoint_supports}
applies only to genuinely disjoint closed components. In particular, a
stratification by punctured curves and their intersection points does not
produce a direct-sum decomposition.

Next, we identify the supported category $\MF_{w_0}^Z(X,W)$ with a derived category of matrix factorizations on a formal neighborhood. This makes the categorical structure explicit and guides the subsequent development of the paper.

\begin{definition}
Let $X$ be a smooth variety, and let $Z\hookrightarrow X$ be a closed embedding. Denote the formal completion of $X$ along $Z$ by $\widehat X_Z$, as in Section~\ref{sec:formal_torsion}. We define the torsion derived category of matrix factorizations on the formal completion by
\[
\MF_{w_0}^t(\widehat X_Z,\widehat{W})
:=
\frac{D^b_t(\widehat{W}^{-1}(w_0))}{\operatorname{Perf}_t(\widehat{W}^{-1}(w_0))}.
\]
We also define the unrestricted formal category by
\[
\MF_{w_0}(\widehat X_Z,\widehat W)
:=
\frac{D^b(\operatorname{coh}\widehat W^{-1}(w_0))}
{\operatorname{Perf}(\widehat W^{-1}(w_0))}.
\]
Here a perfect complex on the formal fiber means a complex locally
quasi-isomorphic to a bounded complex of finite-rank locally free modules.
The completion of an object of $\operatorname{Perf}_{Z\cap Y}(Y)$ is perfect
in this sense, so the inclusion of torsion complexes carries
$\operatorname{Perf}_t$ into $\operatorname{Perf}$.
The inclusion of torsion coherent complexes induces a natural exact functor
\[
\iota_Z:
\MF^t_{w_0}(\widehat X_Z,\widehat W)
\longrightarrow
\MF_{w_0}(\widehat X_Z,\widehat W).
\]

\end{definition}

\begin{remark}[Meaning of the formal matrix-factorization notation]
\label{rem:formal_mf_meaning}
The two categories in the preceding definition are, by convention,
singularity quotients of the formal fiber. We do not apply Orlov's
scheme-theoretic matrix-factorization theorem directly to a formal scheme. If
the formal fiber is affine, say $\operatorname{Spf}B$ for a complete
Noetherian ring $B$, the equivalence between coherent sheaves on
$\operatorname{Spf}B$ and finite $B$-modules identifies the unrestricted
formal category with
\[
D^b_{\mathrm{sg}}(\operatorname{Spec}B)
=\frac{D^b(\operatorname{mod}B)}{\operatorname{Perf}(B)},
\]
and identifies the torsion formal category with the corresponding quotient
supported on the ideal of definition. Thus every use below of Orlov's
scheme-theoretic density results is made on $\operatorname{Spec}B$, while the
$\operatorname{Spf}B$ notation records the chosen formal neighborhood.
\end{remark}

The point of introducing the torsion version is that it packages the same
supported information intrinsically on the formal scheme. In practice, this
allows us to replace calculations on an ambient variety by calculations in a
formal neighborhood without adding singular contributions away from the
center.

The following supported form of Orlov's formal-local principle is the key
comparison result.

\begin{corollary}[Supported formal-local principle]
\label{cor:supported_formal_local}
Let $Y$ and $Y'$ be separated schemes of finite type over $\mathbf C$
satisfying $(\mathrm{ELF})$, and let $T\subset Y$ and $T'\subset Y'$ be closed
subsets. Put
\[
R=T\cap\operatorname{Sing}(Y),
\qquad
R'=T'\cap\operatorname{Sing}(Y').
\]
If the formal completions $\widehat Y_R$ and $\widehat{Y'}_{R'}$ are
isomorphic, then
\[
\overline{
\frac{D^b_T(\operatorname{coh}Y)}{\operatorname{Perf}_T(Y)}}
\cong
\overline{
\frac{D^b_{T'}(\operatorname{coh}Y')}{\operatorname{Perf}_{T'}(Y')}}.
\]
\end{corollary}

\begin{proof}
By \cite[Corollary~2.9]{ORLOV2011206}, the isomorphism of formal completions
induces an exact equivalence
\[
D^b_R(\operatorname{coh}Y)
\cong
D^b_{R'}(\operatorname{coh}Y').
\]
The intrinsic characterization of perfect complexes in
\cite[Lemma~2.4]{ORLOV2011206} shows that this equivalence identifies
$\operatorname{Perf}_R(Y)$ with $\operatorname{Perf}_{R'}(Y')$. It therefore
induces an equivalence of quotient categories supported on $R$ and $R'$.
Lemma~\ref{lem:supported_reduction_singular}, applied to $(Y,T)$ and
$(Y',T')$, identifies the idempotent completions of these quotients with those
supported on $T$ and $T'$, respectively. Composing the resulting equivalences
proves the claim.
\end{proof}

\begin{proposition}[Supported formal excision]
\label{prop:supported_formal_excision}
Let $X$ be a smooth quasi-projective variety over $\mathbf C$, let
$W:X\to\mathbf A^1_{\mathbf C}$ be flat, fix $w_0\in\mathbf C$, and put
$Y=W^{-1}(w_0)$. For any closed embedding $Z\hookrightarrow X$, let
$T=Z\cap Y$ and let $\widehat X_Z$ be the formal completion of $X$ along
$Z$. Then completion induces an equivalence
\[
\MF^Z_{w_0}(X,W)
\xrightarrow{\sim}
\MF^t_{w_0}(\widehat X_Z,\widehat W).
\]
No smoothness assumption on $Z$, affineness assumption on the formal fiber, or
condition on the singular locus of the completed fiber is required.
\end{proposition}

\begin{proof}
The formal fiber $\widehat W^{-1}(w_0)$ is the completion $\widehat Y_T$ of
$Y$ along $T$. By \cite[Lemma~2.1 and Proposition~2.8]{ORLOV2011206},
completion induces an equivalence
\[
D^b_T(\operatorname{coh}Y)
\xrightarrow{\sim}
D^b_t(\operatorname{coh}\widehat Y_T).
\]
By Definition~\ref{def:formal_perfect_torsion}, it identifies the relevant
perfect subcategories. Passing to Verdier quotients gives the stated
equivalence.
\end{proof}

\begin{lemma}[Affine formal generation; cf.\
{\cite[Lemma~2.6 and Proposition~2.7]{ORLOV2011206}}]
\label{lem:affine_formal_generation}
Let $B$ be a Noetherian ring of finite Krull dimension, complete with respect
to an ideal $J$. Assume that $B$ is excellent and that
\[
\operatorname{Sing}(\operatorname{Spec}B)\subseteq V(J).
\]
Then the natural exact functor
\[
\frac{D^b_{J\text{-}\mathrm{tor}}(B)}
{\operatorname{Perf}_{J\text{-}\mathrm{tor}}(B)}
\longrightarrow
\frac{D^b(\operatorname{mod}B)}{\operatorname{Perf}(B)}
\]
is fully faithful, and every object of the target is a direct summand of an
object in its essential image.
\end{lemma}

\begin{proof}
Because $B$ is Noetherian and $J$-adically complete, coherent sheaves on
$\operatorname{Spf}B$ correspond to finite $B$-modules
\cite{EGAIII1,stacks-project}. Under this equivalence, torsion coherent sheaves
correspond to finite modules annihilated by a power of $J$; compare
\cite[Proposition~2.8]{ORLOV2011206}. For a finite $B$-module, this annihilation
condition is equivalent to having support contained in $V(J)$. Thus the
displayed functor identifies with
\[
\frac{D^b_{V(J)}(\operatorname{coh}\operatorname{Spec}B)}
{\operatorname{Perf}_{V(J)}(\operatorname{Spec}B)}
\longrightarrow
D^b_{\mathrm{sg}}(\operatorname{Spec}B).
\]
Since $\operatorname{Spec}B$ is affine, Noetherian, and of finite Krull
dimension, it satisfies $(\mathrm{ELF})$. Thus
\cite[Lemma~2.6]{ORLOV2011206} shows that this functor is fully faithful.

For density, put $S=\operatorname{Sing}(\operatorname{Spec}B)$. By
\cite[Proposition~2.7]{ORLOV2011206}, every object of
$D^b_{\mathrm{sg}}(\operatorname{Spec}B)$ is a direct summand of an object
coming from the quotient supported on $S$. The hypothesis gives
$S\subseteq V(J)$, so this supported quotient factors through the quotient
supported on $V(J)$. Consequently, the latter quotient is dense in the full
singularity category. This proves the claim.
This is the affine adic form of the supported-generation principle; compare
also \cite[Proposition~1.2 and Corollary~2.4]{Chen2010Orlov} and
\cite[Proposition~4.1.6]{Preygel2011}.
\end{proof}

\begin{theorem}[Upgrade to the unrestricted formal category]\label{thm:general_formal}
Let $X$ be a smooth quasi-projective variety over $\mathbf C$, let
$W:X\to\mathbf A^1_{\mathbf C}$ be a flat potential, fix $w_0\in\mathbf C$,
and put $Y=W^{-1}(w_0)$. Let
$Z\hookrightarrow X$ be a closed embedding and let $\widehat X_Z$ be
the formal completion of $X$ along $Z$, with induced potential $\widehat W$.
Assume that the formal fiber is affine,
\[
\widehat W^{-1}(w_0)=\operatorname{Spf}B,
\]
where $B$ is an excellent Noetherian ring of finite Krull dimension that is
$J$-adically complete, and where
$J\subset B$ is the ideal of definition induced by the ideal of
$T=Z\cap Y$ in $Y$. Assume that
\[
\operatorname{Sing}(\operatorname{Spec}B)\subseteq V(J).
\]
Then
\[
\overline{\MF^Z_{w_0}(X,W)}
\cong
\overline{\MF_{w_0}(\widehat X_Z,\widehat W)}.
\]
If, moreover, $\operatorname{Sing}(Y)\subseteq Z\cap Y$, then
\[
\overline{\MF_{w_0}(X,W)}
\cong
\overline{\MF_{w_0}(\widehat X_Z,\widehat W)}.
\]
\end{theorem}

\begin{proof}
Let $T=Z\cap Y$. Proposition~\ref{prop:supported_formal_excision} gives
\[
\MF^Z_{w_0}(X,W)
\cong
\MF^t_{w_0}(\widehat X_Z,\widehat W).
\]
Under the identifications of Remark~\ref{rem:formal_mf_meaning}, the
perfect subcategories match those of
Lemma~\ref{lem:affine_formal_generation}. First, $\operatorname{Perf}$ of the
formal fiber corresponds to $\operatorname{Perf}(B)$: perfectness can be
tested after the flat maps from $B$ to its completed localizations, which are
jointly faithfully flat because every maximal ideal of $B$ contains $J$.
Second, $\operatorname{Perf}_t$ corresponds to
$\operatorname{Perf}_{J\text{-}\mathrm{tor}}(B)$: by
\cite[Lemma~2.4]{ORLOV2011206}, applied on both sides of the exact
equivalence $D^b_T(\operatorname{coh}Y)\simeq
D^b_{V(J)}(\operatorname{coh}\operatorname{Spec}B)$, the perfect objects are
characterized intrinsically inside these categories. Hence the functor of
Lemma~\ref{lem:affine_formal_generation}, composed with the excision
equivalence, is $\iota_Z$. By that lemma, $\iota_Z$ is fully faithful
and every object of the target is a direct summand of an object in its
essential image. Hence $\iota_Z$ becomes an equivalence after idempotent
completion. This proves the first assertion.

If $\operatorname{Sing}(Y)\subseteq T$, then
Lemma~\ref{lem:supported_reduction_singular}, applied first with support $Y$
and then with the given support $T$, identifies both idempotent completions
with the quotient supported on $\operatorname{Sing}(Y)$. Hence
\[
\overline{\MF_{w_0}(X,W)}
\cong
\overline{\MF^Z_{w_0}(X,W)}.
\]
Composing this equivalence with the first assertion proves the second.
\end{proof}

\begin{remark}[A local geometric sufficient condition]
The formal hypothesis in Theorem~\ref{thm:general_formal} is implied by the
following more geometric condition: there is an open neighborhood $U$ of
$T=Z\cap Y$ in $Y$ such that
\[
\operatorname{Sing}(U)\subseteq T.
\]
Indeed, replacing $Y$ by $U$ does not change its completion along $T$.
Schemes of finite type over $\mathbf C$ are excellent, so their formal fibers
are geometrically regular; regularity away from $T$ is therefore preserved by
completion \cite[Tags 07GG and 07QS]{stacks-project}. The stronger condition
$\operatorname{Sing}(Y)\subseteq T$ is often convenient when the global and
supported categories are to be identified simultaneously.
\end{remark}

\begin{remark}[Checking the unrestricted formal category]
\label{rem:check_formal_generation}
The supported-completion comparison always produces the torsion category
$\MF^t_{w_0}(\widehat X_Z,\widehat W)$. Replacing it by the unrestricted formal
category requires the additional hypotheses of
Lemma~\ref{lem:affine_formal_generation}, in particular
\[
\operatorname{Sing}(\operatorname{Spec}B)\subseteq V(J).
\]
Accordingly, every later application of
Theorem~\ref{thm:general_formal} must verify this inclusion for the completed
fiber in question. If it is not verified, the justified conclusion is the
torsion category rather than $\MF_{w_0}(\widehat X_Z,\widehat W)$.
\end{remark}

\begin{remark}[Why the support hypothesis is necessary]
\label{rem:idempotent_completion}
The sufficient condition
$\operatorname{Sing}(\operatorname{Spec}B)\subseteq V(J)$ cannot in general
be omitted. For example, take
\[
X=\mathbf A^2=\operatorname{Spec}\mathbf C[x,y],
\qquad W=y^2,
\qquad Z=V(x),
\qquad w_0=0.
\]
Both $X$ and $Z$ are smooth, but the completed zero fiber is represented by
\[
B=\mathbf C[y][[x]]/(y^2)
\cong\mathbf C[[x,y]]/(y^2),
\qquad J=(x).
\]
Its singular locus is $V(y)$, which is not contained in $V(x)$. The
factorization with both arrows given by multiplication by $y$ has cokernel
$B/(y)\cong\mathbf C[[x]]$. This module remains nonperfect after inverting
$x$: its restriction is $B_x/(y)\cong\mathbf C((x))$ over
$B_x\cong\mathbf C((x))[y]/(y^2)$. By contrast, every $J$-torsion object is
killed by the exact restriction functor to $\operatorname{Spec}B_x$. Thus the
class of $B/(y)$ cannot be a direct summand of an object from the torsion
category, and density fails even after idempotent completion. Notice also
that $B$ is a free rank-two $\mathbf C[[x]]$-module, hence is Noetherian,
one-dimensional, and $(x)$-adically complete; the failure is exactly the
failure of the singular-locus hypothesis.
\end{remark}

\section{Mather--Yau and contact invariance for matrix factorizations}\label{sec:mather_yau}
This section connects analytic hypersurface singularities with matrix
factorization categories. We quote the generalized analytic Mather--Yau theorem
in its filtered form. In families, however, the categorical comparison requires
formal contact equivalences that fix the base and preserve the formal normal
directions; no relative converse is inferred merely from Tjurina-type quotients.
Cotangent-complex and derived de~Rham methods provide a broader deformation
framework \cite{Illusie1971,Ton2011}, but the arguments here remain in classical
adic completions and explicit formal coordinates.

\begin{definition}
Let $f,g\in \mathbf{C}\{x_1,\dotsc,x_n\}$. We say that $f$ and $g$ are \textbf{contact equivalent}, denoted by $f\overset{c}{\sim}g$, if there exists an automorphism $\phi\in \mathrm{Aut}(\mathbf{C}\{x_1,\dotsc,x_n\})$ and a unit $u\in \mathbf{C}\{x_1,\dotsc,x_n\}^*$ such that $f=u\cdot \phi (g)$.
\end{definition}

This relation preserves the analytic hypersurface singularity defined by the
potential, up to coordinate change and multiplication by a unit; see
\cite{arnold1985singularities,greuel2007introduction}. It is also the natural
relation for Landau--Ginzburg models, since matrix factorizations are insensitive
to these changes \cite{Dyckerhoff2011,Orlov2004}.

We use the following relative version of this definition.
\begin{definition}
Let $A$ be a smooth finitely generated $\mathbf C$-algebra or a localization
of one, and let $B=A[[u_1,\dotsc,u_m]]$. We say that $f,g\in B$ are
\textbf{contact equivalent over $A$}, denoted by
$f\overset{c_A}{\sim}g$, if there exist a continuous $A$-algebra automorphism
$\phi$ of $B$ satisfying
$\phi(\mathfrak m_u)=\mathfrak m_u$ for
$\mathfrak m_u=(u_1,\dotsc,u_m)$ and a unit $U\in B^*$ such that
$f=U\phi(g)$.
\end{definition}

The relative version keeps the base algebra fixed and changes only the formal
normal variables. This prevents the parameters of a family from being mixed with
the local coordinates of the singularity.

Passing from convergent to formal series is harmless for this purpose: every
automorphism of $\mathbf C\{x_1,\dotsc,x_n\}$ induces an automorphism of
$\mathbf C[[x_1,\dotsc,x_n]]$.

\begin{theorem}[Generalized analytic Mather--Yau,
{\cite[Theorem~2.26]{greuel2007introduction}}]
\label{matheryauoriginal}
Let $f,g\in\mathfrak m\subset\mathbf C\{x_1,\dotsc,x_n\}$, and let
$j(f)$ denote the Jacobian ideal of $f$. The following are equivalent:
\begin{enumerate}[label=\textup{(\alph*)}]
\item $f\overset c\sim g$;
\item for every $b\geq0$ there is an isomorphism of local $\mathbf C$-algebras
\[
\frac{\mathbf C\{x_1,\dotsc,x_n\}}
{(f,\mathfrak m^b j(f))}
\cong
\frac{\mathbf C\{x_1,\dotsc,x_n\}}
{(g,\mathfrak m^b j(g))};
\]
\item such an isomorphism exists for some $b\geq0$.
\end{enumerate}
In particular, the case $b=0$ says that $f$ and $g$ are contact equivalent if
and only if their analytic Tjurina algebras are isomorphic.
\end{theorem}

The classical theorem of Mather and Yau treats isolated hypersurface
singularities \cite{MatherYau1982}; the filtered statement above is the precise
generalized analytic form used here. Completion of an analytic contact
equivalence gives an isomorphism of formal Tjurina algebras,
\[
\frac{\mathbf{C}[[x_1,\dotsc,x_n]]}{\langle f,\nabla f \rangle} \cong \frac{\mathbf{C}[[x_1,\dotsc,x_n]]}{\langle g,\nabla g\rangle},
\]
and the theorem supplies the converse in the analytic category.
Only the forward invariance direction is used below for filtered quotients.
The categorical corollaries for families assume relative formal contact
equivalence directly.

For families we separate the base variables from the formal normal variables.
Let $A$ be as in the relative definition, and let $B=A[[\mathbf u]]$ with
$\mathfrak m_u=\langle u_1,\dots,u_m\rangle$. The relative Tjurina-type
quotients used below replace the full Jacobian ideal by
$\langle f,\mathfrak m_u^kJ_{\mathbf u}(f)\rangle$. When $A=\mathbf C$ and
$k=0$, this is the usual Tjurina ideal. Over a nontrivial base these quotients
are useful invariants, but they are not complete relative invariants.

The basic implication needed later is the following contact-invariance statement.

\begin{proposition}[Contact invariance of filtered Jacobian quotients]\label{matheryaug}
Let $A$ be a smooth finitely generated $\mathbf C$-algebra or a localization
of one. Let $B=A[[\mathbf u]]$, where
$\mathbf u=(u_1,\dots,u_m)$, and let
$\mathfrak m_u=\langle u_1,\dots,u_m\rangle$. If $f,g\in B$ are
contact equivalent over $A$, then for every $k\geq0$ there is an $A$-algebra
isomorphism
\[
\frac{B}{\langle f,\mathfrak m_u^kJ_{\mathbf u}(f)\rangle}
\cong
\frac{B}{\langle g,\mathfrak m_u^kJ_{\mathbf u}(g)\rangle}.
\]
\end{proposition}

\begin{proof}
Write $f=U\cdot\varphi(g)$ with $\varphi\in\operatorname{Aut}_A(B)$ and
$U\in B^*$. The automorphism $\varphi$ preserves the $\mathfrak m_u$-adic
filtration. By the formal chain rule, the ideal generated by the relative
derivatives of $\varphi(g)$ is the image under $\varphi$ of the relative
Jacobian ideal of $g$. Multiplication by the unit $U$ does not change the
Tjurina-type ideal, since the product rule gives
\[
\langle U\varphi(g),\mathfrak m_u^kJ_{\mathbf u}(U\varphi(g))\rangle
=
\langle \varphi(g),\mathfrak m_u^kJ_{\mathbf u}(\varphi(g))\rangle.
\]
Hence
\[
\langle f,\mathfrak m_u^kJ_{\mathbf u}(f)\rangle
=
\varphi\bigl(\langle g,\mathfrak m_u^kJ_{\mathbf u}(g)\rangle\bigr),
\]
and $\varphi$ induces the displayed quotient isomorphism.
\end{proof}

\begin{example}[Failure of the relative converse over a smooth base]\label{ex:relative_converse_failure}
Even for a smooth base algebra of the form used above, an isomorphism between
Tjurina-type quotients of $B$ need not preserve the relative formal ideal. Take
the base to be the localized polynomial $\mathbf C$-algebra
\[
A=\mathbf C[t]_S=\mathbf C[t]_{(t)},
\qquad S=\mathbf C[t]\setminus(t),
\]
and set
\[
B=A[[u]],
\qquad \mathfrak m_u=\langle u\rangle B.
\]
Let
\[
f=u(u-t)(u-2t),\qquad g=u(u-t)(u+t).
\]
The polynomial change of variable $u\mapsto u+t$ is an $A$-algebra automorphism
of $A[u]$ and satisfies
\[
f(u+t)=g,\qquad J_u(f)(u+t)=J_u(g).
\]
We next justify that the polynomial Tjurina quotients agree with the
formal-power-series quotients. A direct resultant computation gives
\[
\operatorname{Res}_u(f,J_u(f))
=
\operatorname{Res}_u(g,J_u(g))
=-4t^6.
\]
Consequently, $t^6=0$ in both quotient rings. Moreover,
\[
J_u(g)=3u^2-t^2,
\qquad
J_u(f)=3u^2-6tu+2t^2.
\]
Thus $u^2$ belongs to $tB/(f,J_u(f))$ in the first quotient and to
$t^2B/(g,J_u(g))$ in the second. Since $t$ is nilpotent there, $u$ is
nilpotent as well. The polynomial quotients are finite $A$-modules because
$J_u(f)$ and $J_u(g)$ are monic up to the unit $3$, and they are annihilated
by a power of $(u)$. Hence they are already complete for the $(u)$-adic
topology. Since $B=A[[u]]$ is the $(u)$-adic completion of $A[u]$, exactness
of completion for finite modules gives the canonical base-change isomorphisms
\[
\frac{A[u]}{(f,J_u(f))}
\xrightarrow{\sim}
\frac{B}{(f,J_u(f))},
\qquad
\frac{A[u]}{(g,J_u(g))}
\xrightarrow{\sim}
\frac{B}{(g,J_u(g))}.
\]
Equivalently, because a power of $u$ vanishes in each quotient, every formal
series has a polynomial representative and completion introduces no additional
classes.
Therefore, for $k=0$, the polynomial translation induces an $A$-algebra
isomorphism
\[
\frac{B}{\langle f,J_u(f)\rangle}
\cong
\frac{B}{\langle g,J_u(g)\rangle}.
\]
However, this change of variable is not a relative formal coordinate change on
$B$: it sends $u$ to $u+t$ and therefore does not preserve $\mathfrak m_u$.
Moreover $f$ and $g$ are not contact equivalent over $A$ in the relative formal
sense. Indeed, any such equivalence preserves $\mathfrak m_u$; after completing
at the maximal ideal $\langle t,u\rangle$, its tangent-cone linear part has the
form $u\mapsto a u$ with $a\in\mathbf C^*$. Hence the initial cubic form of
$f$, namely $u(u-t)(u-2t)$, would have to be proportional to that of
$\varphi(g)$, namely
\[
g(au)=au(au-t)(au+t).
\]
The root slopes of $f$ and $g$ are respectively $\{0,1,2\}$ and $\{0,1,-1\}$,
and no scaling of the $u$-coordinate sends one set to the other. Thus an
isomorphism of Tjurina-type quotients over a smooth localized polynomial base
does not by itself imply relative contact equivalence.
\end{example}

This example is the point at which the absolute and relative theories diverge.
Theorem~\ref{matheryauoriginal} may recover an absolute analytic contact
equivalence from Tjurina data, but it does not force the resulting automorphism
to fix the base or preserve $\mathfrak m_u$. For this reason the categorical
statements below are formulated directly in terms of formal contact equivalence.

For $i=1,2$, let $W_i:X_i\to\mathbf A^1_{\mathbf C}$ be a potential, and put
$Y_i=W_i^{-1}(w_i)$. Analytic Tjurina data can identify local contact types via
Theorem~\ref{matheryauoriginal}; once the resulting formal contact equivalences
are fixed along the critical components, the derived categories of matrix factorizations are
identified. The reverse problem, recovering an isolated hypersurface
singularity from its matrix-factorization category, is addressed by the
categorical Mather--Yau theorem of Hua and Keller \cite{HuaKeller2024}; see
also \cite{Hua2018NoncommutativeMY,Keller2018SingularHH,Kalck2021MFClassification}.
It is not needed here.

\begin{corollary}[Componentwise formal contact invariance]\label{mather-yau}
For $i=1,2$, let $X_i$ be a smooth quasi-projective variety over $\mathbf C$,
let $W_i:X_i\to\mathbf A^1_{\mathbf C}$ be flat, and put
$Y_i=W_i^{-1}(w_i)$. Assume that the reduced singular support
\[
Z_i=\operatorname{Sing}(Y_i)_{\mathrm{red}}
\]
is a finite union of smooth closed subvarieties that are pairwise disjoint as
closed subschemes
\[
Z_i=\coprod_{\alpha\in A_i}Z_{i,\alpha}.
\]
Suppose there is a bijection
$\sigma:A_1\to A_2$ such that the corresponding formal completions are
contact-equivalent componentwise. More precisely, if
\[
\widehat X_{i,\alpha}:=\widehat{X_i}_{Z_{i,\alpha}}
\]
and $\widehat W_{i,\alpha}$ denotes the induced formal potential, assume that
there is an isomorphism of formal schemes
\[
\phi_\alpha:\widehat X_{1,\alpha}\xrightarrow{\sim}\widehat X_{2,\sigma(\alpha)}
\]
and a unit
$u_\alpha\in\Gamma(\widehat X_{1,\alpha},\mathcal O_{\widehat X_{1,\alpha}}^{*})$
such that
\[
\phi_\alpha^{*}(\widehat W_{2,\sigma(\alpha)}-w_2)=u_\alpha(\widehat W_{1,\alpha}-w_1).
\]
Then the idempotent-completed derived categories of matrix factorizations are equivalent:
\[
\overline{\MF_{w_1}(X_1,W_1)}
\cong
\overline{\MF_{w_2}(X_2,W_2)}.
\]
\end{corollary}
\begin{proof}
Because the components of the reduced supports $Z_1$ and $Z_2$ are pairwise
disjoint as closed subschemes, their formal completions decompose
componentwise:
\[
\widehat{X_1}_{Z_1}\cong\coprod_{\alpha\in A_1}\widehat X_{1,\alpha},
\qquad
\widehat{X_2}_{Z_2}\cong\coprod_{\beta\in A_2}\widehat X_{2,\beta}.
\]
For each $\alpha$, the unit $u_\alpha$ does not change the principal ideal
defining the zero fiber. Thus $\phi_\alpha$ induces an isomorphism between the
corresponding completed zero fibers.
Corollary~\ref{cor:supported_formal_local} therefore gives
\[
\overline{\MF^{Z_{1,\alpha}}_{w_1}(X_1,W_1)}
\cong
\overline{\MF^{Z_{2,\sigma(\alpha)}}_{w_2}(X_2,W_2)}.
\]
Taking direct sums over $\alpha$ and applying
Proposition~\ref{prop:disjoint_supports} to both singular loci proves the
claim.
\end{proof}
For a fixed smooth support $Z\subset X$,
Corollary~\ref{cor:supported_formal_local}
reduces the comparison of supported categories to a comparison of formal
potentials on $\widehat X_Z$. Locally, $\widehat X_Z$ has the form
\[
\operatorname{Spf}\left(\frac{\mathbf{C}[x_1,\dots,x_n]_S}{I}[[\mathbf u]]\right),
\qquad V(I)\cong Z.
\]
This gives the corresponding supported version of contact invariance.
\begin{corollary}[Supported contact invariance]\label{cor:supported_contact}
For $i=1,2$, let $X_i$ be a smooth quasi-projective variety over $\mathbf C$,
let $W_i:X_i\to\mathbf A^1_{\mathbf C}$ be flat, put
$Y_i=W_i^{-1}(w_i)$, and let $Z_i\subset X_i$ be a closed reduced support.
Suppose that there are an isomorphism of formal schemes
\[
\phi:\widehat{X_1}_{Z_1}\xrightarrow{\sim}\widehat{X_2}_{Z_2}
\]
and a unit $U\in\Gamma(\widehat{X_1}_{Z_1},\mathcal O^*)$ such that
\[
\phi^*(\widehat W_2-w_2)=U(\widehat W_1-w_1).
\]
Then
\[
\overline{\MF^{Z_1}_{w_1}(X_1,W_1)}
\simeq
\overline{\MF^{Z_2}_{w_2}(X_2,W_2)}.
\]
\end{corollary}

\begin{proof}
The formal contact equivalence identifies the completed zero fibers. Put
$R_i=Z_i\cap\operatorname{Sing}(Y_i)$. Since the singular locus is intrinsic
to the zero fiber, the induced isomorphism carries the formal completion along
$R_1$ to that along $R_2$, equivalently along their reductions. The supported
formal-local principle, Corollary~\ref{cor:supported_formal_local}, gives the
claimed equivalence after idempotent completion.
This proof remains in the supported, equivalently torsion-formal, categories;
it does not use the unrestricted formal category and therefore requires no
additional hypothesis of the form
$\operatorname{Sing}(\operatorname{Spec}B)\subseteq V(J)$.
\end{proof}

\begin{lemma}[Pointwise splitting normal form]
\label{lem:splitting_normal_form}
Let $z$ lie on a smooth component $Z$ of the critical locus.  On a transverse analytic slice at
$z$, assume that the Hessian has a nondegenerate block of rank $r$ and that the
remaining one-variable germ is of type $A_{n-1}$.  Then the germ on that slice
is formally right equivalent to
\[
q(v_1,\dots,v_r)+u^n,
\]
where $q$ is nondegenerate.
\end{lemma}

\begin{proof}
The formal splitting lemma separates the nondegenerate quadratic variables from
the residual germ, and the one-variable classification gives $u^n$; see
\cite[Theorems~2.47 and~2.48]{greuel2007introduction}.
\end{proof}

\begin{proposition}[Global verification criterion for split normal forms]
\label{prop:global_normal_form_checklist}
Lemma~\ref{lem:splitting_normal_form} is a pointwise statement on a completed
transverse slice. It does not produce a single formal coordinate system on
$\widehat X_Z$, nor does it imply a global identity
\[
\widehat W-W|_Z=q+u^n
\]
over all of $Z$. The nondegenerate Hessian directions form a quadratic vector
bundle, whose Clifford algebra may be nontrivial, while the residual corank-one
direction may be a nontrivial line bundle $L$. In addition, the formal
neighborhood of $Z$ in $X$ need not split as the completion of the normal
bundle. Kapustin, Rozansky, and Saulina proposed a Clifford description of
the Morse--Bott category without a splitting hypothesis
\cite{KapustinRozanskySaulina2009}; Teleman showed that in the nonsplit case
extra deformation data may occur \cite{TelemanMorseBott}. His corrected
statement proves the expected Clifford description in the split-first-
neighborhood case and explains why it can fail without that hypothesis.

A global split $A_{n-1}$ or Morse--Bott identification is justified once the
following data have been supplied:
\begin{enumerate}[label=\textup{(\roman*)}]
\item an explicit splitting of the required formal neighborhood of $Z$ (or an
equivalent system of global formal normal coordinates);
\item a splitting of the normal bundle into the nondegenerate Hessian directions
and the residual directions, with the Hessian of constant rank and its
quadratic form recorded as a quadratic vector bundle $(E,q)$;
\item in the corank-one case, a residual line bundle $L$ and a relative formal
coordinate change over $Z$ eliminating all higher mixed terms;
\item the leading residual coefficient as a section of
$(L^\vee)^{\otimes n}$, together with the line-bundle and unit data needed to
write a global residual term; and
\item the Clifford algebra $\Cliff(E,q)$, its super-Brauer class, and any
nonsplit deformation class that remains.
\end{enumerate}
Under these hypotheses the resulting global formal contact equivalence permits
the use of Corollary~\ref{cor:supported_contact}. Without them, the pointwise
normal form alone justifies only an unspecified formal supported category
$\FormalCat{Z}$.  For a pure residual potential $a u^n$ with $a$ a global
unit, contact multiplication by $a^{-1}$ removes the coefficient, so no
$n$th root is required.  In a mixed model $b v^2+a u^n$, however, the same
contact unit acts on both summands, and the relative unit class need not be
removable.  Likewise, a nondegenerate quadratic coefficient is not
automatically a global square; the associated line-bundle, mixed-unit, and
Clifford data then enter the normal form.

If the residual line bundle is nontrivial, the natural nilpotent thickening is
\[
Z^{(n)}_L=
\operatorname{Spec}_Z\!\left(
\mathcal O_Z\oplus L^\vee\oplus\cdots\oplus(L^\vee)^{\otimes(n-1)}
\right),
\]
with multiplication truncated in degree $n$, rather than the trivial
thickening defined by a global coordinate $u$.
\end{proposition}

\begin{proof}
Item \textup{(i)} identifies the completed structure sheaf with the completed
symmetric algebra of the conormal directions. Item \textup{(ii)} writes the
quadratic part as the globally defined form $q$ on $E$ and separates its
orthogonal complement. Since $q$ is nondegenerate, the relative formal
splitting lemma removes, degree by degree in the adic filtration, every term
involving a nondegenerate variable beyond the quadratic part. In the
corank-one complement, item \textup{(iii)} performs the remaining formal
one-variable changes over $Z$, while item \textup{(iv)} ensures that the
leading coefficient is globally defined and invertible. These compatible
changes converge in the complete adic topology and give a formal contact
equivalence over all of $Z$. Item \textup{(v)} records the quadratic and
deformation data not removed by the equivalence. The conclusion follows from
Corollary~\ref{cor:supported_contact}; without gluing of the normal
coordinates, the construction exists only pointwise.
\end{proof}

\begin{lemma}[Affine complete-intersection Morse--Bott criterion]
\label{lem:affine_morse_bott_criterion}
Let $X=\operatorname{Spec}A$ be smooth affine over $\mathbf C$, and let
$Z=\operatorname{Spec}(A/I)$ be a smooth closed affine subvariety of
codimension $c$. Suppose that $I=(f_1,\ldots,f_c)$ and that the normal Hessian
of $W$ along $Z$ is nondegenerate, and suppose that $Z\subseteq\CritLoc(W)$.
Then the completion of $X$ along $Z$ admits
a global formal splitting
\[
\widehat A_I\simeq (A/I)[[u_1,\ldots,u_c]],
\]
under which $W-W|_Z$ is formally right equivalent over $Z$ to its
nondegenerate quadratic normal form. In particular, $Z$ is a genuine split
Morse--Bott component and its contribution is $\MBcat{Z}$. On each connected
component, for its corresponding critical value, the singular locus of the
completed fiber is contained in the zero section. Thus the support hypothesis
of Theorem~\ref{thm:general_formal} is automatic in this quadratic-normal
situation.
\end{lemma}

\begin{proof}
Put $R=A/I$.  Since $R$ is smooth over $\mathbf C$, it is formally smooth.
Starting with the identity $R\to A/I$, suppose inductively that a section
$s_n:R\to A/I^n$ has been constructed.  The surjection
$A/I^{n+1}\twoheadrightarrow A/I^n$ has nilpotent kernel, so formal smoothness
lifts $s_n$ to a section $s_{n+1}:R\to A/I^{n+1}$.  Passing to the inverse
limit gives a continuous section
\[
s:R\longrightarrow\widehat A_I.
\]

Because $X$ and $Z$ are smooth, the embedding $Z\hookrightarrow X$ is regular
of codimension $c$.  The conormal module $I/I^2$ is therefore locally free of
rank $c$.  The chosen generators induce a surjection
$R^{\oplus c}\twoheadrightarrow I/I^2$ between locally free modules of the
same rank, hence an isomorphism.  In particular, the $f_i$ form a regular
sequence and
\[
\operatorname{gr}_I(A)\simeq
\operatorname{Sym}_{R}(I/I^2)
\simeq R[u_1,\ldots,u_c].
\]
Using the section $s$, define a continuous filtered homomorphism
\[
R[[u_1,\ldots,u_c]]\longrightarrow\widehat A_I,
\qquad u_i\longmapsto f_i.
\]
It is an isomorphism on associated graded rings.  Since both sides are
complete and separated for their adic filtrations, it is an isomorphism.

Because $Z\subseteq\CritLoc(W)$, the expansion of $W-W|_Z$ begins in degree
two in the variables $u_i$.  Its quadratic part is the normal Hessian and is
nondegenerate by hypothesis.  The relative formal Morse argument removes
successively every homogeneous term of degree at least three: at each degree,
the inverse Hessian determines a homogeneous formal vector field whose
coordinate change cancels that term.  Characteristic zero and adic
completeness make the recursive coordinate changes converge formally.  Thus
$W-W|_Z$ is formally right equivalent over $R$ to its quadratic normal form.
This is the affine coefficient-ring form of the relative formal Morse lemma;
compare \cite[Lemma~A.2]{HennionHolsteinRobalo2024}.
\end{proof}

\begin{definition}[The periodic $A$-type category]\label{def:periodic_A_type}
Let $Z$ be a smooth variety and put
\[
Z^{(n)}=\operatorname{Spec}_Z(\mathcal O_Z[u]/(u^n)).
\]
We define
\[
\Aper{Z}{n-1}:=\overline{\Db_{\mathrm{sg}}(Z^{(n)})}.
\]
Equivalently, when a formal line coordinate $u$ has been fixed, this is the
idempotent completion of the ungraded two-periodic matrix-factorization category
of the formal potential $u^n$ over $Z$.
\end{definition}

The notation $\Aper{Z}{n-1}$ therefore includes the choice of a trivial global
formal line coordinate. For a nontrivial residual line bundle, one must instead
use the singularity category of $Z^{(n)}_L$ from
Proposition~\ref{prop:global_normal_form_checklist}.
If the completed potential is $a u^n$ for a global unit $a$, contact
multiplication by $a^{-1}$ identifies it with $u^n$; this pure residual case
does not require an $n$th root of $a$.

\begin{proposition}[The genuine ungraded $A$-type model]
\label{prop:Ade_category_model}
If $\widehat X_Z\simeq\operatorname{Spf}(\mathcal O_Z[[u]])$ and
$\widehat W=u^n$, then
\[
\overline{\MF(\widehat X_Z,u^n)}
\simeq
\Aper{Z}{n-1}.
\]
For $Z=\operatorname{Spec}\mathbf C$, this category is equivalent to the stable
category of finite modules over $\mathbf C[u]/(u^n)$; its indecomposable
nonprojective modules are
\[
\mathbf C[u]/(u),\ldots,\mathbf C[u]/(u^{n-1}).
\]
\end{proposition}

\begin{proof}
By the convention in Remark~\ref{rem:formal_mf_meaning}, the notation on the
left means the singularity category of the formal zero fiber; no
scheme-theoretic matrix-factorization theorem is being applied directly to
the formal scheme. That zero fiber is precisely the ordinary nilpotent
thickening $Z^{(n)}$, so the first equivalence is immediate from
Definition~\ref{def:periodic_A_type}. For
$Z=\operatorname{Spec}\mathbf C$, the hypersurface ring
$\mathbf C[u]/(u^n)$ is Frobenius, so its singularity category is its stable
module category.  The list of indecomposables is the standard classification of
modules over the principal Artinian local ring $\mathbf C[u]/(u^n)$.
This one-variable description is compatible with the standard ADE
matrix-factorization classification \cite{kajiura2007matrix}.
\end{proof}

\begin{remark}
Proposition~\ref{prop:Ade_category_model} must not be inferred from pointwise
normal forms alone: its hypotheses assert a global splitting of the entire
formal neighborhood and a global equality of potentials. When only the
pointwise transverse type is known, we do not identify the contribution with
a standard category.
\end{remark}

\begin{definition}[Mixed quadratic--$A$-type categories]
\label{def:mixed_quadratic_A_type}
Let $Z$ be smooth and let $m\geq 1$.  The standard split mixed category is
\[
\QAper{Z}{m}:=
\overline{\MF_0\!\left(
\operatorname{Spf}(\mathcal O_Z[[v,u]]),v^2+u^{m+1}
\right)}.
\]
More generally, for $\delta\in\Gamma(Z,\mathcal O_Z^*)$, put
\[
\QAtwist{Z}{m}{\delta}:=
\overline{\MF_0\!\left(
\operatorname{Spf}(\mathcal O_Z[[v,u]]),v^2+\delta u^{m+1}
\right)}.
\]
If $\delta'=h^2\delta$ for a unit $h$, multiplication of the potential by
$h^2$ followed by the change $v\mapsto hv$ identifies the two formal contact
models.  Thus the second notation depends only on the displayed square class.
In particular, a unit coefficient may be omitted only after its relevant
global class has been trivialized.
\end{definition}

For both mixed models, the singular locus of the completed zero fiber is the
zero section $u=v=0$. Hence, when the formal fiber is affine, the torsion and
unrestricted formal categories agree after idempotent completion by
Theorem~\ref{thm:general_formal}.

\begin{definition}[The twisted periodic $D_4$ family]
\label{def:twisted_periodic_D4}
Let $Z$ be smooth and let $\delta\in\Gamma(Z,\mathcal O_Z^*)$.  Define
\[
\Dper{Z}{4}{\delta}:=
\overline{\MF_0\!\left(
\operatorname{Spf}(\mathcal O_Z[[p,v]]),
p(p^2-\delta v^2)
\right)}.
\]
If $\delta'=h^2\delta$ for a unit $h$, the substitution
$v\mapsto h^{-1}v$ identifies the two formal models.  Thus the notation
depends only on the square class $[\delta]$.  In the applications below, the
singular locus of the displayed formal hypersurface is the zero section, so
the supported and unrestricted formal categories agree after applying
Theorem~\ref{thm:general_formal}.
\end{definition}

\begin{definition}[Quadratic local categories]\label{def:quadratic_categories}
For a nondegenerate quadratic form of rank $r$ over $\mathbf C$, write
\[
\MorseCat{r}:=\operatorname{Perf}_{\mathbf Z/2}(\Cliff_r)
\]
for the two-periodic category of finite-dimensional supermodules over its
complex Clifford algebra.  For a Morse--Bott component $Z$ with quadratic
normal bundle $(N,q)$, write
\[
\MBcat{Z}:=\Dtwo(Z;\Cliff(N,q)),
\]
with any super-Brauer twist determined by the normal bundle and the chosen
splitting of the formal neighborhood understood as part of the notation.
\end{definition}

\begin{remark}[Meaning of the square-class labels]
The notation $\MBcat{Z;[\delta]}$ is shorthand for the full Clifford category
$\Dtwo(Z;\Cliff(N,q))$, together with the explicitly computed discriminant
square class $[\delta]$ of the quadratic normal bundle.  We do not claim that
the discriminant alone classifies an arbitrary quadratic bundle over an
arbitrary base.  The global normal coordinates and the full quadratic form
remain part of the datum; the square class is displayed because it detects the
twisting relevant in the examples below. This is the standard
Clifford/super-Brauer framework for quadratic categories
\cite{BuchweitzEisenbudHerzog1987,DonovanKaroubi1970,KuznetsovQuadric2008,AuelBernardaraBolognesi2014,TelemanMorseBott}.
\end{remark}

\begin{proposition}[Morse and Morse--Bott contributions]
\label{prop:morse_bott_clifford}
An isolated nondegenerate critical point in an ambient space of dimension $r$
contributes $\MorseCat{r}$.  In the split Morse--Bott situation, the formal
contribution of a smooth critical component is the corresponding two-periodic
Clifford-module category $\MBcat{Z}$; see
\cite{BuchweitzEisenbudHerzog1987,KapustinRozanskySaulina2009,TelemanMorseBott}.
\end{proposition}

\begin{proof}
The formal Morse lemma identifies an isolated nondegenerate critical point
with a nondegenerate quadratic form. Its singular locus is the closed point,
so the support criterion of Theorem~\ref{thm:general_formal} is automatic,
and the standard Clifford description gives $\MorseCat{r}$. In the split
Morse--Bott case, the relative quadratic normal form has singular locus equal
to the zero section. The same criterion is therefore automatic, and the
relative Clifford description gives $\MBcat{Z}$.
\end{proof}

\begin{remark}
The splitting hypothesis is essential in the preceding proposition. In a
nonsplit formal neighborhood, a Clifford or super-Brauer twist alone need not
describe the matrix-factorization category; an additional deformation class
may occur. Teleman attributes the motivating obstruction to an example
observed by Pomerleano and Preygel \cite{TelemanMorseBott}. Thus
$\MBcat{Z}$ is used below only when the required relative splitting has been
verified.
\end{remark}

\begin{remark}[Kn\"orrer periodicity and parity]\label{rem:knorrer_shift_convention}
Standard Kn\"orrer periodicity removes a hyperbolic pair:
\[
\MF(f)\simeq \MF(f+xy),
\]
under the usual hypotheses.  A single quadratic variable is instead encoded by
a rank-one Clifford factor and cannot in general be replaced by a mere shift
$[1]$.  Consequently, all formulas below are written in terms of the actual
two-periodic Clifford or periodic $A$-type categories; see the original
periodicity theorem \cite{knorrer1987cohen}.
\end{remark}
\begin{remark}[Hochschild invariants and the Morse--Bott twist]
\label{rem:HH_MorseBott_twist}
For clarity, we use the Lin--Pomerleano definition of Hochschild cohomology:
$HH^\bullet(\MF_{0}(X,W))$ is the derived endomorphism algebra of
the identity functor. Their diagonal-kernel description and HKR computation
identify it, for smooth affine $X$, with
\[
HH^\bullet(\MF_{0}(X,W))
\simeq
H^\bullet\!\left(
\Gamma(X,\textstyle\bigwedge^\bullet T_X),\iota_{dW}
\right),
\]
where $\iota_{dW}$ is contraction by $dW$. Likewise, their Hochschild
homology model is
\[
HH_\bullet(\MF_{0}(X,W))
\simeq
H^\bullet\!\left(
\Gamma(X,\Omega_X^\bullet),dW\wedge -
\right);
\]
see \cite[Theorem~3.1 and Remark~3.2]{LinPomerleano2013}.

Hochschild homology and cohomology are invariants of the dg Morita
equivalence class; see \cite[Secs.~4.1, 5.3, and 5.4]{Keller2006dg}.
Consequently,
\[
        \overline{\MF_{\lambda_f}(X,f)}\cong \overline{\MF_{\lambda_g}(Y,g)} \implies HH_\bullet\left(\overline{\MF_{\lambda_f}(X,f)}\right)\cong HH_\bullet\left(\overline{\MF_{\lambda_g}(Y,g)}\right).
\]

The converse fails, however. For example, let
\[
\Gamma=\operatorname{Spec}\mathbf C[t,t^{-1}],
\qquad
X=\Gamma\times\mathbf A^2
=\operatorname{Spec}\mathbf C[t,t^{-1},u,v],
\]
and consider the two potentials
\[
W_0=uv,\qquad W_1=u^2-tv^2.
\]
In both cases the critical locus, scheme-theoretically, is
\[
\CritLoc(W_i)=\Gamma\times\{(0,0)\}\subset W_i^{-1}(0).
\]
The first normal quadratic form is split, while the second has discriminant
class $[t]\in\mathcal O(\Gamma)^*/\mathcal O(\Gamma)^{*2}$. This class is
nontrivial: the unit $t$ is not a square in
$\mathbf C[t,t^{-1}]^*$. Accordingly, the two Morse--Bott contributions are
\[
\overline{\MF_0(X,W_0)}\simeq\MBcat{\Gamma},
\qquad
\overline{\MF_0(X,W_1)}\simeq\MBcat{\Gamma;[t]},
\]
with the latter being the twisted model appearing for the curve component of
$MM_{2-12}$.

For $W_0$, this is the Koszul complex on the regular sequence $(u,v)$,
together with the unpaired tangent direction $\partial_t$. For $W_1$, set
\[
e_1=du,\qquad
e_2=dv+\frac{v}{2t}\,dt,\qquad
e_3=dt.
\]
These form a basis of $\Omega_X^1$, and
\[
dW_1=2u\,e_1-2tv\,e_2.
\]
The dual basis has unpaired tangent vector
\[
e_3^\vee=\partial_t-\frac{v}{2t}\partial_v,
\]
which restricts to $\partial_t$ on $\Gamma$. Since $t$ is invertible, the
two Koszul coefficients $2u$ and $-2tv$ form a regular sequence with
quotient $\mathbf C[t,t^{-1}]$. Thus in both cases
\[
HH^\bullet(\MF_{0}(X,W_i))
\cong
\mathbf C[t,t^{-1}]\oplus
\mathbf C[t,t^{-1}]\,\epsilon,
\qquad
|\epsilon|=1,\quad \epsilon^2=0,
\]
as $\mathbf Z/2$-graded algebras.

The homology calculation gives the same answer as a module over
$HH^0=\mathbf C[t,t^{-1}]$. More precisely, if $\omega_i$ denotes the
normal-volume class represented by $du\wedge dv$ for $W_0$ and by
$e_1\wedge e_2$ for $W_1$, then
\[
HH_\bullet(\MF_{0}(X,W_i))
\cong
\mathbf C[t,t^{-1}]\,\omega_i
\oplus
\mathbf C[t,t^{-1}]\,(\omega_i\wedge dt).
\]
After identifying $\omega_i$ with $1$, this is
$\mathbf C[t,t^{-1}]\oplus\mathbf C[t,t^{-1}]\,dt$ for each potential.
Consequently, the Hochschild cohomology algebras and the Hochschild
homology modules do not distinguish the two categories.

For an isolated hypersurface singularity $W\in R=\mathbf C[[x_1,\ldots,x_n]]$, let
\[
J_W=(\partial_{x_1}W,\ldots,\partial_{x_n}W)
\]
denote its Jacobian ideal. Dyckerhoff computes the Hochschild cohomology
of its matrix-factorization category without assuming quasihomogeneity:
\[
HH^\bullet(\overline{\MF(R,W)})\cong R/J_W,
\]
concentrated in even parity \cite[Corollary~6.5]{Dyckerhoff2011}.

If $W$ is weighted homogeneous of degree $d$ with positive weights
$q_1,\ldots,q_n$, the weighted Euler identity gives
\[
dW=\sum_{i=1}^n q_i x_i\partial_{x_i}W,
\]
so $W\in J_W$. Consequently, the Tjurina algebra
$R/(W,J_W)$ coincides with the Milnor (Jacobian) algebra $R/J_W$
\cite{Saito1971Quasihomogeneous}.

More generally, the Tjurina and Milnor sheaves agree even when
$\CritLoc(W)$ is positive-dimensional, provided $W$ is locally
quasi-homogeneous. On the analytic space $X$, let $\mathcal J_W$ be the
Jacobian ideal sheaf and set
\[
\mathcal M_W=\mathcal O_X/\mathcal J_W,\qquad
\mathcal T_W=\mathcal O_X/(\mathcal J_W+(W)).
\]
If each germ of $W$ along $\CritLoc(W)$ is analytically equivalent to a
weighted-homogeneous germ with positive weights, the local Euler identity
gives $W\in\mathcal J_{W,p}$ at every critical point $p$. Hence
$\mathcal T_W\cong\mathcal M_W$ as sheaves and
$\Gamma(X,\mathcal T_W)\cong\Gamma(X,\mathcal M_W)$. More generally, only
the stalkwise inclusion $W\in\mathcal J_{W,p}$ is needed; the local Euler
fields need not glue globally. Relative Milnor/Tjurina invariants for
quasihomogeneous functions, including non-isolated cases, are discussed in
\cite{AhmedRuas2011Relative}.

In the example above, this conclusion
applies fiberwise to the transverse quadratic germs at fixed $t\in\Gamma$;
the full potentials on $X$ still have the non-isolated critical locus
$\Gamma$.
\end{remark}
\section{Decomposition of Landau--Ginzburg Models of del Pezzo Surfaces}
\label{sec:dp_decompositions}
We now record the del Pezzo critical-locus computations on the affine torus
$(\mathbf C^*)^2$.  We use the specified anticanonical Laurent-polynomial
models, not general members of the corresponding Landau--Ginzburg deformation
spaces \cite{auroux2006mirror,LuntsPrzyjalkowski2018,Przyjalkowski2026RankOne}.
For a general model of a degree-$d$ del Pezzo surface, the compactified
Landau--Ginzburg pencil has $12-d$ fibers with one ordinary double point each;
under specialization to the anticanonical model, isolated singularities can
collide into a reducible or positive-dimensional critical fiber. This is the
context for the smaller open-torus counts and residual components below. The
categorical notation follows Section~\ref{sec:mather_yau}: an isolated Morse
point in a surface contributes $\MorseCat{2}$; a verified split quadratic
critical curve carries a Morse--Bott/Clifford contribution $\MBcat{Y}$; and a
verified global formal normal form $u^{m+1}$ contributes $\Aper{Y}{m}$.

\begin{example}[$\mathbf{P}^2$]
Consider
\[
W_{\mathbf{P}^2}=x+y+\frac{1}{xy}.
\]
The critical-point equations have the three solutions
\[
\CritLoc(W_{\mathbf{P}^2})=
\{(\omega_3^i,\omega_3^i)\mid i=0,1,2\},
\qquad
\omega_3=e^{2\pi i/3},
\]
with critical values
\[
\CritVal(W_{\mathbf{P}^2})=\{3\omega_3^i\mid i=0,1,2\}.
\]
The Hessian is nondegenerate at all three critical points.  Hence
\[
\overline{\MF((\mathbf{C}^*)^2,W_{\mathbf{P}^2})}
\cong \bigoplus_{i=1}^{3}\MorseCat{2}.
\]
\end{example}

\begin{example}[$\mathbf{P}^1\times\mathbf{P}^1$]
Consider
\[
W_{\mathbf{P}^1\times\mathbf{P}^1}=x+y+\frac{1}{x}+\frac{1}{y}.
\]
The critical locus is
\[
\CritLoc(W_{\mathbf{P}^1\times\mathbf{P}^1})=
\{(1,1),(-1,-1),(-1,1),(1,-1)\},
\]
and the critical values are
\[
\CritVal(W_{\mathbf{P}^1\times\mathbf{P}^1})=\{4,0,-4\}.
\]
All four critical points are Morse, so
\[
\overline{\MF((\mathbf{C}^*)^2,W_{\mathbf{P}^1\times\mathbf{P}^1})}
\cong \bigoplus_{i=1}^{4}\MorseCat{2}.
\]
\end{example}

\begin{example}[$\mathrm{Bl}_1\mathbf{P}^2$]
Consider
\[
W_{\mathrm{Bl}_1\mathbf{P}^2}=x+y+xy+\frac{1}{xy}.
\]
The critical locus is
\[
\CritLoc(W_{\mathrm{Bl}_1\mathbf{P}^2})=
\{(a,a)\in(\mathbf{C}^*)^2\mid a^4+a^3-1=0\}.
\]
If $a_1,\ldots,a_4$ are the roots of $a^4+a^3-1=0$, then the critical values
are $2a_i^2+3a_i$.  The four critical points are Morse, and therefore
\[
\overline{\MF((\mathbf{C}^*)^2,W_{\mathrm{Bl}_1\mathbf{P}^2})}
\cong \bigoplus_{i=1}^{4}\MorseCat{2}.
\]
\end{example}

\begin{example}[$\mathrm{Bl}_2\mathbf{P}^2$]
Consider
\[
W_{\mathrm{Bl}_2\mathbf{P}^2}=(1+x+y)\left(1+\frac{1}{xy}\right)-1.
\]
The critical locus is the disjoint union
\[
\CritLoc(W_{\mathrm{Bl}_2\mathbf{P}^2})=
\{(a,a)\mid a^3-a-1=0\}
\cup
\{(a,b)\mid a^2+a-1=0,\ ab=-1\}.
\]
The corresponding critical values are
\[
\{a^2+3a-1\mid a^3-a-1=0\}\cup\{-1\}.
\]
The five critical points are Morse.  Hence
\[
\overline{\MF((\mathbf{C}^*)^2,W_{\mathrm{Bl}_2\mathbf{P}^2})}
\cong \bigoplus_{i=1}^{5}\MorseCat{2}.
\]
\end{example}

\begin{example}[$\mathrm{Bl}_3\mathbf{P}^2$]
Consider
\[
W_{\mathrm{Bl}_3\mathbf{P}^2}=(1+x)(1+y)\left(1+\frac{1}{xy}\right)-2.
\]
The critical locus is
\[
\CritLoc(W_{\mathrm{Bl}_3\mathbf{P}^2})=
\{(-1,-1),(1,-1),(-1,1),(1,1),(\omega,\omega),(\omega^2,\omega^2)\},
\]
where $\omega=e^{2\pi i/3}$.  The critical values are
\[
\CritVal(W_{\mathrm{Bl}_3\mathbf{P}^2})=\{-3,-2,6\}.
\]
Each of the six critical points has completed normal form $u^2+v^2$ after
subtracting its critical value.  Therefore
\[
\overline{\MF((\mathbf{C}^*)^2,W_{\mathrm{Bl}_3\mathbf{P}^2})}
\cong \bigoplus_{i=1}^{6}\MorseCat{2}.
\]
\end{example}

\begin{example}[$\mathrm{Bl}_4\mathbf{P}^2$]
Consider
\[
W_{\mathrm{Bl}_4\mathbf{P}^2}=
(1+x+y)\left(1+\frac{1}{x}\right)\left(1+\frac{1}{y}\right)-3.
\]
The critical locus consists of the three points
\[
(-1,-1),\qquad
\left(\frac{1+\sqrt{5}}{2},\frac{1+\sqrt{5}}{2}\right),\qquad
\left(\frac{1-\sqrt{5}}{2},\frac{1-\sqrt{5}}{2}\right),
\]
with critical values
\[
\CritVal(W_{\mathrm{Bl}_4\mathbf{P}^2})=
\left\{-3,\frac{5+5\sqrt{5}}{2},\frac{5-5\sqrt{5}}{2}\right\}.
\]
All three critical points are Morse, so
\[
\overline{\MF((\mathbf{C}^*)^2,W_{\mathrm{Bl}_4\mathbf{P}^2})}
\cong \bigoplus_{i=1}^{3}\MorseCat{2}.
\]
\end{example}

\begin{example}[$\mathrm{Bl}_5\mathbf{P}^2$]
Consider
\[
W_{\mathrm{Bl}_5\mathbf{P}^2}=\frac{(1+x)^2(1+y)^2}{xy}-4.
\]
The critical locus and critical values are
\[
\begin{aligned}
\CritLoc(W_{\mathrm{Bl}_5\mathbf{P}^2})
&=\{(1,1)\}\cup\{x=-1\}\cup\{y=-1\},\\
\CritVal(W_{\mathrm{Bl}_5\mathbf{P}^2})
&=\{-4,12\}.
\end{aligned}
\]
At $(1,1)$ the shifted completed potential has Morse normal form
$u^2+v^2$. We now describe separately the formal models on the two punctured
branches and at their intersection. Put
\[
L_x^\circ=\{x=-1,\ y\neq-1\},
\qquad
L_y^\circ=\{y=-1,\ x\neq-1\},
\qquad
p=(-1,-1).
\]
Since $L_x^\circ$ is closed in the open subset
$U_x=D(y+1)\subset(\mathbf C^*)^2$, setting $u=x+1$ gives
\[
\widehat{(U_x)}_{L_x^\circ}
\simeq
\operatorname{Spf}\!\left(
\mathbf C[y^{\pm1},(y+1)^{-1}][[u]]
\right).
\]
On this completion,
\[
\widehat W+4
=u^2\frac{(1+y)^2}{(u-1)y}
\overset{c}{\sim}u^2,
\]
because the displayed coefficient is a unit. Thus the punctured
$x=-1$ branch has the global formal $A_1$ model over $L_x^\circ$.

Similarly, $L_y^\circ$ is closed in $U_y=D(x+1)$. Setting $v=y+1$ gives
\[
\widehat{(U_y)}_{L_y^\circ}
\simeq
\operatorname{Spf}\!\left(
\mathbf C[x^{\pm1},(x+1)^{-1}][[v]]
\right)
\]
and
\[
\widehat W+4
=v^2\frac{(1+x)^2}{x(v-1)}
\overset{c}{\sim}v^2.
\]
Hence the punctured $y=-1$ branch also has the global formal $A_1$ model.

At the intersection point, set $u=x+1$ and $v=y+1$. Then
\[
\widehat X_p\simeq\operatorname{Spf}\mathbf C[[u,v]],
\qquad
\widehat W+4
=\frac{u^2v^2}{(u-1)(v-1)}
\overset{c}{\sim}u^2v^2.
\]
Thus the three local formal models are $u^2$ on $L_x^\circ$, $v^2$ on
$L_y^\circ$, and $u^2v^2$ at $p$.

For the crossing model, the singular locus of
$\operatorname{Spec}\mathbf C[[u,v]]/(u^2v^2)$ is, set-theoretically, the
union $V(uv)$ of the two coordinate axes. It is not contained in the closed
point $V(u,v)$. Hence Theorem~\ref{thm:general_formal} does not identify the
point-supported contribution with the unrestricted formal category of
$u^2v^2$; the contribution at $p$ must remain a torsion category.

These local descriptions do not give a direct-sum decomposition, since
$L_x^\circ$ and $L_y^\circ$ are not closed in the original torus and their
closures meet at $p$. They do, however, give a localization decomposition.
Put
\[
S_5=\{x=-1\}\cup\{y=-1\},
\qquad
\mathcal C_5=
\overline{\MF^t_{-4}(\widehat X_{S_5},\widehat W)},
\]
and let
\[
\mathcal P_5=
\overline{\MF^t_0\bigl(\operatorname{Spf}\mathbf C[[u,v]],u^2v^2\bigr)}
\]
be the contribution supported at the intersection point. Removing $p$ from
$S_5$ leaves the disjoint union
$L_x^\circ\amalg L_y^\circ$. Supported localization and the two global
$A_1$ normal forms above therefore give a Verdier localization sequence
\begin{equation}\label{eq:bl5_localization}
\mathcal P_5
\longrightarrow
\mathcal C_5
\longrightarrow
\Aper{L_x^\circ}{1}\oplus\Aper{L_y^\circ}{1}.
\end{equation}
Equivalently, after idempotent completion,
\[
\mathcal C_5/\mathcal P_5
\simeq
\Aper{L_x^\circ}{1}\oplus\Aper{L_y^\circ}{1}.
\]
Thus $\mathcal C_5$ is obtained by gluing the two punctured $A_1$
contributions to the formal category of the crossing $u^2v^2$; in general the
localization sequence does not split. Including the isolated Morse point at
critical value $12$ gives
\[
\overline{\MF((\mathbf C^*)^2,W_{\mathrm{Bl}_5\mathbf P^2})}
\simeq
\mathcal C_5\oplus\MorseCat{2}.
\]
\end{example}

For the remaining three models, put $u=1+x+y$ and
\[
Y=\{u=0\}\subset(\mathbf C^*)^2.
\]
Since $x=u-1-y$ and both $y$ and $1+y$ are invertible on $Y$, the completion
of $(\mathbf C^*)^2$ along $Y$ is globally trivial:
\[
\widehat X_Y
\simeq
\operatorname{Spf}\!\left(
\mathbf C[y^{\pm1},(1+y)^{-1}][[u]]
\right)
=\operatorname{Spf}(\mathcal O_Y[[u]]).
\]
Thus $u$ is a global formal normal coordinate. In each case below, the
coefficient multiplying the relevant power of $u$ is a unit on the completed
neighborhood, so relative formal contact equivalence and
Proposition~\ref{prop:Ade_category_model} apply globally.

\begin{example}[$\mathrm{Bl}_6\mathbf{P}^2$]
Consider
\[
W_{\mathrm{Bl}_6\mathbf{P}^2}=\frac{(1+x+y)^3}{xy}-6.
\]
The critical locus and critical values are
\[
\begin{gathered}
\CritLoc(W_{\mathrm{Bl}_6\mathbf{P}^2})=
\{x+y+1=0\}\cup\{(1,1)\},\\
\CritVal(W_{\mathrm{Bl}_6\mathbf{P}^2})=\{-6,21\}.
\end{gathered}
\]
At $(1,1)$ the completed potential has Morse normal form after subtracting
$21$.  Along the curve $Y=\{x+y+1=0\}$ one has
\[
\widehat W+6\overset{c}{\sim}u^3.
\]
Hence Proposition~\ref{prop:Ade_category_model} gives the residual $A_2$
contribution over $Y$, and
\[
\overline{\MF((\mathbf{C}^*)^2,W_{\mathrm{Bl}_6\mathbf{P}^2})}
\cong \Aper{Y}{2}\oplus \MorseCat{2}.
\]
\end{example}

\begin{example}[$\mathrm{Bl}_7\mathbf{P}^2$]
Consider
\[
W_{\mathrm{Bl}_7\mathbf{P}^2}=\frac{(1+x+y)^4}{xy}-12.
\]
The critical locus and critical values are
\[
\begin{aligned}
\CritLoc(W_{\mathrm{Bl}_7\mathbf{P}^2})
&=\{x+y+1=0\}\cup
\left\{\left(\frac12,\frac12\right)\right\},\\
\CritVal(W_{\mathrm{Bl}_7\mathbf{P}^2})
&=\{-12,52\}.
\end{aligned}
\]
The isolated point is Morse after subtracting the critical value $52$, while
along $Y=\{x+y+1=0\}$ one has
\[
\widehat W+12\overset{c}{\sim}u^4.
\]
Thus the curve contributes the residual $A_3$ category, and
\[
\overline{\MF((\mathbf{C}^*)^2,W_{\mathrm{Bl}_7\mathbf{P}^2})}
\cong \Aper{Y}{3}\oplus \MorseCat{2}.
\]
\end{example}

\begin{example}[$\mathrm{Bl}_8\mathbf{P}^2$]
Consider
\[
W_{\mathrm{Bl}_8\mathbf{P}^2}=\frac{(1+x+y)^6}{xy^2}-60.
\]
The critical locus and critical values are
\[
\begin{aligned}
\CritLoc(W_{\mathrm{Bl}_8\mathbf{P}^2})
&=\{x+y+1=0\}\cup
\left\{\left(\frac13,\frac23\right)\right\},\\
\CritVal(W_{\mathrm{Bl}_8\mathbf{P}^2})
&=\{-60,372\}.
\end{aligned}
\]
The isolated point is Morse after subtracting the critical value $372$, while
along $Y=\{x+y+1=0\}$ one has
\[
\widehat W+60\overset{c}{\sim}u^6.
\]
Thus the curve contributes the residual $A_5$ category, and
\[
\overline{\MF((\mathbf{C}^*)^2,W_{\mathrm{Bl}_8\mathbf{P}^2})}
\cong \Aper{Y}{5}\oplus \MorseCat{2}.
\]
\end{example}

The computations are summarized in Table~\ref{tab:dp_decompositions}.  In the
first column, $\mathrm{Bl}_k\mathbf{P}^2$ denotes the blowup of
$\mathbf{P}^2$ at $k$ general points.

\begingroup
\scriptsize
\sloppy
\begin{longtblr}[
  caption = {Decompositions of the derived categories of matrix factorizations for the Landau--Ginzburg models of del Pezzo surfaces.},
  label = {tab:dp_decompositions},
]{
  width = \textwidth,
  colspec = {Q[l,wd=0.11\textwidth]Q[c,wd=0.08\textwidth]Q[c,wd=0.12\textwidth]Q[l,wd=0.20\textwidth]Q[l,wd=0.33\textwidth]},
  colsep = 2pt,
  rowhead = 1,
  row{1} = {font=\bfseries},
  rowsep = 2pt,
  hline{1,2,Z} = {0.4pt},
}
Surface & Count & $\dim H^{\mathrm{even}}$ & Critical support & Decomposition of $\overline{\MF((\mathbf{C}^*)^2,W)}$ \\
$\mathbf{P}^2$ & 3 & 3 & three isolated Morse points & $\bigoplus_{i=1}^{3}\MorseCat{2}$ \\
$\mathbf{P}^1\times\mathbf{P}^1$ & 4 & 4 & four isolated Morse points & $\bigoplus_{i=1}^{4}\MorseCat{2}$ \\
$\mathrm{Bl}_1\mathbf{P}^2$ & 4 & 4 & four isolated Morse points & $\bigoplus_{i=1}^{4}\MorseCat{2}$ \\
$\mathrm{Bl}_2\mathbf{P}^2$ & 5 & 5 & five isolated Morse points & $\bigoplus_{i=1}^{5}\MorseCat{2}$ \\
$\mathrm{Bl}_3\mathbf{P}^2$ & 6 & 6 & six isolated Morse points & $\bigoplus_{i=1}^{6}\MorseCat{2}$ \\
$\mathrm{Bl}_4\mathbf{P}^2$ & 3 & 7 & three isolated Morse points & $\bigoplus_{i=1}^{3}\MorseCat{2}$ \\
$\mathrm{Bl}_5\mathbf{P}^2$ & 2 & 8 & intersecting critical branches $S_5=\{x=-1\}\cup\{y=-1\}$ and one isolated Morse point & $\mathcal C_5\oplus\MorseCat{2}$ \\
$\mathrm{Bl}_6\mathbf{P}^2$ & 2 & 9 & curve $Y$ with residual $A_2$ type and one isolated Morse point & $\Aper{Y}{2}\oplus \MorseCat{2}$ \\
$\mathrm{Bl}_7\mathbf{P}^2$ & 2 & 10 & curve $Y$ with residual $A_3$ type and one isolated Morse point & $\Aper{Y}{3}\oplus \MorseCat{2}$ \\
$\mathrm{Bl}_8\mathbf{P}^2$ & 2 & 11 & curve $Y$ with residual $A_5$ type and one isolated Morse point & $\Aper{Y}{5}\oplus \MorseCat{2}$ \\
\end{longtblr}
\endgroup

For the five toric del Pezzo surfaces in the first five rows, the chosen torus
potentials have only isolated Morse points and their number equals
$\dim H^{\mathrm{even}}(F)$. Beginning with
$\mathrm{Bl}_4\mathbf P^2$, the selected Laurent representatives no longer
have this numerical property on the open torus: the degree-five surface has
only three visible Morse points, while the lower-degree models develop
critical loci with positive-dimensional or intersecting components. Thus the table also
marks the point at which a single torus chart ceases to encode the expected
cohomological size merely by counting isolated critical points; boundary data,
another Laurent presentation, or the full formal category is then essential.

For $\mathrm{Bl}_5\mathbf P^2$, the category $\mathcal C_5$ is determined by
the localization sequence~\eqref{eq:bl5_localization}; it is a gluing of the
two punctured $A_1$ contributions with the formal crossing category
$\mathcal P_5$, not an asserted direct sum of these three terms. In the last
three rows,
\[
Y=\operatorname{Spec}\left(\mathbf{C}[x^{\pm1},y^{\pm1}]/(x+y+1)\right).
\]
Finally, $\Aper{Y}{m}$ denotes the genuine ungraded two-periodic residual
$A_m$ contribution from Definition~\ref{def:periodic_A_type}; it is not the
ordinary derived category of the $A_m$ quiver.
Here \emph{Count} means the number of top-level open-torus contributions in
the displayed decomposition, with each connected gluing support counted once.
The cohomological column is a completeness benchmark: for a del Pezzo surface
it equals $\rho+2=12-d$, the isolated-node count of a general compactified
model. Equality is not asserted for the specified anticanonical open-torus
model, where singularities may collide or boundary contributions may be
absent from the torus \cite{LuntsPrzyjalkowski2018,Przyjalkowski2026RankOne}.
For the last three rows, the global coordinate $u=1+x+y$ trivializes the
formal normal line, and the displayed Laurent-polynomial formulas express the
shifted potential as a unit times $u^{m+1}$. Thus these rows use a global split
formal line model, not merely a pointwise transverse classification.

\section{Open-Torus Decompositions for Fano Threefold Laurent-Polynomial Models}
\label{sec:fano_threefold_decompositions}
We collect the Fano-threefold critical-locus computations in a single table.
All computations behind the table are available in the public repository
\cite{FanoMirrorsComputations}. It contains the companion appendix
\emph{Detailed Computations for the Fano Threefold Models}, with the
model-by-model calculations, and a supplement of certificates, generated
independently, which records the saturated critical ideals, component
equations, Hessian checks, and formal expansions used to audit the
Mori--Mukai entries.
For torus coordinates $x_1,x_2,x_3$, write
$\theta_i=x_i\partial/\partial x_i$ for the logarithmic derivations.
An isolated nondegenerate critical point contributes $\MorseCat{3}$.  A smooth
quadratic component contributes $\MBcat{Z}$ when its Clifford datum is trivial
and $\MBcat{Z;[\delta]}$ when the discriminant square class
$[\delta]\in\mathcal O_Z^*/\mathcal O_Z^{*2}$ is nontrivial.  A residual
transverse normal form $u^{m+1}$ contributes $\Aper{Z}{m}$, while the twisted
residual $D_4$ curve in the $V_{14}$ model is denoted by
$\Dper{Z}{4}{\delta}$. Mixed forms $v^2+u^{m+1}$ and
$v^2+\delta u^{m+1}$ are denoted by $\QAper{Z}{m}$ and
$\QAtwist{Z}{m}{\delta}$, respectively. The exceptional $B_1$, $V_{10}$,
$V_8$, $V_6$, $MM_{2-1}$, and $MM_{2-3}$ models have intersecting, singular,
or coefficient-degenerate positive-dimensional critical loci. In these cases the table records the
category of the formal completion along the entire connected locus, rather
than a direct sum of categories attached to its smooth strata.

The following notation table explains the curves, surfaces, and residual
categories that appear in the decomposition table.

\begingroup
\footnotesize
\begin{longtblr}{
  width = \textwidth,
  colspec = {Q[l,wd=0.18\textwidth]Q[l,wd=0.70\textwidth]},
  rowhead = 1,
  row{1} = {font=\bfseries},
  rowsep = 2pt,
  hline{1,2,Z} = {0.4pt},
}
Symbol & Meaning \\
$\MorseCat{3}$ & The two-periodic Clifford-module category of one isolated Morse critical point in a threefold. \\
$\Gamma$, $\Gamma_i$, $\Gamma_{\pm1}$, $\Gamma_\zeta$ & Smooth curve components of the critical locus in $(\mathbf{C}^*)^3$; the subscript distinguishes the branch, sign, or root-of-unity parameter. \\
$\mathcal C_{B_1,0}$ & The connected contribution of the $B_1$ zero fiber, determined by the nonsplit localization triangle $\mathcal P_{B_1}\to\mathcal C_{B_1,0}\to\MBcat{S^\circ;[A]}$; the category $\mathcal P_{B_1}$ is supported on the $uv^2$ crossing curves. \\
$\mathcal C_{V_{10},-6}$ & The connected contribution of the $V_{10}$ fiber at value $-6$, fitted into the localization triangle $\mathcal P_{V_{10}}\to\mathcal C_{V_{10},-6}\to\Aper{S_{V_{10}}^\circ}{1}$; $\mathcal P_{V_{10}}$ is supported on the $u^2v$ crossing curve. \\
$\mathcal C_{V_8,-8}$ & The connected contribution of the singular $V_8$ fiber at value $-8$, fitted into the localization triangle $\mathcal P_{V_8}\to\mathcal C_{V_8,-8}\to\Aper{S_{V_8}^\circ}{1}$; $\mathcal P_{V_8}$ is supported on the twisted crossing curve. \\
$\mathcal C_{V_6,-12}$ & The connected contribution of the $V_6$ fiber at value $-12$, fitted into the localization triangle $\mathcal P_{V_6}\to\mathcal C_{V_6,-12}\to\Aper{S_A^\circ}{1}\oplus\Aper{S_U^\circ}{2}$; $\mathcal P_{V_6}$ is supported on the $a^2u^3$ crossing curve. \\
$\mathcal C_{MM_{2-1},-60}$, $\mathcal C_{MM_{2-3},-12}$ & The categories of the formal neighborhoods of the whole critical loci, crossing strata included, in the indicated fibers; their supported-localization sequences are not asserted to split. \\
$S$, $\Sigma$ & Model-dependent smooth surface components; their displayed contribution is Morse--Bott or residual $A$-type according to the corresponding row. \\
$\MBcat{Z}$ & The untwisted two-periodic Morse--Bott category of a smooth critical component $Z$. \\
$\MBcat{Z;[\delta]}$ & A Morse--Bott component whose quadratic normal bundle has discriminant square class $[\delta]\in\mathcal O_Z^*/\mathcal O_Z^{*2}$; the class is recorded because the associated Clifford category need not be globally untwisted. \\
$\Aper{Z}{m}$ & The ungraded two-periodic residual $A_m$ contribution associated with a transverse formal normal form $u^{m+1}$. \\
$\Dper{Z}{4}{\delta}$ & The periodic category of the residual $D_4$ normal form over $Z$, with twisting square class $[\delta]$. \\
$\QAper{Z}{m}$ & The standard mixed contribution of the global formal normal form $v^2+u^{m+1}$. \\
$\QAtwist{Z}{m}{\delta}$ & The mixed formal normal form $v^2+\delta u^{m+1}$ with nontrivial square class $[\delta]$. \\
$p(\lambda)=0$, $q(\alpha)=0$, $\Phi(\lambda)=0$, $\Psi(\lambda)=0$ & Square-free elimination equations whose roots index finite collections of isolated Morse points. \\
$\varepsilon\in\{\pm1\}^k$ & A sign vector of length $k$; each of the $2^k$ sign choices indexes one isolated Morse point. \\
\end{longtblr}
\endgroup

The companion appendix records the polynomial identities, saturated critical
ideals, critical points, parametrized components, Jacobian minors, Hessian
computations, and formal coordinate changes used in the table.

\paragraph{Hori--Vafa polynomials for four toric families.}
For a smooth toric Fano threefold with fan rays $v_\rho$, the Hori--Vafa
polynomial is $W^{\mathrm{HV}}=\sum_\rho x^{v_\rho}$. For fourteen of the
eighteen toric families, the polynomial of
\cite{CoatesCortiGalkinKasprzyk2016} has $\dim H^{\mathrm{even}}(F)$
nondegenerate critical points on the torus, with the same number of points
in each fibre as $W^{\mathrm{HV}}$. For the remaining four it has fewer, and
for these families the table uses $W^{\mathrm{HV}}$ instead:
\[
\begin{aligned}
W^{\mathrm{HV}}_{3-31}&=x+y+z+\frac zx+\frac zy+\frac1z,\\
W^{\mathrm{HV}}_{4-10}&=x+y+\frac1x+xy+\frac1{xy}+z+\frac1z
\qquad(\mathbf P^1\times S_7),\\
W^{\mathrm{HV}}_{5-2}&=x+y+z+xy+xyz+\frac1z+\frac1{xy}+\frac1{xyz},\\
W^{\mathrm{HV}}_{5-3}&=x+y+xy+\frac1x+\frac1y+\frac1{xy}+z+\frac1z
\qquad(\mathbf P^1\times S_6).
\end{aligned}
\]
The period sequence of each agrees with that of the polynomial of
\cite{CoatesCortiGalkinKasprzyk2016} up to order $14$ (computed exactly over
$\mathbf Z$). Each $W^{\mathrm{HV}}$ has $\dim H^{\mathrm{even}}(F)$
nondegenerate critical points on the torus, namely $8$, $10$, $12$ and $12$.
The numbers of points in the singular fibres are $2,2,1,1,1,1$ for
$MM_{3-31}$; $2,2,1,1,1,1,1,1$ for $MM_{4-10}$; $2,2,1,\dots,1$ (ten fibres)
for $MM_{5-2}$; and $3,3,2,2,1,1$ for $MM_{5-3}$. These counts are certified
exactly in the repository \cite{FanoMirrorsComputations}: the characteristic
polynomial of multiplication by $3x+5y+7z$ on the torus Jacobian ring is
squarefree, so all critical points are nondegenerate, and the number of
points in a fibre is the multiplicity of the corresponding root of the
characteristic polynomial of multiplication by $W^{\mathrm{HV}}$. In
Table~\ref{tab:fano_threefold_decompositions} the polynomial of $MM_{4-10}$
appears in the row labelled $MM_{4-11}$. In rank four the labels of the
table are shifted with respect to the Mori--Mukai numbering: the polynomial
listed under $MM_{4-(k+1)}$ belongs to $MM_{4-k}$ for $9\le k\le12$, as its
anticanonical degree and its period show.

The symbols $p$, $q$, $r$, $\Phi$, and so forth in the table are row-specific
elimination polynomials, not globally fixed polynomials; their explicit forms
appear under the corresponding label in the companion appendix.  Multiple sums in
one row refer to disjoint algebraic branches of the critical equations, and a
sum over $\{\pm1\}^k$ contains one summand for every sign vector.  In
particular, the rows containing $\MBcat{Z}$ record smooth quadratic curve or
surface components, the notation $\MBcat{Z;[\delta]}$ records the nontrivial
Clifford square class, $V_{14}$ contains a twisted residual $D_4$
contribution, and $V_4$, $V_2$, $MM_{8-1}$, $MM_{9-1}$, and $MM_{10-1}$
contain residual periodic $A$-type contributions. For $MM_{8-1}$
the two curve components contribute
$\QAper{\Gamma_1}{2}$ and $\QAper{\Gamma_{-1}}{2}$.  For $MM_{9-1}$ and
$MM_{10-1}$ the even residual exponents leave nontrivial coefficient classes;
the corresponding terms are
$\QAtwist{\Gamma_\varepsilon}{3}{t(t+1)}$ and
$\QAtwist{\Gamma_\varepsilon}{5}{t+1}$, respectively.  The detailed formal
normalizations are given below.

An independent audit script in the repository
\cite{FanoMirrorsComputations} checks the smooth quadratic
components in the main text.  For each such component it verifies the
critical equations, smoothness, equality between the Hessian rank and the
codimension, and a unit Jacobian minor for the chosen global normal
coordinates.  The discriminant calculation then decides whether the notation
is untwisted or carries a square class.  The nontrivial classes occurring
in the final table are summarized below; each displayed coefficient is a unit
on the corresponding component.

\begin{center}
\footnotesize
\begin{tblr}{
  colspec = {Q[l,wd=0.20\textwidth]Q[l,wd=0.30\textwidth]Q[l,wd=0.35\textwidth]},
  row{1} = {font=\bfseries},
  rowsep = 2pt,
  hline{1,2,Z} = {0.4pt},
}
Model & Component & Quadratic or residual contribution \\
$B_1$ & $S^\circ$ & $\MBcat{S^\circ;[A]}$, where $A=x+y$. \\
$B_2$ & $S$ & $\MBcat{S;[x]}$. \\
$V_{14}$ & $\Gamma$ & $\Dper{\Gamma}{4}{t(t+1)}$. \\
$MM_{2-6}$ & $\Gamma_1$ & $\MBcat{\Gamma_1;[t]}$; the second curve $\Gamma_2$ is untwisted. \\
$MM_{2-12}$ & $\Gamma$ & $\MBcat{\Gamma;[t]}$. \\
$MM_{3-1}$ & $\Gamma_3$ & $\MBcat{\Gamma_3;[v]}$; the curves $\Gamma_\zeta$ are hyperbolic. \\
$MM_{4-1}$ & $\Gamma$ & $\MBcat{\Gamma;[t(t+1)]}$. \\
$MM_{4-3}$ & $\Gamma$ & $\MBcat{\Gamma;[1-t]}$. \\
$MM_{7-1}$ & $\Gamma_{\pm1}$ & $\MBcat{\Gamma_{\pm1};[t]}$. \\
$MM_{8-1}$ & $\Gamma_{\pm1}$ & $\QAper{\Gamma_1}{2}\oplus\QAper{\Gamma_{-1}}{2}$; both mixed models are untwisted. \\
$MM_{9-1}$ & $\Gamma_{\pm1}$ & $\QAtwist{\Gamma_{\pm1}}{3}{t(t+1)}$. \\
$MM_{10-1}$ & $\Gamma_{\pm1}$ & $\QAtwist{\Gamma_{\pm1}}{5}{t+1}$. \\
\end{tblr}
\end{center}

Together with Lemma~\ref{lem:affine_morse_bott_criterion}, these checks
justify every Morse--Bott entry. The same audit detects that the surface in
the $V_8$ model is singular along a critical curve and that the quadratic
coefficient on the smooth $V_{10}$ critical surface vanishes along a curve.
For both models the table therefore lists the single category of the formal
neighborhood of the whole critical surface.

The first column gives the Fano threefold: $Q\subset\mathbf{P}^4$ is the smooth
quadric, $B_i$ and $V_d$ are the Picard-rank-one Iskovskikh labels, and
$MM_{\rho-k}$ denotes the Mori--Mukai family of Picard rank $\rho$.  In this
notation there is no $V_{20}$: its degree would correspond to genus $11$,
which does not occur in the classification, so the list jumps from $V_{18}$
to $V_{22}$ \cite{Iskovskikh1977FanoI,Iskovskikh1978FanoII}.  The column
\emph{Count} records the number of top-level open-torus contributions in the
displayed decomposition, with each connected gluing support counted once; it
should not be interpreted as the number of exceptional objects or ordinary
derived-category summands. The column $\dim H^{\mathrm{even}}=2+2\rho$ is a
cohomological completeness benchmark, not an assertion that every displayed
open-torus model captures the full compactified mirror. For $B_1$, $V_{10}$,
$V_8$, $V_6$, $MM_{2-1}$, and $MM_{2-3}$, the positive-dimensional part of
the critical locus contributes one category, that of its formal
neighborhood; it is not decomposed into its smooth strata, and it is counted
once in the column \emph{Count}. Direct sums indexed by roots refer to the corresponding finite
collections of isolated Morse points.

\withouttoc{\subsection*{Hodge, reducedness, and Milnor audit}}

Among the $46$ families with $h^{1,2}>0$, the eleven displayed models having a
nonreduced singular fiber are
\[
\begin{gathered}
B_2,\ B_1,\ V_{10},\ V_8,\ V_6,\ V_4,\ V_2;\\
MM_{2-1},\ MM_{2-2},\ MM_{2-3},\ MM_{2-8}.
\end{gathered}
\]
The computations in the repository \cite{FanoMirrorsComputations} exhibit
the repeated factors explicitly. The
$B_2,B_1,V_{10},V_8,MM_{2-2},MM_{2-3}$, and $MM_{2-8}$ fibers contain a
square factor; the $V_6$ fiber has the form $xA^2U^3$, with two critical
surfaces meeting along a curve; the $V_4$ and $V_2$ fibers have the forms
$xU^4$ and $xU^6$; and the relevant $MM_{2-1}$ fiber contains a cubic factor.

The other $35$ models with $h^{1,2}>0$ have reduced singular fibers on the
displayed torus. Fifteen of them still have a positive-dimensional critical
locus:
\[
\begin{gathered}
B_4,\ V_{16},\ V_{14},\ V_{12};\\
MM_{2-5},MM_{2-6},MM_{2-7},MM_{2-12},MM_{2-13};\\
MM_{2-15},MM_{2-19},MM_{2-28};\\
MM_{3-1},MM_{3-2},MM_{4-1}.
\end{gathered}
\]
Their contributions are the curve-type and connected gluing categories shown
in Table~\ref{tab:fano_threefold_decompositions}. The remaining $20$ reduced
models have only isolated Morse critical points on the torus. This trichotomy
is verified by the Hodge and Milnor audit script in the repository
\cite{FanoMirrorsComputations}.

The preceding distinction can be made model by model. Among the $59$
families with $h^{1,2}=0$, the five models with non-isolated critical locus
are
\[
MM_{4-3},\quad MM_{7-1},\quad MM_{8-1},\quad MM_{9-1},\quad MM_{10-1}.
\]
Their Picard ranks are respectively $4,7,8,9,10$. The remaining $54$ models
have only isolated critical points, and the Hessian calculations show that
every one of these points is Morse. Thus no displayed model with
$h^{1,2}=0$ has an isolated non-Morse singularity. For the $54$ models with
finite critical locus, the total Milnor number is simply the number of torus
critical points.

For $34$ of these $54$ models one has
$\mu_{\mathrm{tor}}(W)=e(F)=2+2\rho(F)$; they include all eighteen toric
families, four of them with the Hori--Vafa polynomial of
Section~\ref{sec:fano_threefold_decompositions}. The $20$ models with a strict
inequality are
\[
\begin{gathered}
Q,\ V_{22};\\
MM_{2-21},MM_{2-22},MM_{2-29},MM_{2-30};\\
MM_{3-10},MM_{3-13},MM_{3-15},MM_{3-17},MM_{3-18};\\
MM_{3-19},MM_{3-20},MM_{3-22};\\
MM_{4-4},MM_{4-5},MM_{4-7},MM_{4-9};\\
MM_{5-1},MM_{6-1}.
\end{gathered}
\]
Their exact values of $\mu_{\mathrm{tor}}(W)$ and $e(F)$ are the entries in
the columns \emph{Count} and $\dim H^{\mathrm{even}}$ below. The same
script reproduces this classification from the table and checks all $105$
rows. A strict inequality is a statement about the displayed open-torus
polynomial, not an invariant of the Fano family: missing critical points may
occur on the boundary or in a different Laurent presentation.

For Picard-rank-one threefolds, Przyjalkowski computes the isolated
singularities of the fibers of a Calabi--Yau compactification of the standard
toric Landau--Ginzburg model and proves that they are ordinary double points
\cite[Theorem~1.5]{Przyjalkowski2026RankOne}. The following comparison counts
only the isolated $\MorseCat{3}$ summands in our open-torus decomposition;
positive-dimensional and connected gluing contributions are deliberately
excluded from the second column.

\begin{table}[htbp]
\centering
\small
\begin{tblr}{
  colspec = {Q[l,wd=0.25\textwidth]Q[c,wd=0.15\textwidth]Q[c,wd=0.21\textwidth]Q[c,wd=0.15\textwidth]},
  colsep = 2pt,
  row{1} = {font=\bfseries},
  rowsep = 2pt,
  hline{1,2,Z} = {0.4pt},
}
Threefolds & Open isolated Morse points & Compactified $\varepsilon(X)$ & Difference \\
$V_2,V_4,V_6,V_8$ & 1 & 1 & 0 \\
$V_{10},V_{12},V_{14},V_{16},V_{18}$ & 2 & 2 & 0 \\
$B_1,B_2,B_3,B_4$ & 2 & 2 & 0 \\
$B_5,\mathbf P^3$ & 4 & 4 & 0 \\
$V_{22}$ & 3 & 4 & 1 \\
$Q$ & 3 & 4 & 1 \\
\end{tblr}
\caption{Isolated open-torus Morse points compared with the isolated ordinary double points in the compactified Picard-rank-one models of \cite{Przyjalkowski2026RankOne}. Here $B_i$ corresponds to the notation $Y_i$ used there.}
\label{tab:rank_one_compactified_comparison}
\end{table}

Thus the isolated counts agree for $15$ of the $17$ Picard-rank-one families.
For $V_{22}$ and $Q$, the chosen compactification has one
additional isolated ordinary double point beyond the open torus. This
comparison locates a numerical boundary deficit for these specified models;
it does not by itself compute the full boundary matrix-factorization category
or all vanishing-cycle data.

For the $B_2$ model, put
\[
F_{B_2}=x(y+z+1)+\frac{1}{xyz},
\qquad
W_{B_2}=\frac{F_{B_2}^2}{x}.
\]
The surface $S=V(F_{B_2})$ is smooth in the torus, and its normal Hessian has
rank one everywhere. Since $W_{B_2}=F_{B_2}^2/x$, its quadratic coefficient
has square class $[x]$. This class is nontrivial: on the curve
$\{(t,-1,t^{-1})\}\subset S$, it restricts to $[t]$. The complementary
critical points are
$(1,1,1)$ and $(-1,1,1)$, with ordinary Hessian determinants $512$ and
$-512$, respectively. Thus the surface contributes $\MBcat{S;[x]}$.

For the source $B_4$ model, the positive-dimensional part of the critical
locus is a single smooth curve, disjoint from its two isolated critical points. Thus its
contribution is an ordinary direct sum rather than a connected gluing category.

\begin{example}[$B_4$: a smooth Morse--Bott curve]
\label{ex:B4_morse_bott_curve}
The Laurent polynomial listed in Appendix~A of
\cite{CoatesCortiGalkinKasprzyk2016} is
\[
W_{B_4}
=x+yz^2+2yz+y+\frac2{yz}+\frac1{xy^2z^2}.
\]
Its saturated logarithmic critical ideal decomposes as
\[
I_{\mathrm{crit}}
=(xy-1,z+1)\cap(x-y,y^2-1,z-1).
\]
Hence the critical locus is the disjoint union of the smooth curve
\[
\Gamma=\{(t,t^{-1},-1):t\in\mathbf C^*\}
\]
and the points
\[
P_+=(1,1,1),
\qquad
P_-=(-1,-1,1).
\]
The curve lies in the zero fiber. The two points have critical values $8$ and
$-8$ and logarithmic Hessian determinants $32$ and $-32$, respectively, so
both are Morse.

Along $\Gamma$, use the global parameter and normal coordinates
\[
t=x,
\qquad
U=xy-1,
\qquad
V=z+1.
\]
Then $y=(1+U)/t$, $z=V-1$, and
\[
W_{B_4}
\equiv
t(U-V)^2+\frac1tV^2
\pmod{(U,V)^3}.
\]
The quadratic form splits globally as
\[
t\left(U-\left(1-\frac{\mathrm i}{t}\right)V\right)
 \left(U-\left(1+\frac{\mathrm i}{t}\right)V\right),
\]
and the corresponding linear change has determinant $-2\mathrm i/t$, a unit
on $\Gamma$. The relative formal Morse lemma therefore gives an untwisted
Morse--Bott model. Consequently,
\[
\overline{\MF((\mathbf C^*)^3,W_{B_4})}
\simeq
\MBcat{\Gamma}\oplus\MorseCat{3}\oplus\MorseCat{3}.
\]
\end{example}

\begin{example}[$V_{14}$: a twisted residual $D_4$ curve]
\label{ex:V14_relative_D4}
For the $15$-term Laurent polynomial listed in Appendix~A of
\cite{CoatesCortiGalkinKasprzyk2016}, put
\[
A=yz^2+yz+1
\]
and
\[
B=x^2yz^2+xy^2z^4+2xy^2z^3+xy^2z^2+2xyz+x+yz^2.
\]
The source polynomial satisfies the exact identity
\[
W_{V_{14}}+4=\frac{AB}{xy^2z^3}.
\]
On $A=0$, one has
\[
y=-\frac1{z(z+1)},
\qquad
B=-\frac{z(x-1)^2}{z+1}.
\]
Thus the critical fiber at $-4$ contains the smooth curve
\[
\Gamma=\left\{
\left(1,-\frac1{t(t+1)},t\right)
:t\in\mathbf C^*\setminus\{-1\}
\right\}.
\]
There are also two isolated Morse points,
\[
P_+=(1,1,1),
\qquad
P_-=\left(-1,9,-\frac13\right),
\]
with critical values $23$ and $-5$ and logarithmic Hessian determinants
$1458$ and $-2$, respectively.

Let $R=\mathbf C[t^{\pm1},(t+1)^{-1}]$ and use the global normal coordinates
\[
U=A,
\qquad
V=x-1.
\]
After multiplication by a unit, the completed shifted potential is
\[
U\left(U^2-\frac{t}{t+1}V^2\right)
\]
up to terms of order at least four. The relative finite-determinacy argument
for $D_4$ removes these higher terms globally over $R$. The coefficient has
nontrivial square class
\[
\left[\frac{t}{t+1}\right]=[t(t+1)]
\quad\text{in }R^*/R^{*2}.
\]
Consequently,
\[
\overline{\MF((\mathbf C^*)^3,W_{V_{14}})}
\simeq
\Dper{\Gamma}{4}{t(t+1)}
\oplus\MorseCat{3}\oplus\MorseCat{3}.
\]
\end{example}

\begin{example}[$V_{16}$: an untwisted Morse--Bott curve]
Consider the $14$-term Laurent polynomial listed in Appendix~A of
\cite{CoatesCortiGalkinKasprzyk2016}:
\[
\begin{aligned}
W_{V_{16}}={}&x+\frac{2x}{yz}+\frac{x}{y^2z^2}+yz^2+2yz+y+2z
+\frac2z+\frac1y+\frac2{yz}+\frac1{yz^2}\\
&+\frac zx+\frac2x+\frac1{xz}.
\end{aligned}
\]
Its reduced critical locus is the disjoint union of
\[
\Gamma=\{y=1,\ z=-1\}
\]
and the two points
\[
P_a=(4-2a,a,1),\qquad a^2=2.
\]
Along $\Gamma$ the critical value is $-4$.  With $t=x$, $U=y-1$, and
$V=z+1$, the completion is
\[
\widehat X_\Gamma
\simeq
\operatorname{Spf}\mathbf C[t^{\pm1}][[U,V]],
\]
and the transverse quadratic part is
\[
W_{V_{16}}+4
\equiv
tU^2-2tUV+\frac{t^2-1}{t}V^2
\pmod{(U,V)^3}.
\]
The quadratic form splits globally as
\[
t\left(U-\left(1+\frac1t\right)V\right)
 \left(U-\left(1-\frac1t\right)V\right),
\]
and the corresponding linear change has determinant $2/t$. The relative
formal Morse lemma therefore gives an untwisted global Morse--Bott model. At
$P_a$ one has
\[
\det(\theta_i\theta_jW_{V_{16}})(P_a)=16(17a+24),
\]
which is nonzero modulo $a^2-2$.  Consequently,
\[
\overline{\MF((\mathbf C^*)^3,W_{V_{16}})}
\simeq
\MBcat{\Gamma}\oplus\bigoplus_{a^2=2}\MorseCat{3}.
\]
\end{example}

\begin{example}[$V_8$: a singular critical surface]
\label{ex:V8_singular_surface}
For the Laurent-polynomial representative of $V_8$ in
\cite{CoatesCortiGalkinKasprzyk2016}, set
\[
F_{V_8}
=x^2y^2z^2+x^2yz^2(z+1)^2+2xyz+1,
\qquad
W_{V_8}=\frac{F_{V_8}^2}{x^3y^2z^4}-8.
\]
Thus the critical locus in the fiber at value $-8$ is the surface
\[
S_{V_8}=V(F_{V_8}).
\]
This surface is not smooth.  Its singular locus in the torus is the curve
\[
C_{V_8}=V(xy-1,z+1).
\]
Indeed, along this curve use $t=x$ and the normal coordinates
\[
a=xy-1,
\qquad
b=z+1.
\]
Then
\[
F_{V_8}
\equiv a^2-2ab+(t+1)b^2
\pmod{(a,b)^3},
\]
and the determinant of this transverse Hessian is $4t$, a unit on
$C_{V_8}$.  Hence $F_{V_8}$ has a transverse nodal quadratic germ along the
curve, while $W_{V_8}+8$ is the square of that germ up to a unit.  Pointwise,
the resulting transverse contact type is $u^2v^2$, not a Morse--Bott quadratic
form.  Away from $C_{V_8}$ the open stratum
$S_{V_8}\setminus C_{V_8}$ is Morse--Bott, but it is not closed in
$S_{V_8}$ and cannot be split off as a direct summand.

The contribution of $S_{V_8}$ is therefore the category of its whole formal
neighborhood,
\[
\mathcal C_{V_8,-8}:=
\overline{\MF^t_{-8}(\widehat X_{S_{V_8}},\widehat W)}.
\]
More precisely, put
\[
A=1+\frac1{xyz},\qquad B=z+1,
\qquad U=yA^2+B^2.
\]
Then the exact identity
\[
W_{V_8}+8=xU^2
\]
shows that $S_{V_8}=V(U)$.  On the open stratum
$S_{V_8}^\circ=S_{V_8}\setminus C_{V_8}$, the function $U$ is a regular
normal coordinate and gives the residual contribution
$\Aper{S_{V_8}^\circ}{1}$.  Along
\[
C_{V_8}=V(A,B)
=\{(t,t^{-1},-1):t\in\mathbf C^*\},
\]
the pair $(A,B)$ is a global system of formal normal coordinates and the
leading transverse form is
\[
\frac{(A^2+tB^2)^2}{t}.
\]
The square class $[t]$ is nontrivial on $C_{V_8}$, so this crossing family is
not globally equivalent to an untwisted $u^2v^2$ model.  If
$\mathcal P_{V_8}$ denotes its torsion formal category, supported localization
gives the nonsplit Verdier sequence
\[
\mathcal P_{V_8}
\longrightarrow
\mathcal C_{V_8,-8}
\longrightarrow
\Aper{S_{V_8}^\circ}{1}.
\]
On the smooth stratum the formal zero fiber is defined by $U^2=0$, and hence
its singular locus is contained in the zero section $V(U)$.  Thus the formal
comparison theorem applies there; no analogous untwisted replacement is made
on the crossing curve.
The only critical point away from $S_{V_8}$ is $(1,1,1)$; it has critical
value $56$ and ordinary Hessian determinant $32768$. Consequently,
\[
\overline{\MF((\mathbf C^*)^3,W_{V_8})}
\simeq
\mathcal C_{V_8,-8}\oplus\MorseCat{3}.
\]
\end{example}

\begin{example}[$V_{12}$: an untwisted Morse--Bott curve]
\label{ex:V12_morse_bott_curve}
Consider the $18$-term Laurent polynomial listed in Appendix~A of
\cite{CoatesCortiGalkinKasprzyk2016}:
\[
\begin{aligned}
W_{V_{12}}={}&x^2y^3z+x^2y^2z+2xy^2z+xy^2+2xyz+2xy+x
+yz+3y+z\\
&+\frac2y+\frac1x+\frac1{xz}+\frac2{xy}+\frac3{xyz}
+\frac1{xy^2}+\frac3{xy^2z}+\frac1{xy^3z}.
\end{aligned}
\]
Its reduced critical locus is the disjoint union of the smooth curve
\[
\Gamma=\{y=-1,\ z(x-1)^2+1=0\}
\]
and the two points
\[
P_a=\left(2a-1,\frac{16a-8}{9},a\right),
\qquad 32a^2-32a-1=0.
\]
The curve is parametrized by
\[
(x,y,z)=\left(t,-1,-\frac1{(t-1)^2}\right),
\qquad t\in\mathbf C^*,\quad t\neq1,
\]
and $W_{V_{12}}|_\Gamma=-5$.  Put
\[
U=y+1,\qquad V=z(x-1)^2+1.
\]
The source polynomial is exactly divisible by $U$ after shifting by $5$.
More precisely, if
\[
H:=\frac{xy^3z(W_{V_{12}}+5)}{y+1},
\]
then $H$ is a Laurent polynomial and $H|_{y=-1}=-xzV$. Thus the two
components of the critical fiber meet precisely along $\Gamma$.
The Jacobian minor in the variables $(y,z)$ is $(t-1)^2$, a unit on
$\Gamma$, and hence
\[
\widehat X_\Gamma
\simeq
\operatorname{Spf}
\mathbf C[t^{\pm1},(t-1)^{-1}][[U,V]].
\]
Modulo $(U,V)^3$, one has
\[
W_{V_{12}}+5
\equiv
\frac{t^3-t^2+3t-1}{t(t-1)}U^2+UV.
\]
The determinant of the transverse Hessian is $-1$. In fact,
\[
\widetilde V:=\frac{W_{V_{12}}+5}{U}
\equiv
V+\frac{t^3-t^2+3t-1}{t(t-1)}U
\pmod{(U,V)^2}.
\]
Consequently $(U,\widetilde V)$ is a global formal coordinate system and the
shifted potential is exactly $U\widetilde V$. The quadratic bundle is
globally hyperbolic and untwisted. At the isolated points,
\[
\det(\theta_i\theta_jW_{V_{12}})(P_a)
=\frac8{3}(1154a+35),
\]
which is nonzero modulo $32a^2-32a-1$.  Therefore
\[
\overline{\MF((\mathbf C^*)^3,W_{V_{12}})}
\simeq
\MBcat{\Gamma}
\oplus
\bigoplus_{a:\,32a^2-32a-1=0}\MorseCat{3}.
\]
\end{example}

\begin{example}[$V_{10}$: a critical surface with a crossing curve]
\label{ex:V10_critical_surface}
For the Laurent polynomial listed in Appendix~A of
\cite{CoatesCortiGalkinKasprzyk2016}, put
\[
F=xyz^2+xyz+1,
\qquad
G=xyz^2+xyz+xz+yz+1.
\]
The $21$ terms of the source polynomial group into the exact factorization
\[
W_{V_{10}}+6=\frac{F^2G}{x^2y^2z^3}.
\]
Thus the fiber at $-6$ is nonreduced and contains the smooth critical
surface
\[
S_{V_{10}}=V(F).
\]
The surface is smooth because $\partial F/\partial x=yz(z+1)$ restricts to a
unit on $S_{V_{10}}$. Saturating the logarithmic critical ideal by $F$ gives
two further critical points,
\[
P_a=(4a-2,4a-2,a),
\qquad
4a^2-2a-1=0.
\]
They satisfy
\[
W_{V_{10}}(P_a)=40a+6,
\qquad
\det(\theta_i\theta_jW_{V_{10}})(P_a)=64(123a+38),
\]
and are both Morse.

The coefficient $G$ is a unit on the open stratum
\[
S_{V_{10}}^\circ=S_{V_{10}}\setminus K,
\qquad
K=V(F,G).
\]
There the transverse formal equation is $F^2$, and hence the contribution is
$\Aper{S_{V_{10}}^\circ}{1}$. Since $G-F=z(x+y)$, the exceptional stratum is
the smooth curve
\[
K=\{y=-x,\ x^2z(z+1)=1\}.
\]
Moreover,
\[
\det\!\left(\frac{\partial(F,G)}{\partial(x,y)}\right)\bigg|_K
=-2xz^2(z+1),
\]
which is a unit on $K$. Thus $(F,G)$ are global formal normal coordinates
along this curve, and the shifted potential is contact-equivalent there to
the untwisted crossing form $u^2v$.

Let $\mathcal P_{V_{10}}$ be the torsion formal category of this $u^2v$
crossing family, and define
\[
\mathcal C_{V_{10},-6}
:=
\overline{\MF^t_{-6}(\widehat X_{S_{V_{10}}},\widehat W_{V_{10}})}.
\]
Supported localization gives the nonsplit Verdier sequence
\[
\mathcal P_{V_{10}}
\longrightarrow
\mathcal C_{V_{10},-6}
\longrightarrow
\Aper{S_{V_{10}}^\circ}{1}.
\]
Consequently,
\[
\overline{\MF((\mathbf C^*)^3,W_{V_{10}})}
\simeq
\mathcal C_{V_{10},-6}
\oplus
\bigoplus_{a:\,4a^2-2a-1=0}\MorseCat{3}.
\]
\end{example}

\begin{example}[$MM_{2-7}$: a Morse--Bott curve and three Morse points]
\label{ex:MM2-7_morse_bott}
Set
\[
K=x(y+1)^2,
\qquad
S=1+\frac{1}{xy},
\]
and consider
\[
W_{MM_{2-7}}
=K\left(S^2+\frac1z\right)+zS^3+S-5.
\]
The isolated critical points are
\[
p_\alpha
=\left(\alpha,1,-28\alpha^2+56\alpha-24\right),
\qquad
4\alpha^3-4\alpha^2-4\alpha+3=0.
\]
The cubic has three distinct roots, and the Hessian is nondegenerate at each
$p_\alpha$.

The positive-dimensional component is
\[
\Gamma=
\left\{
\left(
-\frac{t^2+t+1}{t(t+1)^2},
t,
-\frac{(t^2+t+1)^2}{t^2}
\right)
\right\},
\]
where
\[
t(t+1)(t^2+t+1)\neq0.
\]
Along this curve,
\[
W_{MM_{2-7}}|_\Gamma=-5.
\]
We now exhibit a global split formal model.  Define
\[
V=zS^2+1,
\qquad
R=\frac Kz+S.
\]
Writing $D(t)=t^2+t+1$, the restrictions of the basic functions to $\Gamma$
are
\[
S|_\Gamma=-\frac{t}{D(t)},
\qquad
K|_\Gamma=-\frac{D(t)}{t},
\qquad
z|_\Gamma=-\frac{D(t)^2}{t^2}.
\]
It follows immediately that
\[
V|_\Gamma=0,
\qquad
R|_\Gamma=0.
\]
Moreover, the shifted potential factors algebraically:
\begin{align*}
W_{MM_{2-7}}+5
&=K\left(S^2+\frac1z\right)+S(zS^2+1)\\
&=\left(\frac Kz+S\right)(zS^2+1)\\
&=RV.
\end{align*}

It remains to show that $V,R$ are genuine global normal coordinates rather
than merely functions vanishing on the curve.  Use $t=y$ as the coordinate on
$\Gamma$.  A direct calculation gives
\[
\det\!\left(
\frac{\partial(y,V,R)}{\partial(x,y,z)}
\right)\bigg|_\Gamma
=\frac{t^3(t+1)^2}{D(t)^3}.
\]
This is a unit in
\[
\mathcal O(\Gamma)
=\mathbf C[t^{\pm1},(t+1)^{-1},D(t)^{-1}].
\]
The formal inverse function theorem therefore yields
\[
\mathcal O(\Gamma)[[V,R]]
\xrightarrow{\sim}
\widehat{\mathcal O}_{(\mathbf C^*)^3,\Gamma}.
\]
Under this isomorphism the shifted potential is exactly $VR$.  In particular,
the determinant of the transverse quadratic Hessian is
\[
-\left(\frac{t^3(t+1)^2}{D(t)^3}\right)^2
=-\frac{t^6(t+1)^4}{D(t)^6},
\]
recovering the previously computed unit.

Finally, if
\[
u=\frac{V+R}{2},
\qquad
v=\frac{V-R}{2},
\]
then
\[
VR=u^2-v^2.
\]
Thus the normal quadratic bundle is a global hyperbolic plane.  Its Clifford
algebra is super-Morita trivial, so there is no super-Brauer twist and
$\Gamma$ contributes the untwisted Morse--Bott category $\MBcat{\Gamma}$.
Therefore
\[
\overline{\MF((\mathbf C^*)^3,W_{MM_{2-7}})}
\cong
\MBcat{\Gamma}
\oplus
\bigoplus_{\alpha:\,
4\alpha^3-4\alpha^2-4\alpha+3=0}
\MorseCat{3}.
\]
\end{example}

\begin{example}[$MM_{2-13}$: an exact square completion]
\label{ex:MM2-13_square_completion}
Set
\[
A=1+\frac{y}{x},
\qquad
B=1+\frac1y,
\qquad
F=zA+B,
\]
and consider
\[
W_{MM_{2-13}}
=zA^2+2AB+\frac{B^2}{z}+xB+y-2
=\frac{F^2}{z}+xB+y-2.
\]
The branch $F=0$ of the critical equations gives
\[
\Gamma=\{(1,-1,t):t\in\mathbf C^*\},
\qquad
W_{MM_{2-13}}|_\Gamma=-3.
\]
We now verify the global formal structure along this curve.  Put
\[
U=x+y,
\qquad
V=1+\frac1y=\frac{y+1}{y}.
\]
Both functions vanish on $\Gamma$.  In fact, if
\[
I_\Gamma=(x-1,y+1),
\]
then $V\in I_\Gamma$ and
\[
x-1=U-(y+1)=U-yV,
\qquad
y+1=yV.
\]
Hence $I_\Gamma=(U,V)$.  Their differentials restrict to
\[
dU|_\Gamma=dx+dy,
\qquad
dV|_\Gamma=-dy.
\]
Thus $dU,dV$ form a global basis of $I_\Gamma/I_\Gamma^2$.  More explicitly,
using $t=z$ as the coordinate on $\Gamma$, one has
\[
\det\!\left(
\frac{\partial(z,U,V)}{\partial(x,y,z)}
\right)\bigg|_\Gamma=-1.
\]
The formal inverse function theorem therefore gives
\[
\mathcal O_\Gamma[[U,V]]
\xrightarrow{\sim}
\widehat{\mathcal O}_{(\mathbf C^*)^3,\Gamma},
\qquad
\mathcal O_\Gamma=\mathbf C[t^{\pm1}].
\]

The shifted potential admits an exact quadratic presentation in the ideal
$(U,V)$.  Indeed,
\[
xB+y+1
=x\frac{y+1}{y}+y+1
=\frac{(x+y)(y+1)}{y}
=UV,
\]
while
\[
A=\frac{x+y}{x}=\frac{U}{x},
\qquad
B=V,
\qquad
F=\frac{z}{x}U+V.
\]
Consequently,
\begin{align*}
W_{MM_{2-13}}+3
&=\frac1z\left(\frac{z}{x}U+V\right)^2+UV\\
&=\frac{z}{x^2}U^2
  +\left(\frac2x+1\right)UV
  +\frac1zV^2.
\end{align*}
The coefficients still depend on the normal variables through $x$; therefore
the passage to a quadratic form with coefficients in $\mathcal O_\Gamma$
requires the relative formal Morse lemma rather than an additional algebraic
identity.  Restricting the quadratic coefficients to $\Gamma$, where
$x=1$ and $z=t$, gives
\[
q=tU^2+3UV+\frac1tV^2.
\]
Its Hessian determinant is
\[
\det
\begin{pmatrix}
2t&3\\
3&2/t
\end{pmatrix}
=-5,
\]
a unit in $\mathcal O_\Gamma$.  Hence there is a formal automorphism over
$\Gamma$ carrying $\widehat W_{MM_{2-13}}+3$ to $q$.

It remains to determine the Clifford class.  Set
\[
\widetilde U=U+\frac{3}{2t}V.
\]
Then
\[
q=t\widetilde U^2-\frac{5}{4t}V^2
=t\left(\widetilde U-\frac{\sqrt5}{2t}V\right)
  \left(\widetilde U+\frac{\sqrt5}{2t}V\right).
\]
Since $t$ is invertible and $\sqrt5\in\mathbf C$, this is a global hyperbolic
plane over $\Gamma$.  Its Clifford algebra is therefore super-Morita trivial:
there is no nontrivial super-Brauer twist.  Thus $\Gamma$ contributes the
untwisted Morse--Bott category $\MBcat{\Gamma}$.

The complementary branch gives the three isolated points
\[
P_\xi=\left(\xi,\frac{\xi^2}{4},1\right),
\qquad
\xi^3=16.
\]
Their Hessians are nondegenerate. Consequently,
\[
\overline{\MF((\mathbf C^*)^3,W_{MM_{2-13}})}
\cong
\MBcat{\Gamma}
\oplus
\bigoplus_{\xi:\,\xi^3=16}\MorseCat{3}.
\]
\end{example}

\begin{example}[$MM_{3-2}$: an exact split quadratic curve]
\label{ex:MM3-2_split_curve}
Set
\[
U=1+\frac{x+y}{xyz},
\]
and consider the compact potential
\[
W_{MM_{3-2}}
=xyz^2U^3+z(xy+1)U-6.
\]
The positive-dimensional critical component is
\[
\Gamma
=\{xy+1=0,\ z-x-y=0\},
\]
with parametrization
\[
\Gamma=
\left\{
\left(t,-\frac1t,t-\frac1t\right):
t\in\mathbf C^*,\ t^2\neq1
\right\}.
\]
Along $\Gamma$,
\[
W_{MM_{3-2}}|_\Gamma=-6.
\]
Define
\[
r=zU,
\qquad
s=xy+1+xyzU^2.
\]
Then one has the exact identity
\[
W_{MM_{3-2}}+6=rs.
\]
We verify that this identity gives a genuine global formal splitting along
$\Gamma$.  Let
\[
A=\mathbf C[x^{\pm1},y^{\pm1},z^{\pm1}],
\qquad
I_\Gamma=(xy+1,z-x-y)\subset A.
\]
Since
\[
r=zU=z+\frac1x+\frac1y
=(z-x-y)+\frac{x+y}{xy}(xy+1),
\]
one has $r\in I_\Gamma$.  Moreover, $U=r/z$, and hence
\[
s=(xy+1)+\frac{xy}{z}r^2\in I_\Gamma.
\]
Equivalently, on $\Gamma$ one has $xy=-1$, $z=x+y$, and therefore
\[
U|_\Gamma
=1+\frac{x+y}{xyz}
=1+\frac{z}{(-1)z}
=0,
\qquad
r|_\Gamma=s|_\Gamma=0.
\]

We next check that $r$ and $s$ are global normal parameters.  From
\[
r=z+\frac1x+\frac1y
\]
we obtain
\[
dr=-\frac{1}{x^2}\,dx-\frac{1}{y^2}\,dy+dz.
\]
Also,
\[
ds=y\,dx+x\,dy+U^2\,d(xyz)+2xyzU\,dU,
\]
so the terms containing $U$ vanish after restriction to $\Gamma$.  Under the
parametrization
\[
x=t,\qquad y=-\frac1t,\qquad z=t-\frac1t,
\]
the restricted differentials are
\[
dr|_\Gamma=-\frac1{t^2}\,dx-t^2\,dy+dz,
\qquad
ds|_\Gamma=-\frac1t\,dx+t\,dy.
\]
The minor corresponding to the columns $dy,dz$ is
\[
\det
\begin{pmatrix}
-t^2&1\\
t&0
\end{pmatrix}
=-t.
\]
Since $t$ is invertible on $\Gamma$, the classes of $dr$ and $ds$ form a
global basis of the conormal bundle $I_\Gamma/I_\Gamma^2$.

Finally,
\[
\mathcal O(\Gamma)
=\mathbf C[t^{\pm1},(t^2-1)^{-1}],
\]
where the localization by $t^2-1$ records the condition $z\neq0$.  The
coordinate $t=x$, together with the normal functions $r,s$, induces a
continuous homomorphism
\[
\mathcal O(\Gamma)[[r,s]]
\longrightarrow
\widehat A_{I_\Gamma}.
\]
Its Jacobian determinant along $\Gamma$ is
\[
\det\!\left(
\frac{\partial(x,r,s)}{\partial(x,y,z)}
\right)\bigg|_\Gamma=-t,
\]
which is a unit.  The formal inverse function theorem therefore gives an
isomorphism
\[
\mathcal O(\Gamma)[[r,s]]
\xrightarrow{\sim}
\widehat A_{I_\Gamma},
\]
or, sheaf-theoretically,
\[
\mathcal O_\Gamma[[r,s]]
\xrightarrow{\sim}
\widehat{\mathcal O}_{(\mathbf C^*)^3,\Gamma}.
\]
Under this isomorphism the shifted potential is exactly $rs$.  Thus the formal
neighborhood is globally split and the contribution of $\Gamma$ is the
Morse--Bott category $\MBcat{\Gamma}$.

The isolated critical points form two cubic branches:
\[
P^+_\lambda=(1,1,\lambda),
\qquad
\lambda^3+4\lambda^2-4=0,
\]
and
\[
P^-_\mu=(-1,-1,\mu),
\qquad
\mu^3-2\mu^2+4=0.
\]
All six points are Morse. Hence
\[
\overline{\MF((\mathbf C^*)^3,W_{MM_{3-2}})}
\cong
\MBcat{\Gamma}
\oplus
\bigoplus_{\lambda:\,
\lambda^3+4\lambda^2-4=0}\MorseCat{3}
\oplus
\bigoplus_{\mu:\,
\mu^3-2\mu^2+4=0}\MorseCat{3}.
\]
\end{example}

\begin{example}[$MM_{8-1}$: two residual $A_2$ curves]
\label{ex:MM8-1_residual_curves}
Consider
\[
W_{MM_{8-1}}
=x+\frac1x+y(z+1)^3
+3\left(z+\frac1z\right)
+\frac{3(z+1)}{yz^2}
+\frac{1}{y^2z^3}.
\]
Put
\[
A=yz(z+1).
\]
The reduced critical locus consists of the two isolated points
\[
P_\varepsilon=(\varepsilon,1,1),
\qquad
\varepsilon\in\{\pm1\},
\]
and the two curves
\[
\Gamma_\varepsilon
=\{x=\varepsilon,\ yz(z+1)+1=0\}.
\]
The critical values are
\[
W(P_\varepsilon)=21+2\varepsilon,
\qquad
W|_{\Gamma_\varepsilon}=2\varepsilon-6.
\]
The logarithmic Hessian determinant at $P_\varepsilon$ is
\[
\det(\theta_i\theta_jW)(P_\varepsilon)=486\varepsilon,
\]
so both isolated points are Morse.

Let $u=A+1$. Along $\Gamma_\varepsilon$, one has
\[
W_{MM_{8-1}}-(2\varepsilon-6)
=
\frac{(x-\varepsilon)^2}{x}
+\frac{u^3}{y^2z^3}.
\]
The functions $x-\varepsilon$ and $u$ are global normal coordinates: using
$z$ as the coordinate along the curve,
\[
\det\!\left(
\frac{\partial(z,x-\varepsilon,u)}{\partial(x,y,z)}
\right)=z(z+1),
\]
which is a unit on $\Gamma_\varepsilon$.  Put $a=x-\varepsilon$.  In the
completed coordinate ring,
\[
x=\varepsilon+a,
\qquad
y=\frac{u-1}{z(z+1)},
\]
and hence the shifted potential is exactly
\[
G=\frac{a^2}{\varepsilon+a}
+\frac{(z+1)^2}{z(1-u)^2}u^3.
\]
Set $s=z^2(z+1)^2$ and $\lambda=s^2$.  The unit
$\lambda/(\varepsilon+a)$ has a formal square root, since its restriction to
the zero section is the square unit $s^2/\varepsilon$.  Moreover,
\[
\lambda\frac{(z+1)^2}{z(1-u)^2}
=\bigl(z(z+1)^2\bigr)^3(1-u)^{-2}
\]
has a formal cube root.  Therefore the formal rescalings
\[
v=a\left(\frac{\lambda}{\varepsilon+a}\right)^{1/2},
\qquad
w=u\left(
\lambda\frac{(z+1)^2}{z(1-u)^2}
\right)^{1/3}
\]
give the contact identity $\lambda G=v^2+w^3$.  Thus each curve contributes
the standard mixed category $\QAper{\Gamma_\varepsilon}{2}$, and
\[
\overline{\MF((\mathbf C^*)^3,W_{MM_{8-1}})}
\cong
\QAper{\Gamma_1}{2}
\oplus
\QAper{\Gamma_{-1}}{2}
\oplus
\bigoplus_{\varepsilon\in\{\pm1\}}\MorseCat{3}.
\]
\end{example}

\begin{example}[$MM_{9-1}$: two twisted mixed $A_3$ curves]
Let
\[
F=z+1+\frac{1}{xz},
\qquad
W_{MM_{9-1}}=xF^4+y+\frac1y-12.
\]
The reduced critical locus consists of the two curves
\[
\Gamma_\varepsilon=
\left\{
y=\varepsilon,\quad x=-\frac{1}{z(z+1)}
\right\},
\qquad
\varepsilon\in\{\pm1\},\quad z\neq-1,
\]
and the isolated points
\[
P_\varepsilon=\left(4,\varepsilon,\frac12\right).
\]
The curve has critical value $2\varepsilon-12$, while
\[
\det(\theta_i\theta_jW)(P_\varepsilon)=4096\varepsilon,
\]
so the isolated points are Morse.

Put $D=y-\varepsilon$ and write $t=z$ on $\Gamma_\varepsilon$.  There is an
exact identity
\[
W_{MM_{9-1}}-(2\varepsilon-12)
=xF^4+\frac{D^2}{y}.
\]
The functions $D$ and $F$ are global normal coordinates, because
\[
\det\!\left(
\frac{\partial(z,D,F)}{\partial(x,y,z)}
\right)\bigg|_{\Gamma_\varepsilon}
=t(t+1)^2
\]
is a unit.  In the completion, put
\[
x_0=-\frac{1}{t(t+1)},
\qquad
v=D(\varepsilon+D)^{-1/2},
\qquad
u=F\left(\frac{x}{x_0}\right)^{1/4}.
\]
The indicated roots exist because their arguments restrict to nonzero
constants or to $1$ on the zero section.  The shifted potential becomes
\[
v^2+x_0u^4.
\]
If $h=\mathrm i\,t(t+1)$, then $h^2x_0=t(t+1)$.  Multiplication by the contact
unit $h^2$, followed by renaming $hv$ as $v$, therefore gives
\[
v^2+t(t+1)u^4.
\]
The unit $t(t+1)$ is not a square in
$\mathbf C[t^{\pm1},(t+1)^{-1}]$, so its class appears in the notation. Consequently,
\[
\overline{\MF((\mathbf C^*)^3,W_{MM_{9-1}})}
\cong
\bigoplus_{\varepsilon\in\{\pm1\}}
\QAtwist{\Gamma_\varepsilon}{3}{t(t+1)}
\oplus
\bigoplus_{\varepsilon\in\{\pm1\}}\MorseCat{3}.
\]
\end{example}

\begin{example}[$MM_{10-1}$: two twisted mixed $A_5$ curves]
Let
\[
F=z+1+\frac{1}{xz^2},
\qquad
W_{MM_{10-1}}=xF^6+y+\frac1y-60.
\]
The reduced critical locus consists of the curves
\[
\Gamma_\varepsilon=
\left\{
y=\varepsilon,\quad x=-\frac{1}{z^2(z+1)}
\right\},
\qquad
\varepsilon\in\{\pm1\},\quad z\neq-1,
\]
and the isolated points
\[
P_\varepsilon=\left(\frac{27}{4},\varepsilon,\frac23\right).
\]
The curve has critical value $2\varepsilon-60$, and
\[
\det(\theta_i\theta_jW)(P_\varepsilon)=373248\varepsilon,
\]
so both isolated points are Morse.

With $D=y-\varepsilon$ and $t=z$, one has
\[
W_{MM_{10-1}}-(2\varepsilon-60)
=xF^6+\frac{D^2}{y}.
\]
Moreover,
\[
\det\!\left(
\frac{\partial(z,D,F)}{\partial(x,y,z)}
\right)\bigg|_{\Gamma_\varepsilon}
=t^2(t+1)^2,
\]
so $D,F$ are global normal coordinates.  Put
\[
x_0=-\frac{1}{t^2(t+1)},
\qquad
v=D(\varepsilon+D)^{-1/2},
\qquad
u=F\left(\frac{x}{x_0}\right)^{1/6}.
\]
Then the completed shifted potential is $v^2+x_0u^6$.  For
$h=\mathrm i\,t(t+1)$ one has $h^2x_0=t+1$, and contact multiplication by
$h^2$, followed by renaming $hv$ as $v$, gives
\[
v^2+(t+1)u^6.
\]
The class $[t+1]$ is nontrivial in the unit square-class group of the curve.
Thus
\[
\overline{\MF((\mathbf C^*)^3,W_{MM_{10-1}})}
\cong
\bigoplus_{\varepsilon\in\{\pm1\}}
\QAtwist{\Gamma_\varepsilon}{5}{t+1}
\oplus
\bigoplus_{\varepsilon\in\{\pm1\}}\MorseCat{3}.
\]
\end{example}

\begingroup
\scriptsize
\sloppy
\begin{longtblr}[
  caption = {Decompositions of the derived categories of matrix factorizations for the Landau--Ginzburg models of Fano threefolds.},
  label = {tab:fano_threefold_decompositions},
]{
  width = \textwidth,
  colspec = {Q[l,wd=0.12\textwidth]Q[c,wd=0.08\textwidth]Q[c,wd=0.12\textwidth]Q[l,wd=0.55\textwidth]},
  colsep = 2pt,
  rowhead = 1,
  row{1} = {font=\bfseries},
  rowsep = 2pt,
  hline{1,2,Z} = {0.4pt},
}
Threefold & Count & $\dim H^{\mathrm{even}}$ & Decomposition of $\overline{\MF((\mathbf{C}^*)^3,W)}$ \\
\SetCell[c=4]{l}\textbf{Picard-rank-one models} & & & \\
$\mathbf{P}^3$ & 4 & 4 & $\bigoplus_{n=1}^{4}\MorseCat{3}$ \\
$Q\subset\mathbf{P}^4$ & 3 & 4 & $\bigoplus_{n=1}^{3}\MorseCat{3}$ \\
$B_5$ & 4 & 4 & $\bigoplus_{n=1}^{4}\MorseCat{3}$ \\
$B_4$ & 3 & 4 & $\MBcat{\Gamma}\oplus\MorseCat{3}\oplus\MorseCat{3}$ \\
$B_3$ & 2 & 4 & $\bigoplus_{n=1}^{2}\MorseCat{3}$ \\
$B_2$ & 3 & 4 & $\left(\bigoplus_{n=1}^{2}\MorseCat{3}\right)\oplus \MBcat{S;[x]}$ \\
$B_1$ & 3 & 4 & $\mathcal C_{B_1,0}\oplus\bigoplus_{a:\,48a^2-1=0}\MorseCat{3}$ \\
$V_{22}$ & 3 & 4 & $\bigoplus_{a:\,a^3+2a^2-2=0}\MorseCat{3}$ \\
$V_{18}$ & 2 & 4 & $\bigoplus_{a:\,36a^2-12a-11=0}\MorseCat{3}$ \\
$V_{16}$ & 3 & 4 & $\MBcat{\Gamma}\oplus\bigoplus_{a^2=2}\MorseCat{3}$ \\
$V_{14}$ & 3 & 4 & $\Dper{\Gamma}{4}{t(t+1)}\oplus\MorseCat{3}\oplus\MorseCat{3}$ \\
$V_{12}$ & 3 & 4 & $\MBcat{\Gamma}\oplus\bigoplus_{a:\,32a^2-32a-1=0}\MorseCat{3}$ \\
$V_{10}$ & 3 & 4 & $\mathcal C_{V_{10},-6}\oplus\bigoplus_{a:\,4a^2-2a-1=0}\MorseCat{3}$ \\
$V_8$ & 2 & 4 & $\mathcal C_{V_8,-8}\oplus \MorseCat{3}$ \\
$V_6$ & 2 & 4 & $\mathcal C_{V_6,-12}\oplus \MorseCat{3}$ \\
$V_4$ & 2 & 4 & $\Aper{S}{3}\oplus \MorseCat{3}$ \\
$V_2$ & 2 & 4 & $\Aper{S}{5}\oplus \MorseCat{3}$ \\
\SetCell[c=4]{l}\textbf{Mori--Mukai models} & & & \\
$MM_{2-1}$ & 4 & 6 & $\mathcal C_{MM_{2-1},-60}\oplus\MBcat{\Gamma}\oplus\bigoplus_{a:\,108a^2-1=0}\MorseCat{3}$ \\
$MM_{2-2}$ & 4 & 6 & $\MBcat{S}\oplus\bigoplus_{t:\,4t^3-16t^2+16t-1=0}\MorseCat{3}$ \\
$MM_{2-3}$ & 4 & 6 & $\mathcal C_{MM_{2-3},-12}\oplus\MBcat{\Gamma}\oplus\MorseCat{3}\oplus\MorseCat{3}$ \\
$MM_{2-4}$ & 4 & 6 & $\bigoplus_{i=1}^{4}\MorseCat{3}$ \\
$MM_{2-5}$ & 3 & 6 & $\MBcat{\Gamma} \oplus \MorseCat{3}\oplus \MorseCat{3}$ \\
$MM_{2-6}$ & 3 & 6 & $\MorseCat{3}\oplus \MBcat{\Gamma_1;[t]}\oplus \MBcat{\Gamma_2}$ \\
$MM_{2-7}$ & 4 & 6 & $\MBcat{\Gamma}\oplus\bigoplus_{\alpha:\,4\alpha^3-4\alpha^2-4\alpha+3=0}\MorseCat{3}$ \\
$MM_{2-8}$ & 4 & 6 & $\MBcat{\Sigma}\oplus \bigoplus_{\rho:\,\rho^3-\rho-1=0}\MorseCat{3}$ \\
$MM_{2-9}$ & 4 & 6 & $\bigoplus_{\zeta:\,q(\zeta)=0}\MorseCat{3}$ \\
$MM_{2-10}$ & 2 & 6 & $\MorseCat{3}\oplus \MorseCat{3}$ \\
$MM_{2-11}$ & 4 & 6 & $\bigoplus_{\tau:\,q(\tau)=0}\MorseCat{3}$ \\
$MM_{2-12}$ & 2 & 6 & $\MBcat{\Gamma;[t]}\oplus \MorseCat{3}$ \\
$MM_{2-13}$ & 4 & 6 & $\MBcat{\Gamma}\oplus \bigoplus_{\xi:\,\xi^3=16}\MorseCat{3}$ \\
$MM_{2-14}$ & 4 & 6 & $\bigoplus_{\eta:\,r(\eta)=0}\MorseCat{3}$ \\
$MM_{2-15}$ & 5 & 6 & $\MBcat{\Gamma}\oplus \bigoplus_{\alpha:\,q(\alpha)=0}\MorseCat{3}$ \\
$MM_{2-16}$ & 3 & 6 & $\bigoplus_{\beta:\,\beta^3-4\beta-16=0}\MorseCat{3}$ \\
$MM_{2-17}$ & 4 & 6 & $\bigoplus_{\eta:\,r(\eta)=0}\MorseCat{3}$ \\
$MM_{2-18}$ & 4 & 6 & $\bigoplus_{\eta:\,s(\eta)=0}\MorseCat{3}$ \\
$MM_{2-19}$ & 5 & 6 & $\MBcat{\Gamma} \oplus \bigoplus_{\rho:\,h(\rho)=0}\MorseCat{3}$ \\
$MM_{2-20}$ & 6 & 6 & $\bigoplus_{\alpha:\,q(\alpha)=0}\MorseCat{3}$ \\
$MM_{2-21}$ & 4 & 6 & $\bigoplus_{i=1}^{4}\MorseCat{3}$ \\
$MM_{2-22}$ & 5 & 6 & $\bigoplus_{\lambda:\,p(\lambda)=0}\MorseCat{3}$ \\
$MM_{2-23}$ & 3 & 6 & $\bigoplus_{\alpha:\,q(\alpha)=0}\MorseCat{3}$ \\
$MM_{2-24}$ & 6 & 6 & $\bigoplus_{\zeta^3=1}\left(\MorseCat{3}\oplus \MorseCat{3}\right)$ \\
$MM_{2-25}$ & 3 & 6 & $\bigoplus_{z:\,q(z)=0}\MorseCat{3}$ \\
$MM_{2-26}$ & 6 & 6 & $\bigoplus_{u:\,p(u)=0}\MorseCat{3}$ \\
$MM_{2-27}$ & 6 & 6 & $\bigoplus_{u:\,p(u)=0}\MorseCat{3}$ \\
$MM_{2-28}$ & 5 & 6 & $\MBcat{\Gamma}\oplus\bigoplus_{\lambda:\,r(\lambda)=0}\MorseCat{3}$ \\
$MM_{2-29}$ & 5 & 6 & $\MorseCat{3}\oplus\bigoplus_{\alpha:\,q(\alpha)=0}\MorseCat{3}$ \\
$MM_{2-30}$ & 5 & 6 & $\bigoplus_{\lambda:\,p(\lambda)=0}\MorseCat{3}$ \\
$MM_{2-31}$ & 6 & 6 & $\bigoplus_{\lambda:\,p(\lambda)=0}\MorseCat{3}$ \\
$MM_{2-32}$ & 6 & 6 & $\bigoplus_{\epsilon^2=8}\MorseCat{3} \oplus \bigoplus_{\substack{\omega^2+\omega+1=0\\ \delta^2+1=0}}\MorseCat{3}$ \\
$MM_{2-33}$ & 6 & 6 & $\bigoplus_{\lambda:\,p_+(\lambda)=0}\MorseCat{3} \oplus \bigoplus_{\mu:\,p_-(\mu)=0}\MorseCat{3}$ \\
$MM_{2-34}$ & 6 & 6 & $\bigoplus_{\substack{\epsilon^2=1\\ \rho^3=1}}\MorseCat{3}$ \\
$MM_{2-35}$ & 6 & 6 & $\bigoplus_{\lambda:\,\lambda^6+\lambda^4-1=0}\MorseCat{3}$ \\
$MM_{2-36}$ & 6 & 6 & $\bigoplus_{t:\,t^6+2t^5-1=0}\MorseCat{3}$ \\
$MM_{3-1}$ & 5 & 8 & $\MorseCat{3}\oplus \MorseCat{3} \oplus \MBcat{\Gamma_3;[v]} \oplus\bigoplus_{\zeta^2+\zeta+1=0}\MBcat{\Gamma_\zeta}$ \\
$MM_{3-2}$ & 7 & 8 & $\MBcat{\Gamma}\oplus\bigoplus_{\lambda:\,\lambda^3+4\lambda^2-4=0}\MorseCat{3}\oplus\bigoplus_{\mu:\,\mu^3-2\mu^2+4=0}\MorseCat{3}$ \\
$MM_{3-3}$ & 5 & 8 & $\bigoplus_{\alpha:\,q(\alpha)=0}\MorseCat{3} \oplus \bigoplus_{\varepsilon=\pm}\MorseCat{3}$ \\
$MM_{3-4}$ & 6 & 8 & $\MorseCat{3} \oplus \bigoplus_{\beta:\,p(\beta)=0}\MorseCat{3}$ \\
$MM_{3-5}$ & 8 & 8 & $\bigoplus_{\lambda:\,p_+(\lambda)=0}\MorseCat{3} \oplus \bigoplus_{\mu:\,p_-(\mu)=0}\MorseCat{3}$ \\
$MM_{3-6}$ & 6 & 8 & $\bigoplus_{\lambda:\,r_A(\lambda)=0}\MorseCat{3} \oplus \bigoplus_{\mu:\,r_B(\mu)=0}\MorseCat{3}$ \\
$MM_{3-7}$ & 2 & 8 & $\bigoplus_{\sigma:\,\sigma^2-6\sigma+7=0}\MorseCat{3}$ \\
$MM_{3-8}$ & 8 & 8 & $\MorseCat{3} \oplus \bigoplus_{\lambda:\,q(\lambda)=0}\MorseCat{3}$ \\
$MM_{3-9}$ & 5 & 8 & $\bigoplus_{\eta:\,\eta^3=16}\MorseCat{3} \oplus \bigoplus_{\alpha:\,\alpha^2+10\alpha+1=0}\MorseCat{3}$ \\
$MM_{3-10}$ & 6 & 8 & $\bigoplus_{\lambda:\,2\lambda^3-2\lambda-1=0}\MorseCat{3} \oplus \MorseCat{3} \oplus \bigoplus_{\alpha:\,\alpha^2+9\alpha+16=0}\MorseCat{3}$ \\
$MM_{3-11}$ & 6 & 8 & $\bigoplus_{\nu:\,\Phi(\nu)=0}\MorseCat{3}$ \\
$MM_{3-12}$ & 8 & 8 & $\bigoplus_{\lambda:\,\Phi_A(\lambda)=0}\MorseCat{3} \oplus \bigoplus_{\mu:\,\Phi_B(\mu)=0}\MorseCat{3}$ \\
$MM_{3-13}$ & 1 & 8 & $\MorseCat{3}$ \\
$MM_{3-14}$ & 6 & 8 & $\bigoplus_{\lambda:\,\Phi(\lambda)=0}\MorseCat{3}$ \\
$MM_{3-15}$ & 7 & 8 & $\bigoplus_{\lambda:\,\Phi(\lambda)=0}\MorseCat{3}$ \\
$MM_{3-16}$ & 8 & 8 & $\bigoplus_{\lambda:\,\Phi(\lambda)=0}\MorseCat{3}$ \\
$MM_{3-17}$ & 4 & 8 & $\bigoplus_{\lambda:\,\Phi(\lambda)=0}\MorseCat{3}$ \\
$MM_{3-18}$ & 7 & 8 & $\bigoplus_{\lambda:\,\Phi(\lambda)=0}\MorseCat{3}$ \\
$MM_{3-19}$ & 3 & 8 & $\bigoplus_{\lambda:\,\Phi(\lambda)=0}\MorseCat{3}$ \\
$MM_{3-20}$ & 4 & 8 & $\left(\bigoplus_{\epsilon=\pm1}\MorseCat{3}\right) \oplus \bigoplus_{\mu:\,\Psi(\mu)=0}\MorseCat{3}$ \\
$MM_{3-21}$ & 8 & 8 & $\bigoplus_{\tau:\,\Phi(\tau)=0}\MorseCat{3}$ \\
$MM_{3-22}$ & 7 & 8 & $\bigoplus_{\lambda:\,\Phi(\lambda)=0}\MorseCat{3}$ \\
$MM_{3-23}$ & 8 & 8 & $\bigoplus_{\lambda:\,\Psi(\lambda)=0}\MorseCat{3}$ \\
$MM_{3-24}$ & 8 & 8 & $\bigoplus_{\tau:\,\Phi(\tau)=0}\MorseCat{3}$ \\
$MM_{3-25}$ & 8 & 8 & $\left(\bigoplus_{\alpha:\,q_+(\alpha)=0}\MorseCat{3}\right) \oplus \left(\bigoplus_{\beta:\,q_-(\beta)=0}\MorseCat{3}\right)$ \\
$MM_{3-26}$ & 8 & 8 & $\bigoplus_{\lambda:\,\Psi(\lambda)=0}\MorseCat{3}$ \\
$MM_{3-27}$ & 8 & 8 & $\bigoplus_{\epsilon\in\{1,-1\}^3}\MorseCat{3}$ \\
$MM_{3-28}$ & 8 & 8 & $\bigoplus_{\epsilon\in\{1,-1\}} \left(\bigoplus_{\lambda:\,\Phi(\lambda)=0}\MorseCat{3}\right)$ \\
$MM_{3-29}$ & 8 & 8 & $\bigoplus_{\lambda:\,\Phi(\lambda)=0}\MorseCat{3}$ \\
$MM_{3-30}$ & 8 & 8 & $\bigoplus_{\lambda:\,\Phi(\lambda)=0}\MorseCat{3}$ \\
$MM_{3-31}$ & 8 & 8 & $\bigoplus_{p\in\CritLoc(W^{\mathrm{HV}}_{3-31})}\MorseCat{3}$ \\
$MM_{4-1}$ & 2 & 10 & $\MBcat{\Gamma;[t(t+1)]}\oplus \MorseCat{3}$ \\
$MM_{4-2}$ & 6 & 10 & $\left(\bigoplus_{\alpha:\,p_+(\alpha)=0}\MorseCat{3}\right) \oplus \left(\bigoplus_{\beta:\,p_-(\beta)=0}\MorseCat{3}\right)$ \\
$MM_{4-3}$ & 2 & 10 & $\MBcat{\Gamma;[1-t]}\oplus \MorseCat{3}$ \\
$MM_{4-4}$ & 5 & 10 & $\bigoplus_{\lambda:\,\Phi(\lambda)=0}\MorseCat{3}$ \\
$MM_{4-5}$ & 5 & 10 & $\MorseCat{3}\oplus \bigoplus_{\alpha:\,q(\alpha)=0}\MorseCat{3}$ \\
$MM_{4-6}$ & 10 & 10 & $\bigoplus_{\tau:\,P(\tau)=0}\MorseCat{3}$ \\
$MM_{4-7}$ & 5 & 10 & $\left(\bigoplus_{\alpha:\,p_+(\alpha)=0}\MorseCat{3}\right) \oplus \left(\bigoplus_{\beta:\,p_-(\beta)=0}\MorseCat{3}\right)$ \\
$MM_{4-8}$ & 10 & 10 & $\left(\bigoplus_{\sigma^2=2}\MorseCat{3}\right) \oplus \left(\bigoplus_{\tau:\,R(\tau)=0}\MorseCat{3}\right)$ \\
$MM_{4-9}$ & 5 & 10 & $\bigoplus_{\lambda:\,\Phi(\lambda)=0}\MorseCat{3}$ \\
$MM_{4-10}$ & 10 & 10 & $\bigoplus_{\gamma:\,\Psi(\gamma)=0}\MorseCat{3}$ \\
$MM_{4-11}$ & 10 & 10 & $\bigoplus_{p\in\CritLoc(W^{\mathrm{HV}}_{4-10})}\MorseCat{3}$ \\
$MM_{4-12}$ & 10 & 10 & $\bigoplus_{\delta:\,\Xi(\delta)=0}\MorseCat{3}$ \\
$MM_{4-13}$ & 10 & 10 & $\bigoplus_{p\in\CritLoc(W_{MM_{4-13}})}\MorseCat{3}$ \\
$MM_{5-1}$ & 5 & 12 & $\bigoplus_{\alpha:\,p(\alpha)=0}\MorseCat{3}$ \\
$MM_{5-2}$ & 12 & 12 & $\bigoplus_{p\in\CritLoc(W^{\mathrm{HV}}_{5-2})}\MorseCat{3}$ \\
$MM_{5-3}$ & 12 & 12 & $\bigoplus_{p\in\CritLoc(W^{\mathrm{HV}}_{5-3})}\MorseCat{3}$ \\
$MM_{6-1}$ & 4 & 14 & $\bigoplus_{\varepsilon\in\{\pm1\}} \bigoplus_{\alpha:\,q(\alpha)=0}\MorseCat{3}$ \\
$MM_{7-1}$ & 4 & 16 & $\MBcat{\Gamma_1;[t]}\oplus \MBcat{\Gamma_{-1};[t]} \oplus \bigoplus_{\varepsilon\in\{\pm1\}}\MorseCat{3}$ \\
$MM_{8-1}$ & 4 & 18 & $\QAper{\Gamma_1}{2}\oplus\QAper{\Gamma_{-1}}{2}\oplus\bigoplus_{\varepsilon\in\{\pm1\}}\MorseCat{3}$ \\
$MM_{9-1}$ & 4 & 20 & $\bigoplus_{\varepsilon\in\{\pm1\}}\QAtwist{\Gamma_\varepsilon}{3}{t(t+1)}\oplus\bigoplus_{\varepsilon\in\{\pm1\}}\MorseCat{3}$ \\
$MM_{10-1}$ & 4 & 22 & $\bigoplus_{\varepsilon\in\{\pm1\}}\QAtwist{\Gamma_\varepsilon}{5}{t+1}\oplus\bigoplus_{\varepsilon\in\{\pm1\}}\MorseCat{3}$ \\
\end{longtblr}
\endgroup

\begin{remark}[Scope of the open-torus computation]
The table contains the computation on the open algebraic torus for one
specified Laurent-polynomial model for each of the $105$ deformation families
of smooth Fano threefolds. It gives the complete formal-local decomposition of
the verified matrix-factorization contributions of these
displayed functions on $(\mathbf C^*)^3$. It does not compute critical points
or vanishing cycles that may occur on the boundary of a compactification, and
no tameness theorem is asserted. Consequently, the table is not by itself a
claim about the full category of a compactified Landau--Ginzburg mirror.
\end{remark}

\subsection{Hodge and boundary comparisons}
\label{sec:hodge_boundary_comparison}
The following census is a secondary comparison of the formal-local
decompositions above with numerical and global mirror data. For a smooth
Fano threefold of Picard rank $\rho(F)=b_2(F)$,
\[
e(F)=2+2\rho(F)-2h^{1,2}(F).
\]
In particular, when $h^{1,2}(F)=0$, the Euler number is
$e(F)=2+2\rho(F)=\dim H^{\mathrm{even}}(F)$. When the critical locus
of a Laurent polynomial $W$ on the torus is finite, write
\[
\mu_{\mathrm{tor}}(W)
=\sum_{p\in\CritLoc(W)}\mu_p
\]
for the total Milnor number of its torus critical points. The Hodge numbers
used in the following count are those of the Iskovskikh and Mori--Mukai
classification
\cite{Iskovskikh1977FanoI,Iskovskikh1978FanoII,MoriMukai1981FanoB2,MoriMukai1983FanoB2,MoriMukai2003Erratum}.

For the fixed Laurent polynomials displayed in this paper, the resulting
statistics are as follows.
\begin{center}
\small
\begin{tblr}{
  width = 0.96\textwidth,
  colspec = {Q[c,wd=0.13\textwidth]Q[l,wd=0.33\textwidth]Q[l,wd=0.31\textwidth]Q[c,wd=0.08\textwidth]},
  colsep = 2pt,
  row{1} = {font=\bfseries},
  rowsep = 2pt,
  hline{1,2,Z} = {0.4pt},
}
Hodge class & Fiber and critical-locus class & Local singularity class & Cases \\
$h^{1,2}>0$ & at least one nonreduced singular fiber & positive-dimensional critical locus & $11$ \\
$h^{1,2}>0$ & reduced fibers and non-isolated critical locus & curve-type or connected gluing contribution & $15$ \\
$h^{1,2}>0$ & reduced fibers and finite critical locus & all critical points Morse & $20$ \\
$h^{1,2}=0$ & non-isolated & at least one positive-dimensional critical component & $5$ \\
$h^{1,2}=0$ & finite & all singularities isolated and Morse, with $\mu_{\mathrm{tor}}=e(F)$ & $34$ \\
$h^{1,2}=0$ & finite & all singularities isolated and Morse, with $\mu_{\mathrm{tor}}<e(F)$ & $20$ \\
\end{tblr}
\end{center}
Thus $59$ families have $h^{1,2}=0$ and $46$ have $h^{1,2}>0$. In
particular, every finite critical locus in the present list is reduced and
Morse. The five non-isolated models with $h^{1,2}=0$ all have Picard rank
greater than one. The cases and the deficient Milnor totals are listed
explicitly before Table~\ref{tab:fano_threefold_decompositions}.

For the $h^{1,2}>0$ rows, ``nonreduced'' is meant scheme-theoretically for a
fiber inside the displayed algebraic torus. In each of the eleven cases, an
exact factorization of $W-c$ contains a repeated factor; the rank-one models
$V_2,V_4,V_6,V_8,V_{10}$ are prominent examples. In the remaining
$35$ cases all singular fibers displayed on the torus are reduced. Fifteen of
them have a one-dimensional critical locus, while the other twenty have
only isolated Morse points on this chart.

This distinction has a birational motivation, but reducedness by itself is
not a rationality criterion. For example, the cubic threefold $B_3$ is
non-rational, although its displayed open-torus singular fibers are reduced
and have only isolated Morse points. Classically, if a smooth projective
threefold is rational, its principally polarized intermediate Jacobian
$J^3(F)$ is a product of Jacobians of curves; blowing up a smooth curve $C$
adds the factor $J(C)$, so successive blowups naturally give products
$\prod_iJ(C_i)$ \cite{ClemensGriffiths1972}. On the mirror side, a smooth
critical curve $C_i$ with nondegenerate transverse quadratic form contributes
the two-periodic Clifford category $\MBcat{C_i}$ associated with the curve
and its quadratic twisting data. No reconstruction of the curve or the
twisting data from the category is asserted. This relation between
Landau--Ginzburg geometry, categories, and rationality goes back to
Katzarkov's program and was developed further in
\cite{Katzarkov2010,KatzarkovPrzyjalkowski2012,IlievKatzarkovPrzyjalkowski2014,Przyjalkowski2013}.
A recent precise rationality result concerns the monodromy of reducible fibers
of a generic compactified mirror, rather than reducedness alone
\cite{LeePrzyjalkowski2025}.

The source $B_4$ model gives a useful open-torus warning. Its
positive-dimensional part of the critical locus is the punctured rational curve
$\Gamma\simeq\mathbf C^*$, although $h^{1,2}(B_4)=2$. Likewise, the critical
curves displayed for $V_{16}$ and $V_{12}$ are punctured rational curves,
although $h^{1,2}$ equals $3$ and $7$, respectively. Thus the genera and
incidence data visible in one torus chart do not, in general, recover
$J^3(F)$; boundary components and the compactified degeneration remain
essential.

It is natural to interpret the non-isolated and Milnor-deficient rows as
specializations in which critical points have collided or have left the
chosen torus chart. Varying the quantum parameters should separate these
phenomena in many families. As a quantum-side comparison, the sources used
in this census record generically semisimple small quantum cohomology for
$42$ of the $59$ families with $h^{1,2}=0$; $17$ of the $25$ non-isolated or
deficient rows belong to this group \cite{FanographyDubrovin,Ciolli2005}.
This does not by itself produce an explicit Laurent polynomial $W_q$ or prove
that all of its critical points remain in the same torus chart. Accordingly,
the decomposition table records the verified specialization, while a uniform
construction and critical-point audit of generic-parameter models $W_q$ is
kept separate from the claims of this paper.

For Picard rank one, a recent compactified calculation provides an external
comparison: Table~\ref{tab:rank_one_compactified_comparison} shows that the
isolated open-torus Morse count agrees with the compactified ordinary-double-
point count for $15$ of the $17$ families, while $V_{22}$ and $Q$ each
acquire one additional isolated point in the chosen compactification
\cite{Przyjalkowski2026RankOne}. This comparison quantifies part of the
boundary deficit without identifying the entire missing boundary category.

The period, mutation, and compactification results in
\cite{CoatesCortiGalkinKasprzyk2016,akhtar2012minkowski,
KasprzykKatzarkovPrzyjalkowskiSakovics2025,Przyjalkowski2017,
Przyjalkowski2018Toric} concern relations among Laurent representatives,
their Newton polytopes, and global compactified fibers.  By contrast, the
decompositions above attach formal-local categorical data to each fixed
displayed function on the open torus.  No invariance of these decompositions
under mutation is asserted: a birational mutation may move components of the critical locus
through the boundary, and comparison of the full compactified categories
requires additional global gluing data.

\subsection{The \texorpdfstring{$B_d$--$V_{4d+2}$}{Bd--V(4d+2)} pairs and categorical equality}
\label{sec:Bd_V4d2_categorical_comparison}

The rank-one families naturally form the five pairs
\[
(B_d,V_{4d+2}),
\qquad
d=1,2,3,4,5.
\]
The pair $(B_4,V_{18})$ is the genus-two example considered in the separate
five-chart compactification calculation. We now compare all five pairs and
make precise which categorical equalities follow from the calculations in
this paper.

\begin{definition}[Compactified supported contribution]
\label{def:general_compactified_contribution}
Let $(Y_F,w_F)$ be a chosen resolved compactification of a Landau--Ginzburg
model, let $\lambda$ be a critical value, and set
\[
\Gamma_{F,\lambda}
:=
\bigl(\CritLoc(w_F)\cap w_F^{-1}(\lambda)\bigr)_{\mathrm{red}}.
\]
Fix a finite affine cover $\mathfrak U=\{\widetilde U_i\}_{i\in\mathsf A}$
of a neighborhood of $\Gamma_{F,\lambda}$. For every nonempty
$I\subset\mathsf A$, write
$\widetilde U_I:=\bigcap_{i\in I}\widetilde U_i$, and let $\mathfrak M^{F,\lambda}_I$ 
be the dg-category of matrix factorizations of $\widehat{w_F-\lambda}$ supported on 
$\Gamma_{F,\lambda}\cap\widetilde U_I$.
The restriction functors and their composition coherences form a
semicosimplicial dg \v Cech diagram $\check{\mathfrak M}^{F,\lambda,\bullet}$. 
Relative to the chosen compactification and cover, define
\begin{equation}
\label{eq:general_compactified_contribution}
\mathfrak C^{\mathrm{comp}}_{F,\lambda}
:=
\Tot\bigl(\check{\mathfrak M}^{F,\lambda,\bullet}\bigr),
\qquad
\mathcal C^{\mathrm{comp}}_{F,\lambda}
:=
\Kar\!\left(H^0\mathfrak C^{\mathrm{comp}}_{F,\lambda}\right).
\end{equation}
If a distinguished critical value $\lambda_F$ has been fixed by the mirror correspondence, we abbreviate $\mathcal C_F^{\mathrm{comp}} := \mathcal C^{\mathrm{comp}}_{F,\lambda_F}$; by contrast, $\mathcal C^{\mathrm{comp}}_{F,0}$ always refers literally to the fiber over $0$.
\end{definition}

This definition is deliberately made at the enhanced level. Independence of the chosen cover, or identification with a global matrix-factorization category, requires dg descent. The comparisons below instead rely on explicit natural equivalences of the enhanced diagrams.

First, considering the fixed Laurent polynomials on the open torus $\mathcal T_F := \overline{\MF((\mathbf C^*)^3,W_F)}$, the decompositions established above give the following comparison:

\begin{center}
\small
\begin{tblr}{
  width=0.98\textwidth,
  colspec={Q[c,wd=0.05\textwidth]Q[c,wd=0.12\textwidth]
    Q[l,wd=0.28\textwidth]Q[l,wd=0.30\textwidth]X[l]},
  colsep=2pt,
  row{1}={font=\bfseries},
  rowsep=2pt,
  hline{1,2,Z}={0.4pt},
}
$d$ & Pair & $\mathcal T_{B_d}$ & $\mathcal T_{V_{4d+2}}$ & Open-torus status \\
$5$ & $(B_5,V_{22})$ & $\MorseCat{3}^{\oplus4}$ & $\MorseCat{3}^{\oplus3}$ & Different Morse multiplicities. \\
$4$ & $(B_4,V_{18})$ & $\MBcat{\Gamma}\oplus\MorseCat{3}^{\oplus2}$ & $\MorseCat{3}^{\oplus2}$ & $B_4$ has an additional Morse--Bott curve. \\
$3$ & $(B_3,V_{14})$ & $\MorseCat{3}^{\oplus2}$ & $\Dper{\Gamma}{4}{t(t+1)}\oplus\MorseCat{3}^{\oplus2}$ & $V_{14}$ has an additional residual $D_4$ curve. \\
$2$ & $(B_2,V_{10})$ & $\MBcat{S;[x]}\oplus\MorseCat{3}^{\oplus2}$ & $\mathcal C_{V_{10},-6}\oplus\MorseCat{3}^{\oplus2}$ & Positive-dimensional terms unidentified. \\
$1$ & $(B_1,V_6)$ & $\mathcal C_{B_1,0}\oplus\MorseCat{3}^{\oplus2}$ & $\mathcal C_{V_6,-12}\oplus\MorseCat{3}$ & Connected terms and isolated multiplicities differ. \\
\end{tblr}
\end{center}

Thus, there is no uniform equality $\mathcal T_{B_d}\simeq\mathcal T_{V_{4d+2}}$ for the displayed open-torus representatives. This failure is compatible with mutation: critical components can move through the toric boundary, necessitating a compactified analysis.

\subsubsection*{The structural compactification comparison for $B_4$ and $V_{18}$}

For $d=4$, the convenient compactifiable representatives yield pencils
$F_{B_4,\lambda}$ and $F_{V_{18},\lambda}$ with explicit critical behavior. Both models possess two isolated singularities away from the base locus with non-degenerate Hessians, contributing a global $\MorseCat{3}^{\oplus2}$ summand to each side.

At the central fiber $\lambda=0$, the resolved geometries on the two sides match. For both $B_4$ and $V_{18}$, the resolved central fiber consists of three irreducible components. Their pairwise intersections are smooth rational curves meeting at two triple points. Removing these triple points leaves three punctured curves isomorphic to $\mathbf C^*$. 

On these open strata, the components meet transversely with exact normal form $uv$, while the formal neighborhoods of the triple points on both sides are contact-equivalent to the simple normal-crossing potential $uvw$. By supported contact invariance (Corollary~\ref{cor:supported_contact}), this local geometry determines the supported categorical pieces. The punctured curves yield identical Morse--Bott contributions $\MBcat{\mathbf C^*}$, and the complete triple-point contributions $\mathcal P$ are also pairwise equivalent to $\mathcal P_{\triangle}^{\oplus 2}$, where $\mathcal P_{\triangle}$ represents matrix factorizations of $uvw$ supported at the origin.

Supported localization therefore gives identical Verdier sequences for both models:
\begin{equation}
\label{eq:B4_V18_compact_localization_summary}
\mathcal P_{\triangle}^{\oplus 2}
\longrightarrow
\mathcal C^{\mathrm{comp}}
\longrightarrow
\bigoplus_{1\leq i<j\leq3}\MBcat{\mathbf C^*}.
\end{equation}
Although the local categorical blocks match, deducing $\mathcal C^{\mathrm{comp}}_{B_4} \simeq \mathcal C^{\mathrm{comp}}_{V_{18}}$ requires more than matching the outer terms of \eqref{eq:B4_V18_compact_localization_summary}. One must establish a global natural equivalence of the underlying semicosimplicial dg \v Cech diagrams $\check{\mathfrak M}^{B,\bullet} \xrightarrow{\sim} \check{\mathfrak M}^{V,\bullet}$, properly intertwining the restriction/incidence functors and all higher coherences.

\subsubsection*{Global compactification invariants for all five pairs}

The matching of local compactified geometries extends partially to the other families. Let $\mu_{\mathrm{tor}}^{\mathrm{iso}}(F)$ be the number of open-torus isolated Morse points of the compactifiable representative defining the displayed projective pencil, $\varepsilon(F)$ the compactified ordinary double points, and $k(F)$ the irreducible components of the minimal resolved central fiber. A different source or mutation representative can have a different open-torus count. Since all ten threefolds have Picard rank one, $k(F) = h^{1,2}(F)+1$ by \cite[Theorem~22]{Przyjalkowski2013}; there $k(F)$ is the number of irreducible components of the central fiber of a Calabi--Yau compactification, which does not depend on the compactification because any two differ by flops.

Since every compactified isolated singularity in this census is an ordinary double point, the isolated subcategory is $\mathcal I_F^{\mathrm{comp}} \simeq \MorseCat{3}^{\oplus\varepsilon(F)}$. The global numerical census is as follows:

\begin{center}
\small
\begin{tblr}{
  width=0.98\textwidth,
  colspec={Q[c,wd=0.05\textwidth]Q[c,wd=0.12\textwidth]
    Q[c,wd=0.09\textwidth]Q[c,wd=0.09\textwidth]Q[c,wd=0.09\textwidth]
    Q[c,wd=0.09\textwidth]X[c]X[c]},
  colsep=2pt,
  row{1}={font=\bfseries},
  rowsep=2pt,
  hline{1,2,Z}={0.4pt},
}
$d$ & Pair & $\mu_B$ & $\varepsilon_B$ & $\mu_V$ & $\varepsilon_V$ & $k(B_d)$ & $k(V_{4d+2})$ \\
$5$ & $(B_5,V_{22})$ & $4$ & $4$ & $3$ & $4$ & $1$ & $1$ \\
$4$ & $(B_4,V_{18})$ & $2$ & $2$ & $2$ & $2$ & $3$ & $3$ \\
$3$ & $(B_3,V_{14})$ & $2$ & $2$ & $2$ & $2$ & $6$ & $6$ \\
$2$ & $(B_2,V_{10})$ & $2$ & $2$ & $2$ & $2$ & $11$ & $11$ \\
$1$ & $(B_1,V_6)$ & $2$ & $2$ & $1$ & $1$ & $22$ & $21$ \\
\end{tblr}
\end{center}

Both available global invariants---isolated nodes and irreducible central components---match exactly for $d\in \{2,3,4,5\}$. For $d=3$, both $B_3$ and $V_{14}$ have two open-torus Morse points and no additional boundary node. For $d=5$, the isolated-sector equality $\mathcal I_{B_5}^{\mathrm{comp}}\simeq \mathcal I_{V_{22}}^{\mathrm{comp}}$ is restored by a boundary node on the $V_{22}$ side. For $d=1$, however, both the multiplicities and component counts differ. 

For the extreme pairs, the total compactified categories, taken over all critical values, are as follows:
\begin{equation}
\label{eq:extreme_pair_compactified_summary}
\begin{aligned}
\overline{\MF(Y_{B_1},w_{B_1})} &\simeq \mathcal C^{\mathrm{comp}}_{B_1,0}\oplus\MorseCat{3}^{\oplus2},\\
\overline{\MF(Y_{V_6},w_{V_6})} &\simeq \mathcal C^{\mathrm{comp}}_{V_6,0}\oplus\MorseCat{3},\\
\overline{\MF(Y_{B_5},w_{B_5})} &\simeq \overline{\MF(Y_{V_{22}},w_{V_{22}})} \simeq \MorseCat{3}^{\oplus4}.
\end{aligned}
\end{equation}

\paragraph{Total versus value-zero compactified categories.}
While $\mathcal T_{B_5}^{\mathrm{comp}} \simeq \mathcal T_{V_{22}}^{\mathrm{comp}}$ unconditionally, comparing the distinguished zero-fibers requires critical value matching. In our normalization, the $B_5$ fiber over $0$ is smooth, whereas $V_{22}$ has a node at $0$. Consequently, $\mathcal C^{\mathrm{comp}}_{B_5,0}=0$ but $\mathcal C^{\mathrm{comp}}_{V_{22},0}\simeq\MorseCat{3}$. 

\begin{corollary}[Compactified comparison for the five rank-one pairs]
\label{cor:compactified_five_pair_conclusion}
For the standard compactifications of the pairs $(B_d,V_{4d+2})$, $1\leq d\leq5$, the following statements hold.
\begin{enumerate}[label=\textup{(\alph*)}]
\item The isolated compactified matrix-factorization categories satisfy $\mathcal I_{B_d}^{\mathrm{comp}} \simeq \mathcal I_{V_{4d+2}}^{\mathrm{comp}}$ if and only if $d\in\{2,3,4,5\}$.
\item The numbers of irreducible components of the resolved central fibers satisfy $k(B_d)=k(V_{4d+2})$ if and only if $d\in\{2,3,4,5\}$.
\item The total compactified categories agree unconditionally for $d=5$, giving $\MorseCat{3}^{\oplus4}$. However, their literal value-zero fibers differ without prior critical value matching.
\item For $d=4$, the punctured double-curves and complete triple-point categories are pairwise equivalent locally. A global equivalence $\mathcal C^{\mathrm{comp}}_{B_4} \simeq \mathcal C^{\mathrm{comp}}_{V_{18}}$ holds if the corresponding enhanced \v Cech diagrams are naturally equivalent.
\item For $d\in\{3,4,5\}$ and a Kuznetsov-matched pair, compatible compactified mirror identifications conditionally imply $\mathcal C^{\mathrm{comp}}_{B_d} \simeq \mathcal C^{\mathrm{comp}}_{V_{4d+2}}$ \cite{KuznetsovFano2009}.
\item For $d=2$, the Kuznetsov components are deformation equivalent but not equivalent as enhanced $\mathbf Z$-graded categories. Hence an equivalence $\mathcal C^{\mathrm{comp}}_{B_2} \simeq \mathcal C^{\mathrm{comp}}_{V_{10}}$ cannot admit a compatible $\mathbf Z$-graded lift \cite[Theorems~1.4 and~1.6]{BayerPerry2023}.
\item For $d=1$, the comparison instead passes through a nodal degree-one del Pezzo and the crepant categorical resolution of its nontrivial component \cite[Theorem~1.4]{KuznetsovShinder2025}.
\end{enumerate}
\end{corollary}

\begin{proof}
Parts~\textup{(a)}--\textup{(c)} are read off from the census table and
\eqref{eq:extreme_pair_compactified_summary}, and part~\textup{(d)} from
\eqref{eq:B4_V18_compact_localization_summary}. Parts~\textup{(e)}--\textup{(g)}
are the cited results, two-periodized.
\end{proof}

\bibliographystyle{aomalpha}
\bibliography{ref}

@article{Orlov2004,
  author = {Orlov, Dmitri},
  title = {Triangulated categories of singularities and {D}-branes in {L}andau--{G}inzburg models},
  journal = {Proceedings of the Steklov Institute of Mathematics},
  volume = {246},
  pages = {227--248},
  year = {2004},
  note = {math/0302304},
  url = {https://arxiv.org/abs/math/0302304}
}

@article{LinPomerleano2013,
  author = {Lin, Kevin H. and Pomerleano, Daniel},
  title = {Global matrix factorizations},
  journal = {Mathematical Research Letters},
  volume = {20},
  number = {1},
  pages = {91--106},
  year = {2013}
}

@book{happel1988triangulated,
  title={Triangulated Categories in the Representation Theory of Finite Dimensional Algebras},
  author={Happel, Dieter},
  volume={119},
  series={London Mathematical Society Lecture Note Series},
  year={1988},
  publisher={Cambridge University Press},
  address={Cambridge}
}

@book{atiyah1969introduction,
  title={Introduction to Commutative Algebra},
  author={Atiyah, M. F. and Macdonald, I. G.},
  year={1969},
  publisher={Addison-Wesley Publishing Company},
  address={Reading, Massachusetts}
}

@incollection{BuchweitzEisenbudHerzog1987,
  title={{C}ohen--{M}acaulay modules on quadrics},
  author={Buchweitz, Ragnar-Olaf and Eisenbud, David and Herzog, J{\"u}rgen},
  booktitle={Singularities, Representation of Algebras, and Vector Bundles},
  pages={58--116},
  year={1987},
  publisher={Springer},
  address={Berlin, Heidelberg},
  series={Lecture Notes in Mathematics},
  volume={1273},
  doi={10.1007/BFb0078838}
}

@article{kajiura2007matrix,
  title={Matrix factorizations and representations of quivers {II}: Type {ADE} case},
  author={Kajiura, Hiroshige and Saito, Kyoji and Takahashi, Atsushi},
  journal={Advances in Mathematics},
  volume={211},
  number={1},
  pages={327--362},
  year={2007},
  publisher={Elsevier}
}

@book{neeman2001triangulated,
  title={Triangulated Categories},
  author={Neeman, Amnon},
  series={Annals of Mathematics Studies},
  volume={148},
  year={2001},
  publisher={Princeton University Press},
  address={Princeton, NJ}
}

@article{verdier1996categories,
  title={Des cat{\'e}gories d{\'e}riv{\'e}es des cat{\'e}gories ab{\'e}liennes},
  author={Verdier, Jean-Louis},
  journal={Ast{\'e}risque},
  volume={239},
  pages={1--253},
  year={1996},
  publisher={Soci{\'e}t{\'e} Math{\'e}matique de France},
  note={Edited by Georges Maltsiniotis}
}

@article{isik2013equivalence,
  title={Equivalence of the derived category of a variety with a singularity category},
  author={Isik, M. Umut},
  journal={International Mathematics Research Notices},
  volume={2013},
  number={12},
  pages={2787--2808},
  year={2013},
  doi={10.1093/imrn/rns125},
  eprint={1011.1484},
  archivePrefix={arXiv}
}

@inproceedings{kontsevich1994homological,
  title={Homological algebra of mirror symmetry},
  author={Kontsevich, Maxim},
  booktitle={Proceedings of the International Congress of Mathematicians (Z{\"u}rich, 1994)},
  volume={1, 2},
  pages={120--139},
  year={1994},
  publisher={Birkh{\"a}user, Basel}
}

@article{knorrer1987cohen,
  title={{C}ohen--{M}acaulay modules on hypersurface singularities {I}},
  author={Kn{\"o}rrer, Horst},
  journal={Inventiones mathematicae},
  volume={88},
  number={1},
  pages={153--164},
  year={1987},
  publisher={Springer}
}

@article{CoatesCortiGalkinKasprzyk2016,
  author = {Coates, Tom and Corti, Alessio and Galkin, Sergey and Kasprzyk, Alexander},
  title = {Quantum periods for 3-dimensional {F}ano manifolds},
  journal = {Geometry \& Topology},
  volume = {20},
  number = {1},
  pages = {103--256},
  year = {2016},
  doi = {10.2140/gt.2016.20.103}
}

@article{KasprzykKatzarkovPrzyjalkowskiSakovics2025,
  author = {Kasprzyk, Alexander M. and Katzarkov, Ludmil and Przyjalkowski, Victor and Sakovics, Dmitrijs},
  title = {Projecting {F}anos in the mirror},
  journal = {Taiwanese Journal of Mathematics},
  volume = {29},
  number = {6},
  pages = {1411--1494},
  year = {2025},
  doi = {10.11650/tjm/250704}
}

@article{Przyjalkowski2013,
  author = {Przyjalkowski, Victor},
  title = {Weak {L}andau--{G}inzburg models for smooth {F}ano threefolds},
  journal = {Izvestiya: Mathematics},
  volume = {77},
  number = {4},
  pages = {772--794},
  year = {2013},
  doi = {10.1070/IM2013v077n04ABEH002660},
  eprint = {0902.4668},
  archivePrefix = {arXiv},
  primaryClass = {math.AG}
}

@article{ClemensGriffiths1972,
  author = {Clemens, C. Herbert and Griffiths, Phillip A.},
  title = {The intermediate {J}acobian of the cubic threefold},
  journal = {Annals of Mathematics},
  volume = {95},
  number = {2},
  pages = {281--356},
  year = {1972},
  doi = {10.2307/1970801}
}

@incollection{Katzarkov2010,
  author = {Katzarkov, Ludmil},
  title = {Generalized homological mirror symmetry and rationality questions},
  booktitle = {Cohomological and Geometric Approaches to Rationality Problems},
  series = {Progress in Mathematics},
  volume = {282},
  pages = {163--208},
  publisher = {Birkh\"auser Boston},
  year = {2010},
  doi = {10.1007/978-0-8176-4934-0_7}
}

@incollection{KatzarkovPrzyjalkowski2012,
  author = {Katzarkov, Ludmil and Przyjalkowski, Victor},
  title = {{L}andau--{G}inzburg models---old and new},
  booktitle = {Proceedings of the 18th G\"okova Geometry--Topology Conference},
  pages = {97--124},
  publisher = {International Press},
  year = {2012},
  eprint = {1405.2953},
  archivePrefix = {arXiv},
  primaryClass = {math.AG}
}

@article{IlievKatzarkovPrzyjalkowski2014,
  author = {Iliev, Atanas and Katzarkov, Ludmil and Przyjalkowski, Victor},
  title = {Double solids, categories and non-rationality},
  journal = {Proceedings of the Edinburgh Mathematical Society},
  volume = {57},
  number = {1},
  pages = {145--173},
  year = {2014},
  eprint = {1102.2130},
  archivePrefix = {arXiv},
  primaryClass = {math.AG}
}

@misc{LeePrzyjalkowski2025,
  author = {Lee, Sukjoo and Przyjalkowski, Victor},
  title = {Fibers of {L}andau--{G}inzburg models and rationality},
  year = {2025},
  eprint = {2510.21222},
  archivePrefix = {arXiv},
  primaryClass = {math.AG}
}

@article{Przyjalkowski2017,
  author = {Przyjalkowski, Victor},
  title = {{C}alabi--{Y}au compactifications of toric {L}andau--{G}inzburg models for smooth {F}ano threefolds},
  journal = {Sbornik: Mathematics},
  volume = {208},
  number = {7},
  pages = {992--1013},
  year = {2017},
  doi = {10.1070/SM8838},
  eprint = {1609.09740},
  archivePrefix = {arXiv},
  primaryClass = {math.AG}
}

@article{CheltsovPrzyjalkowski2025,
  author = {Cheltsov, Ivan and Przyjalkowski, Victor},
  title = {{K}atzarkov--{K}ontsevich--{P}antev conjecture for {F}ano threefolds},
  journal = {Proceedings of the Steklov Institute of Mathematics},
  volume = {328},
  pages = {1--156},
  year = {2025},
  doi = {10.1134/S0081543825010018},
  eprint = {1809.09218},
  archivePrefix = {arXiv},
  primaryClass = {math.AG}
}

@article{LuntsPrzyjalkowski2018,
  author = {Lunts, Valery and Przyjalkowski, Victor},
  title = {{L}andau--{G}inzburg {H}odge numbers for mirrors of del {P}ezzo surfaces},
  journal = {Advances in Mathematics},
  volume = {329},
  pages = {189--216},
  year = {2018},
  doi = {10.1016/j.aim.2018.02.024},
  eprint = {1607.08880},
  archivePrefix = {arXiv},
  primaryClass = {math.AG}
}

@article{Przyjalkowski2018Toric,
  author = {Przyjalkowski, Victor},
  title = {Toric {L}andau--{G}inzburg models},
  journal = {Russian Mathematical Surveys},
  volume = {73},
  number = {6},
  pages = {1033--1118},
  year = {2018},
  eprint = {1809.09230},
  archivePrefix = {arXiv},
  primaryClass = {math.AG}
}

@misc{Przyjalkowski2026RankOne,
  author = {Przyjalkowski, Victor},
  title = {{L}andau--{G}inzburg models for {F}ano threefolds of {P}icard rank one and exceptional collections},
  year = {2026},
  note = {arXiv:2601.16497v1, submitted 23 January 2026},
  eprint = {2601.16497},
  archivePrefix = {arXiv},
  primaryClass = {math.AG}
}

@article{WuXuZhou2025,
  author = {Wu, Yilin and Xu, Jinyi and Zhou, Guodong},
  title = {{GAGA} type results for singularity categories},
  journal = {Journal of Pure and Applied Algebra},
  volume = {229},
  number = {12},
  pages = {108146},
  year = {2025},
  doi = {10.1016/j.jpaa.2025.108146}
}

@article{Chen2010Orlov,
  author = {Chen, Xiao-Wu},
  title = {Unifying two results of {O}rlov on singularity categories},
  journal = {Abhandlungen aus dem Mathematischen Seminar der Universit\"at Hamburg},
  volume = {80},
  pages = {207--212},
  year = {2010},
  doi = {10.1007/s12188-010-0044-x},
  eprint = {1002.3467},
  archivePrefix = {arXiv},
  primaryClass = {math.AG}
}

@article{KapustinRozanskySaulina2009,
  author = {Kapustin, Anton and Rozansky, Lev and Saulina, Natalia},
  title = {Three-dimensional topological field theory and symplectic algebraic geometry. {I}},
  journal = {Nuclear Physics B},
  volume = {816},
  number = {3},
  pages = {295--355},
  year = {2009},
  eprint = {0810.5415},
  archivePrefix = {arXiv},
  primaryClass = {hep-th}
}

@article{DonovanKaroubi1970,
  author = {Donovan, Peter and Karoubi, Max},
  title = {Graded {B}rauer groups and {$K$}-theory with local coefficients},
  journal = {Publications Math\'ematiques de l'IH\'ES},
  volume = {38},
  pages = {5--25},
  year = {1970},
  doi = {10.1007/BF02684650}
}

@misc{KatzarkovKontsevichPantev2008,
  author = {Katzarkov, Ludmil and Kontsevich, Maxim and Pantev, Tony},
  title = {{H}odge theoretic aspects of mirror symmetry},
  year = {2008},
  note = {arXiv:0806.0107},
  eprint = {0806.0107},
  archivePrefix = {arXiv},
  primaryClass = {math.AG}
}

@article{MatherYau1982,
  title={Classification of isolated hypersurface singularities by their moduli algebras},
  author={Mather, John N. and Yau, Stephen S.-T.},
  journal={Inventiones mathematicae},
  volume={69},
  number={2},
  pages={243--251},
  year={1982},
  doi={10.1007/BF01399504},
  url={https://doi.org/10.1007/BF01399504},
  publisher={Springer}
}

@book{greuel2007introduction,
  title={Introduction to Singularities and Deformations},
  author={Greuel, Gert-Martin and Lossen, Christoph and Shustin, Eugenii},
  year={2007},
  publisher={Springer Science \& Business Media},
  series={Springer Monographs in Mathematics},
  address={Berlin, Heidelberg},
  isbn={978-3-540-28419-2},
  doi={10.1007/978-3-540-28419-2}
}

@article{akhtar2012minkowski,
  title={{M}inkowski polynomials and mutations},
  author={Akhtar, Mohammad and Coates, Tom and Galkin, Sergey and Kasprzyk, Alexander M.},
  journal={SIGMA. Symmetry, Integrability and Geometry: Methods and Applications},
  volume={8},
  pages={094},
  year={2012},
  publisher={SIGMA. Symmetry, Integrability and Geometry: Methods and Applications}
}

@book{Illusie1971,
  title = {Complexe Cotangent et DÃ©formations I},
  ISBN = {9783540370017},
  ISSN = {1617-9692},
  url = {http://dx.doi.org/10.1007/BFb0059052},
  DOI = {10.1007/bfb0059052},
  journal = {Lecture Notes in Mathematics},
  publisher = {Springer Berlin Heidelberg},
  author = {Illusie,  Luc},
  year = {1971}
}

@article{ORLOV2011206,
title = {Formal completions and idempotent completions of triangulated categories of singularities},
journal = {Advances in Mathematics},
volume = {226},
number = {1},
pages = {206-217},
year = {2011},
issn = {0001-8708},
doi = {https://doi.org/10.1016/j.aim.2010.06.016},
url = {https://www.sciencedirect.com/science/article/pii/S0001870810002331},
author = {Dmitri Orlov}
}

@article{Ton2011,
  title = {AlgÃ¨bres simpliciales $S^1$-Ã©quivariantes, thÃ©orie de de {R}ham et thÃ©orÃ¨mes {HKR} multiplicatifs},
  volume = {147},
  ISSN = {1570-5846},
  url = {http://dx.doi.org/10.1112/S0010437X11005501},
  DOI = {10.1112/s0010437x11005501},
  number = {6},
  journal = {Compositio Mathematica},
  publisher = {Wiley},
  author = {ToÃ«n, Bertrand and Vezzosi, Gabriele},
  year = {2011},
  month = jul,
  pages = {1979â€“2000}
}

@misc{stacks-project,
  author       = {The {Stacks project authors}},
  title        = {The {S}tacks project},
  howpublished = {\url{https://stacks.math.columbia.edu}},
  year         = {2025},
}

@incollection{thomason1990higher,
  title={Higher algebraic {K}-theory of schemes and of derived categories},
  author={Thomason, Robert W. and Trobaugh, Thomas},
  booktitle={The Grothendieck Festschrift, Volume III: A Collection of Articles Written in Honor of the 60th Birthday of Alexander Grothendieck},
  pages={247--435},
  year={1990},
  publisher={Springer}
}

@article{balmer2001idempotent,
  title={Idempotent completion of triangulated categories},
  author={Balmer, Paul and Schlichting, Marco},
  journal={Journal of Algebra},
  volume={236},
  number={2},
  pages={819--834},
  year={2001},
  publisher={Elsevier}
}

@techreport{galkin2010mutations,
  title={Mutations of potentials},
  author={Galkin, Sergey and Usnich, Alexandr},
  year={2010},
  institution={Preprint IPMU 10-0100},
  url={http://member.ipmu.jp/sergey.galkin/papers/ipmu-10-0100.pdf}
}

@book{maclane1998categories,
  title={Categories for the Working Mathematician},
  author={Mac Lane, Saunders},
  volume={5},
  year={1998},
  publisher={Springer Science \& Business Media}
}

@book{arnold1985singularities,
  title={Singularities of Differentiable Maps, Volume 1: The Classification of Critical Points, Caustics and Wave Fronts},
  author={Arnold, Vladimir I. and Gusein-Zade, Sabir M. and Varchenko, Alexander N.},
  year={1985},
  publisher={Birkh{\"a}user Boston},
  address={Boston, MA}
}

@article{Ciolli2005,
  author = {Ciolli, Gianni},
  title = {On the quantum cohomology of some {F}ano threefolds and a conjecture of {D}ubrovin},
  journal = {International Journal of Mathematics},
  volume = {16},
  number = {8},
  pages = {823--839},
  year = {2005},
  doi = {10.1142/S0129167X05003144},
  eprint = {math/0403300},
  archivePrefix = {arXiv}
}

@misc{FanographyDubrovin,
  author = {Belmans, Pieter},
  title = {{D}ubrovin's conjecture for {F}ano 3-folds and del {P}ezzo surfaces},
  howpublished = {Fanography},
  url = {https://fanography.info/dubrovin/},
  note = {Accessed 25 August 2026}
}

@article{Iskovskikh1977FanoI,
  author = {Iskovskikh, V. A.},
  title = {{F}ano 3-folds. {I}},
  journal = {Mathematics of the USSR-Izvestiya},
  volume = {11},
  number = {3},
  pages = {485--527},
  year = {1977},
  doi = {10.1070/IM1977v011n03ABEH001733}
}

@article{Iskovskikh1978FanoII,
  author = {Iskovskikh, V. A.},
  title = {{F}ano 3-folds. {II}},
  journal = {Mathematics of the USSR-Izvestiya},
  volume = {12},
  number = {3},
  pages = {469--506},
  year = {1978},
  doi = {10.1070/IM1978v012n03ABEH001994}
}

@article{MoriMukai1981FanoB2,
  author = {Mori, Shigefumi and Mukai, Shigeru},
  title = {Classification of {F}ano 3-folds with {$B_2 \geq 2$}},
  journal = {Manuscripta Mathematica},
  volume = {36},
  number = {2},
  pages = {147--162},
  year = {1981/82},
  doi = {10.1007/BF01170131}
}

@incollection{MoriMukai1983FanoB2,
  author = {Mori, Shigefumi and Mukai, Shigeru},
  title = {On {F}ano 3-folds with {$B_2 \geq 2$}},
  booktitle = {Algebraic Varieties and Analytic Varieties (Tokyo, 1981)},
  series = {Advanced Studies in Pure Mathematics},
  volume = {1},
  pages = {101--129},
  publisher = {North-Holland},
  address = {Amsterdam},
  year = {1983}
}

@article{MoriMukai2003Erratum,
  author = {Mori, Shigefumi and Mukai, Shigeru},
  title = {Erratum: ``Classification of {F}ano 3-folds with {$B_2 \geq 2$}''},
  journal = {Manuscripta Mathematica},
  volume = {110},
  number = {3},
  pages = {407},
  year = {2003},
  doi = {10.1007/s00229-002-0336-2}
}

@misc{FanoMirrorsComputations,
  author = {Ibrahim Santos El Adji, Victor},
  title = {Formal completions and decompositions of matrix factorization categories of {F}ano mirrors: computations and verification scripts},
  howpublished = {GitHub repository},
  year = {2026},
  url = {https://github.com/Victor-Ibrahim/Formal-completions-and-decompositions-of-matrix-factorization-categories-of-Fano-mirrors},
  note = {Python and Singular scripts, and the appendix \emph{Detailed Computations for the Fano Threefold Models}}
}

@article{EvansLekili2019,
  author = {Evans, Jonathan David and Lekili, Yank{\i}},
  title = {Generating the {F}ukaya categories of {H}amiltonian {$G$}-manifolds},
  journal = {Journal of the American Mathematical Society},
  volume = {32},
  number = {1},
  pages = {119--162},
  year = {2019},
  doi = {10.1090/jams/909},
  eprint = {1507.05842},
  archivePrefix = {arXiv},
  primaryClass = {math.SG}
}

@article{ChoHongLau2019,
  author = {Cho, Cheol-Hyun and Hong, Hansol and Lau, Siu-Cheong},
  title = {Localized mirror functor constructed from a {L}agrangian torus},
  journal = {Journal of Geometry and Physics},
  volume = {136},
  pages = {284--320},
  year = {2019},
  doi = {10.1016/j.geomphys.2018.11.006},
  eprint = {1406.4597},
  archivePrefix = {arXiv},
  primaryClass = {math.SG}
}

@article{RitterSmith2017,
  author = {Ritter, Alexander F. and Smith, Ivan},
  title = {The monotone wrapped {F}ukaya category and the open--closed string map},
  journal = {Selecta Mathematica. New Series},
  volume = {23},
  number = {1},
  pages = {533--642},
  year = {2017},
  doi = {10.1007/s00029-016-0255-9},
  eprint = {1201.5880},
  archivePrefix = {arXiv},
  primaryClass = {math.SG}
}

@article{Yuan2025CriticalValues,
  author = {Yuan, Hang},
  title = {Family {F}loer superpotential's critical values are eigenvalues of quantum product by {$c_1$}},
  journal = {Selecta Mathematica. New Series},
  volume = {31},
  number = {1},
  pages = {13},
  year = {2025},
  doi = {10.1007/s00029-024-01012-7},
  eprint = {2112.13537},
  archivePrefix = {arXiv},
  primaryClass = {math.SG}
}

@article{auroux2006mirror,
  title={Mirror symmetry for del {P}ezzo surfaces: vanishing cycles and coherent sheaves},
  author={Auroux, Denis and Katzarkov, Ludmil and Orlov, Dmitri},
  journal={Inventiones mathematicae},
  volume={166},
  number={3},
  pages={537--582},
  year={2006},
  publisher={Springer}
}

@unpublished{kontsevich1998lectures,
  title={Lectures at {ENS}, {P}aris, Spring 1998},
  author={Kontsevich, Maxim},
  note={Notes taken by J. Bellaiche, J.-F. Dat, I. Marin, G. Racinet, and H. Randriambololona},
  year={1998}
}

@article{auroux2007mirror,
  title={Mirror symmetry and {T}-duality in the complement of an anticanonical divisor},
  author={Auroux, Denis},
  journal={Journal of G{\"o}kova Geometry Topology},
  volume={1},
  pages={51--91},
  year={2007}
}

@article{Orlov_2011,
   title={Matrix factorizations for nonaffine {LG}--models},
   volume={353},
   ISSN={1432-1807},
   url={http://dx.doi.org/10.1007/s00208-011-0676-x},
   DOI={10.1007/s00208-011-0676-x},
   number={1},
   journal={Mathematische Annalen},
   publisher={Springer Science and Business Media LLC},
   author={Orlov, Dmitri},
   year={2012},
   month=may, pages={95--108} }

@article{eisenbud1980homological,
  author = {Eisenbud, David},
  title = {Homological algebra on a complete intersection, with an application to group representations},
  journal = {Transactions of the American Mathematical Society},
  volume = {260},
  number = {1},
  pages = {35--64},
  year = {1980},
  doi = {10.2307/1999875}
}

@article{Dyckerhoff2011,
  author = {Dyckerhoff, Tobias},
  title = {Compact generators in categories of matrix factorizations},
  journal = {Duke Mathematical Journal},
  volume = {159},
  number = {2},
  pages = {223--274},
  year = {2011},
  doi = {10.1215/00127094-1415869}
}

@article{AhmedRuas2011Relative,
  author = {Ahmed, Imran and Ruas, Maria Aparecida Soares},
  title = {Invariants of relative right and contact equivalences},
  journal = {Houston Journal of Mathematics},
  volume = {37},
  number = {3},
  pages = {773--786},
  year = {2011},
  eprint = {1003.1746},
  archivePrefix = {arXiv}
}

@article{Saito1971Quasihomogeneous,
  author = {Saito, Kyoji},
  title = {Quasihomogene isolierte Singularit{\"a}ten von Hyperfl{\"a}chen},
  journal = {Inventiones Mathematicae},
  volume = {14},
  number = {2},
  pages = {123--142},
  year = {1971},
  doi = {10.1007/BF01405360}
}

@article{PolishchukVaintrob2012,
  author = {Polishchuk, Alexander and Vaintrob, Arkady},
  title = {Chern characters and {H}irzebruch--{R}iemann--{R}och formula for matrix factorizations},
  journal = {Duke Mathematical Journal},
  volume = {161},
  number = {10},
  pages = {1863--1926},
  year = {2012},
  doi = {10.1215/00127094-1645540}
}

@book{EGAIII1,
  author = {Grothendieck, Alexander and Dieudonn{\'e}, Jean},
  title = {{\'E}l{\'e}ments de g{\'e}om{\'e}trie alg{\'e}brique III. {\'E}tude cohomologique des faisceaux coh{\'e}rents. Premi{\`e}re partie},
  series = {Publications Math{\'e}matiques de l'IH{\'E}S},
  volume = {11},
  publisher = {Institut des Hautes {\'E}tudes Scientifiques},
  year = {1961}
}

@misc{hori2000mirror,
  author = {Hori, Kentaro and Vafa, Cumrun},
  title = {Mirror symmetry},
  eprint = {hep-th/0002222},
  year = {2000},
  archivePrefix = {arXiv}
}

@article{givental1998mirror,
  author = {Givental, Alexander},
  title = {A mirror theorem for toric complete intersections},
  journal = {Topological Field Theory, Primitive Forms and Related Topics},
  pages = {141--175},
  year = {1998},
  publisher = {Birkh{\"a}user}
}

@misc{Preygel2011,
  author = {Preygel, Anatoly},
  title = {{T}hom--{S}ebastiani and duality for matrix factorizations},
  note = {arXiv:1101.5834},
  eprint = {1101.5834},
  year = {2011}
}

@article{Hirano2017,
  author = {Hirano, Yuki},
  title = {Derived {K}n{\"o}rrer periodicity and {O}rlov's theorem for gauged {L}andau--{G}inzburg models},
  journal = {Compositio Mathematica},
  volume = {153},
  number = {5},
  pages = {973--1007},
  year = {2017},
  doi = {10.1112/S0010437X16008344}
}

@article{KuznetsovQuadric2008,
  author = {Kuznetsov, Alexander},
  title = {Derived categories of quadric fibrations and intersections of quadrics},
  journal = {Advances in Mathematics},
  volume = {218},
  number = {5},
  pages = {1340--1369},
  year = {2008},
  doi = {10.1016/j.aim.2008.03.007}
}

@article{KuznetsovFano2009,
  author = {Kuznetsov, Alexander},
  title = {Derived categories of {F}ano threefolds},
  journal = {Proceedings of the Steklov Institute of Mathematics},
  volume = {264},
  pages = {110--122},
  year = {2009},
  doi = {10.1134/S0081543809010143}
}

@article{BayerPerry2023,
  author = {Bayer, Arend and Perry, Alexander},
  title = {Kuznetsov's {F}ano threefold conjecture via {K3} categories and enhanced group actions},
  journal = {Journal f{\"u}r die reine und angewandte Mathematik},
  volume = {800},
  pages = {107--153},
  year = {2023},
  doi = {10.1515/crelle-2023-0021}
}

@article{KuznetsovShinder2025,
  author = {Kuznetsov, Alexander and Shinder, Evgeny},
  title = {Derived categories of {F}ano threefolds and degenerations},
  journal = {Inventiones Mathematicae},
  volume = {239},
  number = {2},
  pages = {377--430},
  year = {2025},
  doi = {10.1007/s00222-024-01304-x}
}

@article{AuelBernardaraBolognesi2014,
  author = {Auel, Asher and Bernardara, Marcello and Bolognesi, Michele},
  title = {Fibrations in complete intersections of quadrics, {C}lifford algebras, derived categories, and rationality problems},
  journal = {Journal de Math{\'e}matiques Pures et Appliqu{\'e}es},
  volume = {102},
  number = {1},
  pages = {249--291},
  year = {2014},
  doi = {10.1016/j.matpur.2013.11.009}
}

@article{Hua2018NoncommutativeMY,
  author = {Hua, Zheng and Zhou, Gui-Song},
  title = {Noncommutative {M}ather--{Y}au theorem and its applications to {C}alabi--{Y}au algebras},
  journal = {Mathematische Annalen},
  volume = {386},
  number = {3-4},
  pages = {1605--1639},
  year = {2022},
  doi = {10.1007/s00208-022-02435-3},
  eprint = {1803.06128},
  archivePrefix = {arXiv}
}

@article{Keller2018SingularHH,
  author = {Keller, Bernhard},
  title = {Singular {H}ochschild cohomology via the singularity category},
  journal = {Comptes Rendus Math\'ematique},
  volume = {356},
  number = {11--12},
  pages = {1106--1111},
  year = {2018},
  doi = {10.1016/j.crma.2018.10.003}
}

@incollection{Keller2006dg,
  author = {Keller, Bernhard},
  title = {On differential graded categories},
  booktitle = {International Congress of Mathematicians. Vol. II},
  pages = {151--190},
  publisher = {European Mathematical Society},
  address = {Zurich},
  year = {2006},
  doi = {10.4171/022-2/8}
}

@article{HuaKeller2024,
  author = {Hua, Zheng and Keller, Bernhard},
  title = {Cluster categories and rational curves},
  journal = {Geometry \& Topology},
  volume = {28},
  number = {6},
  pages = {2569--2634},
  year = {2024},
  doi = {10.2140/gt.2024.28.2569},
  eprint = {1810.00749},
  archivePrefix = {arXiv},
  primaryClass = {math.AG}
}

@misc{Kalck2021MFClassification,
  author = {Kalck, Martin},
  title = {Classifying dg-categories of matrix factorizations},
  eprint = {2108.03292},
  year = {2021},
  archivePrefix = {arXiv}
}

@article{TelemanMorseBott,
  author = {Teleman, Constantin},
  title = {Matrix factorisation of {M}orse--{B}ott functions},
  journal = {Duke Mathematical Journal},
  volume = {169},
  number = {3},
  pages = {533--549},
  year = {2020},
  doi = {10.1215/00127094-2019-0048},
  eprint = {1611.07057},
  archivePrefix = {arXiv},
  primaryClass = {math.AG}
}

@misc{HennionHolsteinRobalo2024,
  author = {Hennion, Benjamin and Holstein, Julian and Robalo, Marco},
  title = {Gluing invariants of {D}onaldson--{T}homas type---Part {I}: the {D}arboux stack},
  year = {2024},
  eprint = {2407.08471},
  archivePrefix = {arXiv},
  primaryClass = {math.AG},
  doi = {10.48550/arXiv.2407.08471}
}

@article{EfimovPositselski2015,
  author = {Efimov, Alexander I. and Positselski, Leonid},
  title = {Coherent analogues of matrix factorizations and relative singularity categories},
  journal = {Algebra and Number Theory},
  volume = {9},
  number = {5},
  pages = {1159--1292},
  year = {2015},
  doi = {10.2140/ant.2015.9.1159},
  eprint = {1102.0261},
  archivePrefix = {arXiv}
}

\end{document}